\documentclass[11pt]{article}
\usepackage{amscd}
\usepackage{amssymb,amsfonts,amsmath,psfrag,amscd,cite}
\usepackage{amsmath}
  \usepackage{paralist}
  \usepackage{graphics}
  \usepackage{epsfig}
  \usepackage[colorlinks=true]{hyperref}
  \usepackage[all]{xy}
\usepackage{graphicx}
\newtheorem{theo}{Theorem}
\newtheorem{prop}[theo]{Proposition}

\newtheorem{defi}[theo]{Definition}
\newtheorem{lemm}[theo]{Lemma}
\newtheorem{rema}[theo]{Remark}

\makeatletter \@addtoreset{equation}{section}
\@addtoreset{theo}{section} \makeatother

\begin{document}
\date{}
\title{Hamilton-Jacobi Theorems for Regular Controlled Hamiltonian System and Its Reduced Systems}
\author{Hong Wang  \\
School of Mathematical Sciences and LPMC,\\
Nankai University, Tianjin 300071, P.R.China\\
E-mail: hongwang@nankai.edu.cn\\
May 21, 2020} \maketitle

{\bf Abstract.}
In this paper, we give precisely the geometric constraint conditions of
canonical symplectic form and regular reduced symplectic forms for the
dynamical vector fields of a regular controlled Hamiltonian
(RCH) system and its regular reduced systems.
which are called the Type I and Type II of Hamilton-Jacobi equations.
At first, we first prove two types of
Hamilton-Jacobi theorem for an RCH system on cotangent
bundle of a configuration manifold, by using the canonical symplectic form
and the dynamical vector field, which
are the development of the two types of geometric version of
Hamilton-Jacobi theorem for a Hamiltonian system given in Wang
\cite{wa17}. Moreover, we also prove
two types of Hamilton-Jacobi theorem for a controlled magnetic Hamiltonian (CMH) system,
by using the magnetic symplectic form and the magnetic Hamiltonian vector field.
Secondly, we generalize the above results for a regular
reducible RCH system with symmetry and momentum map,
and prove two types of Hamilton-Jacobi theorems
for the regular point reduced RCH system and the regular orbit reduced RCH system.
Thirdly, we prove that the
RCH-equivalence for the RCH system, and RpCH-equivalence and
RoCH-equivalence for the regular reducible RCH systems with symmetries, leave the
solutions of corresponding Hamilton-Jacobi equations invariant.
Finally, as an application of the theoretical results, we
show the Type I and Type II of Hamilton-Jacobi equations
for the $R_p$-reduced controlled rigid body-rotor system and
the $R_p$-reduced controlled heavy top-rotor
system on the generalizations of rotation group
$\textmd{SO}(3)$ and Euclidean group $\textmd{SE}(3)$, respectively.
These researches reveal the deeply internal
relationships of the geometrical structures of phase spaces, the dynamical
vector fields and controls of the RCH system.\\


\tableofcontents

\section{Introduction}

It is well-known Hamilton-Jacobi theory is an important research subject
in mathematics and analytical mechanics.
On the one hand, it provides a characterization
of the generating functions of certain time-dependent canonical
transformations, such that a given Hamiltonian system in such a form
that its solutions are extremely easy to find by reduction to the
equilibrium. On the other hand,
it is possible in many cases that Hamilton-Jacobi equation provides an
immediate way to integrate the equation of motion of a system, even
when the problem of Hamiltonian system itself has not been or cannot
be solved completely, see Abraham and Marsden \cite{abma78}, Arnold
\cite{ar89} and Marsden and Ratiu \cite{mara99}.
In addition, the Hamilton-Jacobi equation is
also fundamental in the study of the quantum-classical relationship
in quantization, and it also plays an important role
in the study of stochastic dynamical systems, see
Woodhouse \cite{wo92}, Ge and Marsden \cite{gema88},
and L\'{a}zaro-Cam\'{i} and Ortega \cite{laor09}.
For these reasons it is described as a useful tool in the study of
Hamiltonian system theory, and has been extensively developed in
past many years and become one of the most active subjects in the
study of modern applied mathematics and analytical mechanics.\\

Just as we have known that Hamilton-Jacobi theory from the
variational point of view is originally developed by Jacobi in 1866,
which state that the integral of Lagrangian of a mechanical system along the
solution of its Euler-Lagrange equation satisfies the
Hamilton-Jacobi equation. The classical description of this problem
from the generating function and the geometrical point of view is
given by Abraham and Marsden in \cite{abma78} as follows:
Let $Q$ be a smooth manifold and $TQ$
the tangent bundle, $T^* Q$ the cotangent bundle with a canonical
symplectic form $\omega$ and the projection $\pi_Q: T^* Q
\rightarrow Q. $
\begin{theo}
Assume that the triple $(T^*Q,\omega,H)$ is a Hamiltonian system
with Hamiltonian vector field $X_H$, and $W: Q\rightarrow
\mathbb{R}$ is a given generating function. Then the following two assertions
are equivalent:\\
\noindent $(\mathrm{i})$ For every curve $\sigma: \mathbb{R}
\rightarrow Q $ satisfying $\dot{\sigma}(t)= T\pi_Q
(X_H(\mathbf{d}W(\sigma(t))))$, $\forall t\in \mathbb{R}$, then
$\mathbf{d}W \cdot \sigma $ is an integral curve of the Hamiltonian
vector field $X_H$.\\
\noindent $(\mathrm{ii})$ $W$ satisfies the Hamilton-Jacobi equation
$H(q^i,\frac{\partial W}{\partial q^i})=E, $ where $E$ is a
constant.
\end{theo}

There are the following three reasons for us to develop the work.
At first, from the proof of the above theorem given in
Abraham and Marsden \cite{abma78}, we know that
the assertion $(\mathrm{i})$ with equivalent to
Hamilton-Jacobi equation by the generating function
gives a geometric constraint condition of the canonical symplectic form
on the cotangent bundle $T^*Q$
for Hamiltonian vector field of the system.
Thus, the Hamilton-Jacobi equation reveals the deeply internal relationships of
the generating function, the canonical symplectic form
and the dynamical vector field of a Hamiltonian system.\\

But, from Marsden et al.\cite{mawazh10} we know that,
the set of Hamiltonian systems with symmetries on a cotangent
bundle is not complete under the Marsden-Weinstein reduction.
Since the symplectic reduced system of a
Hamiltonian system with symmetry defined on the cotangent bundle
$T^*Q$ may not be a Hamiltonian system on a cotangent bundle,
then we cannot give the Hamilton-Jacobi theorem for the Marsden-Weinstein
reduced Hamiltonian system just like same as the above Theorem 1.1.
We have to look for a new way.\\

Secondly, it is worthy of noting that if we take that $\gamma=\mathbf{d}W$ in
the above Theorem 1.1, then $\gamma$ is a closed one-form on $Q$, and
the equation $\mathbf{d}(H \cdot \mathbf{d}W)=0$ is equivalent to
the Hamilton-Jacobi equation $H(q^i,\frac{\partial W}{\partial
q^i})=E$, where $E$ is a constant. This result is used the
formulation of a geometric version of Hamilton-Jacobi theorem for
Hamiltonian system, see Cari\~{n}ena et al \cite{cagrmamamuro06,
cagrmamamuro10}. On the other hand, Theorem 1.1 is
developed in the context of time-dependent Hamiltonian system by
Marsden and Ratiu in \cite{mara99}. The Hamilton-Jacobi equation may
be regarded as a nonlinear partial differential equation for some
generating function $S$, and the problem is become how to choose a
time-dependent canonical transformation $\Psi: T^*Q\times \mathbb{R}
\rightarrow T^*Q\times \mathbb{R}, $ which transforms the dynamical
vector field of a time-dependent Hamiltonian system to equilibrium,
such that the generating function $S$ of $\Psi$ satisfies the
time-dependent Hamilton-Jacobi equation. In particular, for the
time-independent Hamiltonian system, ones may look for a symplectic
map as the canonical transformation. This work offers an important
idea that one can use the dynamical vector field of a Hamiltonian
system to describe Hamilton-Jacobi equation. Moreover, assume that
$\gamma: Q \rightarrow T^*Q$ is a closed one-form on $Q$, and define
that $X_H^\gamma = T\pi_{Q}\cdot X_H \cdot \gamma$, where $X_{H}$ is
the dynamical vector field of Hamiltonian system $(T^*Q,\omega,H)$.
Then the fact that $X_H^\gamma$ and $X_H$ are $\gamma$-related, that
is, $T\gamma\cdot X_H^\gamma= X_H\cdot \gamma$ is equivalent that
$\mathbf{d}(H \cdot \gamma)=0, $ which is given in Cari\~{n}ena et
al \cite{cagrmamamuro06, cagrmamamuro10}.
Motivated by the above research work, Wang in \cite{wa17} prove
an important lemma, which is
a modification for the corresponding result of Abraham and Marsden
in \cite{abma78}, such that we can derive precisely
the geometric constraint conditions of
the regular reduced symplectic forms for the
dynamical vector fields of a regular reducible Hamiltonian
system, which are called the Type I and Type II of Hamilton-Jacobi equation.\\

Thirdly, we know that a regular controlled Hamiltonian (RCH) system is a
Hamiltonian system with external force and control,
which is defined in Marsden et al.\cite{mawazh10}. In general,
under the actions of external force and control,
an RCH system is not Hamiltonian,
however, it is a dynamical system closely related to a
Hamiltonian system, and it can be explored and studied by extending
the methods for external force and control in the study of Hamiltonian systems.
Thus, we can emphasize explicitly the impact of external force
and control in the study for the RCH systems.
In particular, in Marsden et al.\cite{mawazh10}, the authors
we give the regular point reduction and
the regular orbit reduction for an RCH system with
symmetry, by analyzing carefully the geometrical and topological
structures of the phase space and the reduced phase space of the
corresponding Hamiltonian system.These research work not only gave a
variety of reduction methods for the RCH systems, but
also showed a variety of relationships of controlled Hamiltonian
equivalence of these systems.
However, since an RCH system defined on the cotangent bundle
$T^*Q$, in general,  may not be a Hamiltonian system,
and it has yet no generating function,
we cannot give the Hamilton-Jacobi theorem for the RCH system
and its regular reduced systems just like same as the above Theorem 1.1.
Thus, it is a natural problem how to derive precisely the geometric constraint conditions of
canonical symplectic form and regular reduced symplectic forms for the
dynamical vector fields of an RCH system and its regular reduced systems,
and how to describe explicitly the relationship between the RCH-equivalence
and the solutions of corresponding Hamilton-Jacobi equations.
These research are our goal in this paper.\\

A brief of outline of this paper is as follows. In the second
section, we first review some relevant definitions and basic facts
about the RCH systems and RCH-equivalence, then prove a key lemma,
which is an important tool for the proofs of two types of Hamilton-Jacobi
theorems of the RCH system and its regular reduced systems. Then we
derive precisely the geometric constraint conditions of
canonical symplectic form for the dynamical vector fields
of an RCH system on the cotangent bundle of a
configuration manifold, that is,
the two types of Hamilton-Jacobi equation for an RCH system.
Moreover, in the third section,
in order to describe the impact of different geometric structures
for the Hamilton-Jacobi theorem, we prove
two types of Hamilton-Jacobi theorem for a controlled magnetic Hamiltonian (CMH) system,
by using the magnetic symplectic form and the magnetic Hamiltonian vector field.
From the fourth section we begin to discuss the regular
reducible RCH system with symmetry and momentum map, by combining with the
Hamilton-Jacobi theory and regular symplectic reduction theory for
RCH system, and the two types of Hamilton-Jacobi
equations for the $R_p$-reduced and $R_o$-reduced RCH systems
are obtained respectively in the fourth
section and the fifth section, by using the reduced symplectic forms and
the regular reducible dynamical vector fields.
These results are the development of the two types of
Hamilton-Jacobi theorems for a Hamiltonian system
and its reduced systems given in Wang \cite{wa17}.
In particular, we prove that the
RCH-equivalence for the RCH system, and RpCH-equivalence and
RoCH-equivalence for the regular reducible RCH systems with symmetries, leave the
solutions of corresponding Hamilton-Jacobi equations invariant.
Finally, as the applications of the theoretical results,
in sixth section, we consider the regular point
reducible RCH system on the generalization of a Lie group,
and derive precisely the regular point reduction and the
two types of Hamilton-Jacobi equation for the reduced system.
In particular, we show the Type I and Type II of Hamilton-Jacobi equations
for the $R_p$-reduced controlled rigid body-rotor system and
the $R_p$-reduced controlled heavy top-rotor
system on the generalization of rotation group
$\textmd{SO}(3)$ and on the generalization of Euclidean group
$\textmd{SE}(3)$ by calculation in detail, respectively.
These research work reveal the deeply internal
relationships of the geometrical structures of phase spaces, the dynamical
vector fields and controls of the RCH system,
and develop the theory and application of the regular symplectic reduction
and Hamilton-Jacobi theory for the RCH systems
with symmetries, and make us have much deeper understanding and
recognition for the structure of Hamiltonian systems and RCH
systems.

\section{Hamilton-Jacobi Theorem for an RCH System }

In this paper, our goal is to study Hamilton-Jacobi theory for an RCH
system with symplectic structure and symmetry, and prove the two types of
Hamilton-Jacobi theorems for an RCH system and its regular reduced RCH
systems, and describe the relationship between the RCH-equivalence
for the RCH systems and the solutions of corresponding Hamilton-Jacobi
equations. In order to do these, in this section, we first review
some relevant definitions and basic facts about RCH systems and
RCH-equivalence, which will be used in subsequent sections. Then
we give precisely the geometric constraint conditions of canonical symplectic form for the
dynamical vector field of the RCH system, that is,
the two types of Hamilton-Jacobi equation, and state that the solution of Hamilton-Jacobi equation for
the RCH system leaves invariant under the conditions of RCH-equivalence.
We shall follow some of the notations and conventions introduced in Abraham
and Marsden \cite{abma78}, Marsden and Ratiu \cite{mara99}, Marsden \cite{ma92}
Libermann and Marle \cite{lima87}, Ortega and Ratiu \cite{orra04},
Marsden et al \cite{mawazh10} and Wang \cite{wa17}. In this paper,
we assume that all manifolds are real, smooth and finite dimensional
and all actions are smooth left actions.\\

In order to describe uniformly RCH systems defined on a cotangent
bundle and on the regular reduced spaces, in this subsection we
first define an RCH system on a symplectic fiber bundle,
see Marsden et al.\cite{mawazh10}. Then we can
obtain the RCH system on the cotangent bundle of a configuration
manifold as a special case, and discuss RCH-equivalence. In
consequence, we can regard the associated Hamiltonian system on the
cotangent bundle as a spacial case of the RCH system without
external force and control, such that we can study the RCH systems
with symmetries by combining with regular symplectic reduction theory
of Hamiltonian systems. For convenience, we assume that all
controls appearing in this paper are the admissible controls.\\

Let $(E,M,N,\pi,G)$ be a fiber bundle and $(E, \omega_E)$ be a
symplectic fiber bundle. If for any function $H: E \rightarrow
\mathbb{R}$, we have a Hamiltonian vector field $X_H$ by
$i_{X_H}\omega_E=\mathbf{d}H$, then $(E, \omega_E, H )$ is a
Hamiltonian system. Moreover, if considering the external force and
control, we can define a kind of regular controlled Hamiltonian
(RCH) system on the symplectic fiber bundle $E$ as
follows.
\begin{defi}
(RCH System) An RCH system on $E$ is a 5-tuple
$(E, \omega_E, H, F, W)$, where $(E, \omega_E, H )$ is a
Hamiltonian system, and the function $H: E \rightarrow \mathbb{R}$
is called the Hamiltonian, a fiber-preserving map $F: E\rightarrow
E$ is called the (external) force map, and a fiber submanifold $W$
of $E$ is called the control subset.
\end{defi}
Sometimes, $W$ also denotes the set of fiber-preserving maps from
$E$ to $W$. When a feedback control law $u: E\rightarrow W$ is
chosen, the 5-tuple $(E, \omega_E, H, F, u)$ denotes a closed-loop
dynamic system. In particular, when $Q$ is a smooth manifold, and
$T^\ast Q$ its cotangent bundle with a symplectic form $\omega$ (not
necessarily canonical symplectic form), then $(T^\ast Q, \omega )$
is a symplectic vector bundle. If we take that $E= T^* Q$, from
above definition we can obtain an RCH system on the cotangent bundle
$T^\ast Q$, that is, 5-tuple $(T^\ast Q, \omega, H, F, W)$. Where
the fiber-preserving map $F: T^*Q\rightarrow T^*Q$ is the (external)
force map, that is the reason that the fiber-preserving map $F:
E\rightarrow E$ is called an (external) force map in above
definition.\\

In order to describe the dynamics of the RCH system
$(E,\omega_E,H,F,W)$ with a control law $u$, we need to give a good
expression of the dynamical vector field of the RCH system. At first, we
introduce a notations of vertical lift maps of a vector along a
fiber. For a smooth manifold $E$, its tangent bundle $TE$ is a
vector bundle, and for the fiber bundle $\pi: E \rightarrow M$, we
consider the tangent mapping $T\pi: TE \rightarrow TM$ and its
kernel $ker (T\pi)=\{\rho\in TE| T\pi(\rho)=0\}$, which is a vector
subbundle of $TE$. Denote by $VE:= ker(T\pi)$, which is called a
vertical bundle of $E$. Assume that there is a metric on $E$, and we
take a Levi-Civita connection $\mathcal{A}$ on $TE$, and denote by
$HE:= ker(\mathcal{A})$, which is called a horizontal bundle of $E$,
such that $TE= HE \oplus VE. $ For any $x\in M, \; a_x, b_x \in E_x,
$ any tangent vector $\rho(b_x)\in T_{b_x}E$ can be split into
horizontal and vertical parts, that is, $\rho(b_x)=
\rho^h(b_x)\oplus \rho^v(b_x)$, where $\rho^h(b_x)\in H_{b_x}E$ and
$\rho^v(b_x)\in V_{b_x}E$. Let $\gamma$ be a geodesic in $E_x$
connecting $a_x$ and $b_x$, and denote by $\rho^v_\gamma(a_x)$ a
tangent vector at $a_x$, which is a parallel displacement of the
vertical vector $\rho^v(b_x)$ along the geodesic $\gamma$ from $b_x$
to $a_x$. Since the angle between two vectors is invariant under a
parallel displacement along a geodesic, then
$T\pi(\rho^v_\gamma(a_x))=0, $ and hence $\rho^v_\gamma(a_x) \in
V_{a_x}E. $ Now, for $a_x, b_x \in E_x $ and tangent vector
$\rho(b_x)\in T_{b_x}E$, we can define the vertical lift map of a
vector along a fiber given by
$$\textnormal{vlift}: TE_x \times E_x \rightarrow TE_x; \;\;
\textnormal{vlift}(\rho(b_x),a_x) = \rho^v_\gamma(a_x). $$
It is easy to check from the basic fact in differential geometry
that this map does not depend on the choice of $\gamma$. If $F: E
\rightarrow E$ is a fiber-preserving map, for any $x\in M$, we have
that $F_x: E_x \rightarrow E_x$ and $TF_x: TE_x \rightarrow TE_x$,
then for any $a_x \in E_x$ and $\rho\in TE_x$, the vertical lift of
$\rho$ under the action of $F$ along a fiber is defined by
$$(\textnormal{vlift}(F_x)\rho)(a_x)
=\textnormal{vlift}((TF_x\rho)(F_x(a_x)), a_x)
= (TF_x\rho)^v_\gamma(a_x), $$
where $\gamma$ is a geodesic in $E_x$ connecting $F_x(a_x)$ and
$a_x$.\\

In particular, when $\pi: E \rightarrow M$ is a vector bundle, for
any $x\in M$, the fiber $E_x=\pi^{-1}(x)$ is a vector space. In this
case, we can choose the geodesic $\gamma$ to be a straight line, and
the vertical vector is invariant under a parallel displacement along
a straight line, that is, $\rho^v_\gamma(a_x)= \rho^v(b_x).$
Moreover, when $E= T^*Q, \; M=Q $, by using the local trivialization
of $TT^*Q$, we have that $TT^*Q\cong TQ \times T^*Q$. Because of
$\pi: T^*Q \rightarrow Q$, and $T\pi: TT^*Q \rightarrow TQ$, then in
this case, for any $\alpha_x, \; \beta_x \in T^*_x Q, \; x\in Q, $
we know that $(0, \beta_x) \in V_{\beta_x}T^*_x Q, $ and hence we
can get that
$$ \textnormal{vlift}((0, \beta_x)(\beta_x), \alpha_x) = (0, \beta_x)(\alpha_x)
= \left.\frac{\mathrm{d}}{\mathrm{d}s}\right|_{s=0}(\alpha_x+s\beta_x),
$$ which is consistent with the definition of vertical lift map
along fiber in Marsden and Ratiu \cite{mara99}.\\

For a given RCH System $(T^\ast Q, \omega, H, F, W)$, the dynamical
vector field of the associated Hamiltonian system $(T^\ast Q,
\omega, H) $ is  $X_H$, which satisfies the equation
$i_{X_H}\omega=\mathbf{d}H$. If considering the
external force $F: T^*Q \rightarrow T^*Q, $ by using the above
notation of vertical lift map of a vector along a fiber, the change
of $X_H$ under the action of $F$ is that
$$\textnormal{vlift}(F)X_H(\alpha_x)
= \textnormal{vlift}((TFX_H)(F(\alpha_x)), \alpha_x)
= (TFX_H)^v_\gamma(\alpha_x),$$
where $\alpha_x \in T^*_x Q, \; x\in Q $ and $\gamma$ is a straight
line in $T^*_x Q$ connecting $F_x(\alpha_x)$ and $\alpha_x$. In the
same way, when a feedback control law $u: T^\ast Q \rightarrow W$ is
chosen, the change of $X_H$ under the action of $u$ is that
$$\textnormal{vlift}(u)X_H(\alpha_x)
= \textnormal{vlift}((TuX_H)(u(\alpha_x)), \alpha_x)
= (TuX_H)^v_\gamma(\alpha_x).$$
In consequence, we can give an expression of the dynamical vector
field of the RCH system as follows.
\begin{theo}
The dynamical vector field of an RCH system $(T^\ast Q,\omega,H,F,W)$
with a control law $u$ is the synthetic of Hamiltonian vector field
$X_H$ and its changes under the actions of the external force $F$
and control $u$, that is,
\begin{equation}
X_{(T^\ast Q,\omega,H,F,u)}(\alpha_x)
= X_H(\alpha_x)+ \textnormal{vlift}(F)X_H(\alpha_x)
+ \textnormal{vlift}(u)X_H(\alpha_x), \;\; \label{2.1}
\end{equation}
 for any $\alpha_x \in T^*_x
Q, \; x\in Q $. For convenience, it is simply written as
\begin{equation}X_{(T^\ast Q,\omega,H,F,u)}
=X_H +\textnormal{vlift}(F) +\textnormal{vlift}(u). \;\; \label{2.2}
\end{equation}
\end{theo}
We also denote that $\textnormal{vlift}(W)= \bigcup\{\textnormal{vlift}(u)X_H |
\; u\in W\}$. It is worthy of noting that in order to deduce and calculate
easily, we always use the simple expression of dynamical vector
field $X_{(T^\ast Q,\omega,H,F,u)}$. \\

In the following we shall derive precisely the geometric constraint conditions of
canonical symplectic form for the dynamical vector field
of the RCH system, that is, Type I and Type II of Hamilton-Jacobi equation for
the RCH system, and state that the solution of the corresponding
Hamilton-Jacobi equation leaves invariant under the conditions of
RCH-equivalence. In order to do this, in the following we first give
an important notion and prove a key lemma, which is an important
tool for the proofs of two types of
Hamilton-Jacobi theorem for the RCH system.\\

Denote by $\Omega^i(Q)$ the set of all i-forms on $Q$, $i=1,2.$
For any $\gamma \in \Omega^1(Q),\; q\in Q, $ then $\gamma(q)\in T_q^*Q, $
and we can define a map $\gamma: Q \rightarrow T^*Q, \; q \rightarrow (q, \gamma(q)).$
Hence we say often that the map $\gamma: Q
\rightarrow T^*Q$ is an one-form on $Q$. If the one-form $\gamma$ is closed,
then $\mathbf{d}\gamma(x,y)=0, \; \forall\;
x, y \in TQ$. In the following we give a weaker notion.
\begin{defi}
The one-form $\gamma$ is called to be closed with respect to $T\pi_{Q}:
TT^* Q \rightarrow TQ, $ if for any $v, w \in TT^* Q, $ we have
$\mathbf{d}\gamma(T\pi_{Q}(v),T\pi_{Q}(w))=0. $
\end{defi}
From the above definition we know that, if $\gamma$ is a closed one-form,
then it must be closed with respect to $T\pi_{Q}: TT^* Q \rightarrow
TQ. $ Conversely, if $\gamma$ is closed with respect to
$T\pi_{Q}: TT^* Q \rightarrow TQ, $ then it may not be closed. We can
prove a general result as follows.
\begin{prop}
Assume that $\gamma: Q \rightarrow T^*Q$ is an one-form on $Q$ and
it is not closed. we define the set $N$, which is a subset of $TQ$,
such that the one-form $\gamma$ on $N$ satisfies the condition that
for any $x,y \in N, \; \mathbf{d}\gamma(x,y)\neq 0. $ Denote by
$Ker(T\pi_Q)= \{u \in TT^*Q| \; T\pi_Q(u)=0 \}, $ and $T\gamma: TQ
\rightarrow TT^* Q .$ If $T\gamma(N)\subset Ker(T\pi_Q), $ then
$\gamma$ is closed with respect to $T\pi_{Q}: TT^* Q \rightarrow TQ.
$\end{prop}

Now, we prove the following Lemma 2.5. It is worthy of noting that
this lemma is obtained by a careful modification for the
corresponding result of Abraham and Marsden in \cite{abma78}, also
see Wang \cite{wa17}. This lemma is very important for our research,
and we also give its proof here.
\begin{lemm}
Assume that $\gamma: Q \rightarrow T^*Q$ is an one-form on $Q$, and
$\lambda=\gamma \cdot \pi_{Q}: T^* Q \rightarrow T^* Q .$ Then
we have that the following two assertions hold.\\
\noindent $(\mathrm{i})$ For any $x, y \in TQ, \;
\gamma^*\omega(x,y)= -\mathbf{d}\gamma (x,y),$ and for any $v, w \in
TT^* Q, \; \lambda^*\omega(v,w)=\\ -\mathbf{d}\gamma(T\pi_{Q}(v), \;
T\pi_{Q}(w)),$
since $\omega$ is the canonical symplectic form on $T^*Q$; \\
\noindent $(\mathrm{ii})$ For any $v, w \in TT^* Q, \;
\omega(T\lambda \cdot v,w)= \omega(v, w-T\lambda \cdot
w)-\mathbf{d}\gamma(T\pi_{Q}(v), \; T\pi_{Q}(w)). $
\end{lemm}
\noindent{\bf Proof:}
We first prove the assertion $(\mathrm{i})$.
Since $\omega$ is the canonical symplectic form on $T^*Q$, we know
that there is an unique canonical one-form $\theta$, such that
$\omega= -\mathbf{d} \theta. $ From the Proposition 3.2.11 in
Abraham and Marsden \cite{abma78}, we have that for the one-form
$\gamma: Q \rightarrow T^*Q, \; \gamma^* \theta= \gamma. $ Then we
can obtain that
\begin{align*}
\gamma^*\omega(x,y) = \gamma^* (-\mathbf{d} \theta) (x, y) =
-\mathbf{d}(\gamma^* \theta)(x, y)= -\mathbf{d}\gamma (x, y).
\end{align*}
Note that $\lambda=\gamma \cdot \pi_{Q}: T^* Q \rightarrow T^* Q, $
and $\lambda^*= \pi_{Q}^* \cdot \gamma^*: T^*T^* Q \rightarrow
T^*T^* Q, $ then we have that
\begin{align*}
\lambda^*\omega(v,w) &= \lambda^* (-\mathbf{d} \theta) (v, w)
=-\mathbf{d}(\lambda^* \theta)(v, w)= -\mathbf{d}(\pi_{Q}^* \cdot
\gamma^* \theta)(v, w)\\ &= -\mathbf{d}(\pi_{Q}^* \cdot\gamma )(v,
w)= -\mathbf{d}\gamma(T\pi_{Q}(v), \; T\pi_{Q}(w)).
\end{align*}
It follows that the assertion $(\mathrm{i})$ holds.\\

Next, we prove the assertion $(\mathrm{ii})$. For any $v, w \in TT^*
Q,$ note that $v- T(\gamma \cdot \pi_Q)\cdot v$ is vertical, because
$$
T\pi_Q(v- T(\gamma \cdot \pi_Q)\cdot v)=T\pi_Q(v)-T(\pi_Q\cdot
\gamma\cdot \pi_Q)\cdot v= T\pi_Q(v)-T\pi_Q(v)=0,
$$
where we have used the relation $\pi_Q\cdot \gamma\cdot \pi_Q= \pi_Q. $
Thus, $\omega(v- T(\gamma \cdot \pi_Q)\cdot v,w- T(\gamma \cdot
\pi_Q)\cdot w)= 0, $ and hence,
$$\omega(T(\gamma \cdot \pi_Q)\cdot v, \; w)=
\omega(v, \; w-T(\gamma \cdot \pi_Q)\cdot w)+ \omega(T(\gamma \cdot
\pi_Q)\cdot v, \; T(\gamma \cdot \pi_Q)\cdot w). $$ However, the
second term on the right-hand side is given by
$$
\omega(T(\gamma \cdot \pi_Q)\cdot v, \; T(\gamma \cdot \pi_Q)\cdot
w)= \gamma^*\omega(T\pi_Q(v), \; T\pi_Q(w))=
-\mathbf{d}\gamma(T\pi_{Q}(v), \; T\pi_{Q}(w)),
$$
where we have used the assertion $(\mathrm{i})$. It follows that
\begin{align*}
\omega(T\lambda \cdot v,w) &=\omega(T(\gamma \cdot \pi_Q)\cdot v, \;
w)\\ &= \omega(v, \; w-T(\gamma \cdot \pi_Q)\cdot w)-\mathbf{d}\gamma(T\pi_{Q}(v), \; T\pi_{Q}(w))
\\ &= \omega(v,
w-T\lambda \cdot w)-\mathbf{d}\gamma(T\pi_{Q}(v), \; T\pi_{Q}(w)).
\end{align*}
Thus, the assertion $(\mathrm{ii})$ holds.
\hskip 0.3cm $\blacksquare$\\

From the expression of the dynamical vector
field of an RCH system, we know that under the actions of the external force $F$
and control $u$, in general, the dynamical vector
field is not Hamiltonian, and the RCH system is not
yet a Hamiltonian system, and hence, we can not describe the Hamilton-Jacobi equation for an
RCH system from the viewpoint of generating
function just like same as Theorem 1.1
given by Abraham and Marsden in \cite{abma78}.
But, for a given RCH system $(T^*Q,\omega,H,F,W)$ in which $\omega$ is the
canonical symplectic form on $T^*Q$, by using
the above Lemma 2.5 and the dynamical vector field $X_{(T^\ast Q,\omega,H,F,u)}$,
we can derive precisely the following two types of geometric constraint condition of
canonical symplectic form for the dynamical vector field of the RCH system,
that is, the two types of geometric
Hamilton-Jacobi equation for the RCH system.
At first, by using the fact that the one-form $\gamma: Q
\rightarrow T^*Q $ is closed with respect to
$T\pi_Q: TT^* Q \rightarrow TQ, $ we can prove the Type I of geometric
Hamilton-Jacobi theorem for the RCH system. For convenience,
the maps involved in the following theorem and its proof are shown
in Diagram-1.
\begin{center}
\hskip 0cm \xymatrix{ & T^* Q \ar[d]_{X_H}\ar[r]^{\pi_Q}
& Q \ar[d]_{\tilde{X}^\gamma} \ar[r]^{\gamma} & T^*Q \ar[d]^{\tilde{X}} \\
  & T(T^*Q) & TQ \ar[l]_{T\gamma} & T(T^* Q)\ar[l]_{T\pi_Q}}
\end{center}
$$\mbox{Diagram-1}$$
\begin{theo} (Type I of Hamilton-Jacobi Theorem for an RCH System)
For the RCH system $(T^*Q,\omega,H,F,W)$ with the
canonical symplectic form $\omega$ on $T^*Q$, assume that $\gamma: Q
\rightarrow T^*Q$ is an one-form on $Q$, and $\tilde{X}^\gamma = T\pi_{Q}\cdot \tilde{X} \cdot \gamma$,
where $\tilde{X}=X_{(T^\ast Q,\omega,H,F,u)}$ is the dynamical vector field of the RCH system
$(T^*Q,\omega,H,F,W)$ with a control law $u$.
If the one-form $\gamma: Q \rightarrow T^*Q $ is closed with respect to
$T\pi_Q: TT^* Q \rightarrow TQ, $ then $\gamma$ is a solution of the equation
$T\gamma\cdot \tilde{X}^\gamma= X_H\cdot \gamma ,$ where $X_H$ is the Hamiltonian vector field
of the corresponding Hamiltonian system $(T^*Q,\omega,H),$ and the equation is called the Type I of
Hamilton-Jacobi equation for the RCH system $(T^*Q,\omega,H,F,W)$ with a control law $u$.
\end{theo}
\noindent{\bf Proof: } Since $\tilde{X}=X_{(T^\ast Q,\omega,H,F,u)}=X_H
+\textnormal{vlift}(F)+\textnormal{vlift}(u), $ and
$T\pi_{Q}\cdot \textnormal{vlift}(F)=T\pi_{Q}\cdot \textnormal{vlift}(u)=0, $
then we have that $T\pi_{Q}\cdot \tilde{X}\cdot \gamma=T\pi_{Q}\cdot X_H\cdot \gamma. $
If we take that $v= X_H\cdot \gamma \in TT^* Q, $ and for
any $w \in TT^* Q, \; T\pi_{Q}(w)\neq 0, $ from Lemma 2.5(ii) we have that
\begin{align*}
& \omega(T\gamma \cdot \tilde{X}^\gamma, \; w)= \omega(T\gamma \cdot
T\pi_Q\cdot \tilde{X}\cdot \gamma, \; w)= \omega(T\gamma \cdot
T\pi_Q\cdot X_H\cdot \gamma, \; w)\\ &=\omega(T(\gamma \cdot
\pi_Q)\cdot X_H\cdot \gamma, \; w)= \omega(X_H\cdot \gamma, \;
w-T(\gamma \cdot \pi_Q)\cdot
w)-\mathbf{d}\gamma(T\pi_{Q}(X_H\cdot \gamma), \; T\pi_{Q}(w))\\
& =\omega(X_H\cdot \gamma, \; w) - \omega(X_H\cdot \gamma, \;
T\lambda \cdot w)-\mathbf{d}\gamma(T\pi_{Q}(X_H\cdot \gamma), \; T\pi_{Q}(w)).
\end{align*}
Because the one-form $\gamma: Q \rightarrow T^*Q $ is closed with respect to
$T\pi_Q: TT^* Q \rightarrow TQ, $ then we have that
$$
\mathbf{d}\gamma(T\pi_{Q}(X_H\cdot \gamma), \; T\pi_{Q}(w))=0,
$$
and hence
\begin{equation}
\omega(T\gamma \cdot \tilde{X}^\gamma, \; w)- \omega(X_H\cdot \gamma, \; w)
= -\omega(X_H\cdot \gamma, \; T\lambda \cdot w). \;\; \label{2.3}
\end{equation}
If $\gamma$ satisfies the equation $T\gamma\cdot \tilde{X}^\gamma= X_H\cdot \gamma ,$
from Lemma 2.5(i) we can obtain that
\begin{align*}
-\omega(X_H\cdot \gamma, \; T\lambda \cdot w) &
= -\omega(T\gamma \cdot \tilde{X}^\gamma, \; T\lambda \cdot w)\\
& =-\omega(T\gamma \cdot T\pi_{Q} \cdot \tilde{X}\cdot\gamma, \; T\lambda \cdot w)
=-\omega(T\lambda \cdot \tilde{X}\cdot\gamma, \; T\lambda \cdot w)\\
& = -\lambda^*\omega( \tilde{X}\cdot\gamma, \; w)=
\mathbf{d}\gamma(T\pi_{Q}( \tilde{X}\cdot\gamma ), \; T\pi_{Q}(w))=0,
\end{align*}
since the one-form $\gamma: Q \rightarrow T^*Q $ is closed with respect to
$T\pi_Q: TT^* Q \rightarrow TQ. $ But, because the symplectic form $\omega$ is non-degenerate,
the left side of (2.3) equals zero, only when
$\gamma$ satisfies the equation $T\gamma\cdot \tilde{X}^\gamma= X_H\cdot \gamma .$ Thus,
if the one-form $\gamma: Q \rightarrow T^*Q $ is closed with respect to
$T\pi_Q: TT^* Q \rightarrow TQ, $ then $\gamma$ must be a solution of
the Type I of Hamilton-Jacobi equation
$T\gamma\cdot \tilde{X}^\gamma= X_H\cdot \gamma .$
\hskip 0.3cm $\blacksquare$\\

Next, for any symplectic map $\varepsilon: T^* Q \rightarrow T^* Q $,
we can prove the following Type II of geometric
Hamilton-Jacobi theorem for the RCH system $(T^*Q,\omega,H,F,W)$. For convenience,
the maps involved in the following theorem and its proof are shown
in Diagram-2.
\begin{center}
\hskip 0cm \xymatrix{ & T^* Q \ar[r]^{\varepsilon}
& T^*Q \ar[d]_{X_{H\cdot \varepsilon}}
\ar[dr]^{\tilde{X}^\varepsilon} \ar[r]^{\pi_Q}
& Q \ar[r]^{\gamma} & T^*Q \ar[d]^{\tilde{X}} \\
&  & T(T^*Q) & TQ \ar[l]_{T\gamma} & T(T^* Q)\ar[l]_{T\pi_Q}}
\end{center}
$$\mbox{Diagram-2}$$
\begin{theo}
 (Type II of Hamilton-Jacobi Theorem for an RCH System)
For the RCH system $(T^*Q,\omega,H,F,W)$ with the
canonical symplectic form $\omega$ on $T^*Q$, assume that $\gamma: Q
\rightarrow T^*Q$ is an one-form on $Q$, and
$\lambda=\gamma\cdot\pi_{Q}: T^* Q \rightarrow T^* Q $, and for any
symplectic map $\varepsilon: T^* Q \rightarrow T^* Q $, denote by
$\tilde{X}^\varepsilon = T\pi_{Q}\cdot \tilde{X} \cdot \varepsilon$,
where $\tilde{X}=X_{(T^\ast Q,\omega,H,F,u)}$
is the dynamical vector field of the RCH system
$(T^*Q,\omega,H,F,W)$ with a control law $u$.
Then $\varepsilon$ is a solution of the equation
$T\varepsilon\cdot X_{H\cdot\varepsilon}= T\lambda \cdot \tilde{X} \cdot \varepsilon,$
if and only if it is a solution of the equation $T\gamma\cdot \tilde{X}^\varepsilon= X_H\cdot
\varepsilon, $ where $X_H$ and $ X_{H\cdot\varepsilon} \in
TT^*Q $ are the Hamiltonian vector fields of the functions $H$ and $H\cdot\varepsilon:
T^*Q\rightarrow \mathbb{R}, $ respectively.
The equation $T\gamma\cdot \tilde{X}^\varepsilon= X_H\cdot
\varepsilon ,$ is called the Type II of Hamilton-Jacobi equation
for the RCH system $(T^*Q,\omega,H,F,W)$ with a control law $u$.
\end{theo}
\noindent{\bf Proof: } Since $\tilde{X}=X_{(T^\ast Q,\omega,H,F,u)}=X_H
+\textnormal{vlift}(F)+\textnormal{vlift}(u), $ and
$T\pi_{Q}\cdot \textnormal{vlift}(F)=T\pi_{Q}\cdot \textnormal{vlift}(u)=0, $
then we have that $T\pi_{Q}\cdot \tilde{X}\cdot \varepsilon=T\pi_{Q}\cdot X_H\cdot \varepsilon. $
If we take that $v= X_H\cdot \varepsilon \in TT^* Q, $ and for
any $w \in TT^* Q, \; T\lambda(w)\neq 0, $ from Lemma 2.5 we have that
\begin{align*}
&\omega(T\gamma \cdot \tilde{X}^\varepsilon, \; w)
= \omega(T\gamma \cdot
T\pi_Q\cdot \tilde{X}\cdot \varepsilon, \; w)= \omega(T\gamma \cdot
T\pi_Q\cdot X_H\cdot \varepsilon, \; w)\\ &
= \omega(T(\gamma \cdot
\pi_Q)\cdot X_H\cdot \varepsilon, \; w)= \omega(X_H\cdot \varepsilon, \;
w-T(\gamma \cdot \pi_Q)\cdot
w)-\mathbf{d}\gamma(T\pi_{Q}(X_H\cdot \varepsilon), \; T\pi_{Q}(w))\\
& =\omega(X_H\cdot \varepsilon, \; w) - \omega(X_H\cdot \varepsilon, \;
T\lambda \cdot w)+\lambda^*\omega(X_H\cdot \varepsilon, \; w)\\
& =\omega(X_H\cdot \varepsilon, \; w) - \omega(X_H\cdot \varepsilon, \;
T\lambda \cdot w)+ \omega(T\lambda \cdot X_H\cdot \varepsilon, \; T\lambda \cdot w).
\end{align*}
Because $\varepsilon: T^* Q
\rightarrow T^* Q $ is symplectic, and hence $ X_H\cdot \varepsilon=
T\varepsilon \cdot X_{H\cdot\varepsilon}, $ along $\varepsilon$. Note that
$T\lambda \cdot X_H\cdot \varepsilon=T\gamma \cdot
T\pi_Q\cdot X_H\cdot \varepsilon=T\gamma \cdot
T\pi_Q\cdot \tilde{X}\cdot \varepsilon=T\lambda\cdot \tilde{X}\cdot \varepsilon.$
From the above arguments, we can obtain that
\begin{align*}
&\omega(T\gamma \cdot \tilde{X}^\varepsilon, \; w)- \omega(X_H\cdot \varepsilon, \; w)\\
& =- \omega(X_H\cdot \varepsilon, \;
T\lambda \cdot w)+ \omega(T\lambda \cdot X_H\cdot \varepsilon, \; T\lambda \cdot w)\\
& =-\omega(T\varepsilon \cdot X_{H\cdot\varepsilon}, \; T\lambda \cdot w)+ \omega(T\lambda \cdot \tilde{X}\cdot \varepsilon, \; T\lambda \cdot w)\\
& = \omega(T\lambda \cdot \tilde{X}\cdot \varepsilon -T\varepsilon \cdot X_{H\cdot\varepsilon}, \; T\lambda \cdot w).
\end{align*}
Because the symplectic form $\omega$ is non-degenerate,
it follows that $T\gamma\cdot \tilde{X}^\varepsilon= X_H\cdot
\varepsilon ,$ is equivalent to $T\varepsilon \cdot X_{H\cdot\varepsilon} = T\lambda\cdot \tilde{X}\cdot \varepsilon $.
Thus, $\varepsilon$ is a solution of the equation
$T\varepsilon\cdot X_{H\cdot\varepsilon}= T\lambda \cdot \tilde{X} \cdot\varepsilon,$ if and only if it is a solution of
the Type II of Hamilton-Jacobi equation $T\gamma\cdot \tilde{X}^\varepsilon= X_H\cdot
\varepsilon .$
\hskip 0.3cm $\blacksquare$

\begin{rema}
In particular, if both the external force and control of an RCH
system $(T^*Q,\omega,H,F,u)$ are zero, in this case the RCH system
is just a Hamiltonian system $(T^*Q,\omega,H)$, and from the proofs of
the above Theorem 2.6 and Theorem 2.7, we can obtain two types of geometric Hamilton-Jacobi
theorem for the associated Hamiltonian system, which is given in Wang \cite{wa17}.
Thus, Theorem 2.6 and Theorem 2.7 can be regarded as an extension of two types of Hamilton-Jacobi
theorem for a Hamiltonian system to that for the system with external force and
control. \end{rema}

Next, we note that when an RCH system is given, the force map $F$ is
determined, but the feedback control law $u: T^\ast Q\rightarrow W$
could be chosen. In order to describe the feedback control law to
modify the structure of the RCH system, the controlled Hamiltonian matching
conditions and RCH-equivalence are induced as follows.
\begin{defi}
(RCH-equivalence) Suppose that we have two RCH systems $(T^\ast
Q_i,\omega_i,H_i,F_i,W_i),$ $ i= 1,2,$ we say them to be
RCH-equivalent, or simply, $(T^\ast
Q_1,\omega_1,H_1,F_1,W_1)\stackrel{RCH}{\sim}\\ (T^\ast
Q_2,\omega_2,H_2,F_2,W_2)$, if there exists a diffeomorphism
$\varphi: Q_1\rightarrow Q_2$, such that the following controlled Hamiltonian
matching conditions hold:\\
\noindent {\bf RCH-1:} The cotangent lift map of $\varphi$, that is,
$\varphi^\ast= T^\ast \varphi:T^\ast Q_2\rightarrow T^\ast Q_1$ is
symplectic, and $W_1=\varphi^\ast (W_2).$\\
\noindent {\bf RCH-2:}
$Im[(X_{H_1}+\textnormal{vlift}(F_1)-T\varphi^\ast (X_{H_2})-\textnormal{vlift}(\varphi^\ast
F_2\varphi_\ast)]\subset \textnormal{vlift}(W_1)$, where the map
$\varphi_\ast=(\varphi^{-1})^\ast: T^\ast Q_1\rightarrow T^\ast
Q_2$, and $(\varphi^\ast)_\ast=(\varphi_\ast)^\ast=T^\ast
\varphi_\ast: T^\ast T^\ast Q_2\rightarrow T^\ast T^\ast Q_1$, and
$Im$ means the pointwise image of the map in brackets.
\end{defi}

It is worthy of noting that our RCH system is defined by using the
symplectic structure on the cotangent bundle of a configuration
manifold, we must keep with the symplectic structure when we define
the RCH-equivalence, that is, the induced equivalent map $\varphi^*$
is symplectic on the cotangent bundle. Moreover, the following
Theorem 2.10 explains the
significance of the above RCH-equivalence relation, its proof is
given in Marsden et al \cite{mawazh10}.
\begin{theo}
Suppose that two RCH systems $(T^\ast Q_i,\omega_i,H_i,F_i,W_i)$,
$i=1,2,$ are RCH-equivalent, then there exist two control laws $u_i:
T^\ast Q_i \rightarrow W_i, \; i=1,2, $ such that the two
closed-loop systems produce the same equations of motion, that is,
$X_{(T^\ast Q_1,\omega_1,H_1,F_1,u_1)}\cdot \varphi^\ast
=T(\varphi^\ast) \cdot X_{(T^\ast Q_2,\omega_2,H_2,F_2,u_2)}$, where the
map $T(\varphi^\ast):TT^\ast Q_2\rightarrow TT^\ast Q_1$ is the
tangent map of $\varphi^\ast$. Moreover, the explicit relation
between the two control laws $u_i, i=1,2$ is given by
\begin{equation}\textnormal{vlift}(u_1) -\textnormal{vlift}(\varphi^\ast
u_2\varphi_\ast)=-X_{H_1}
-\textnormal{vlift}(F_1)+T\varphi^\ast (X_{H_2})+\textnormal{vlift}(\varphi^\ast F_2
\varphi_\ast). \; \label{2.4}
\end{equation}
\end{theo}

Moreover, if considering the RCH-equivalence of the RCH
systems, by using the above theorems,
we can prove the following Theorem 2.11, which states that
the solutions of two types of Hamilton-Jacobi equation for an RCH system leave
invariant under the conditions of RCH-equivalence.
\begin{theo}
Suppose that two RCH systems $(T^\ast Q_i,\omega_i,H_i,F_i,W_i)$,
$i=1,2,$ are RCH-equivalent with an equivalent map $\varphi: Q_1
\rightarrow Q_2 $, under the hypotheses and notations of Theorem 2.6,
Theorem 2.7 and Theorem 2.10, we have that\\

\noindent $(\mathrm{i})$ If the one-form $\gamma_2: Q_2 \rightarrow T^* Q_2$ is closed with
respect to $T\pi_{Q_2}: TT^* Q_2 \rightarrow TQ_2, $ then $\gamma_1=
\varphi^* \cdot \gamma_2\cdot \varphi: Q_1 \rightarrow T^* Q_1 $ is also closed with respect to $T\pi_{Q_1}:
TT^* Q_1 \rightarrow TQ_1, $ and hence it is
a solution of the Type I of Hamilton-Jacobi equation for the RCH system
$(T^*Q_1,\omega_1,H_1,F_1,W_1). $ Vice versa;\\

\noindent $(\mathrm{ii})$ If the symplectic map $\varepsilon_2: T^*Q_2\rightarrow T^* Q_2$ is
a solution of the Type II of Hamilton-Jacobi equation for the RCH system
$(T^*Q_2,\omega_2,H_2, F_2,W_2)$, then $\varepsilon_1=
\varphi^* \cdot \varepsilon_2\cdot \varphi_*: T^*Q_1 \rightarrow
T^* Q_1 $ is a symplectic map, and it is a solution of the Type II of Hamilton-Jacobi equation for the RCH
system $(T^*Q_1,\omega_1,H_1,F_1,W_1). $ Vice versa.
\end{theo}
\noindent{\bf Proof: }
We first prove the assertion $(\mathrm{i})$.
If two RCH systems $(T^\ast Q_i,\omega_i,H_i,F_i,W_i)$, $i=1,2,$ are
RCH-equivalent with an equivalent map $\varphi: Q_1 \rightarrow Q_2
$, from Theorem 2.10 we know that there exist two control laws $u_i:
T^\ast Q_i \rightarrow W_i, \; i=1,2, $ such that $X_{(T^\ast
Q_1,\omega_1,H_1,F_1,u_1)}\cdot \varphi^\ast =T(\varphi^\ast)
X_{(T^\ast Q_2,\omega_2,H_2,F_2,u_2)}$, that is, $\tilde{X}_1 \cdot
\varphi^\ast =T(\varphi^\ast)\cdot \tilde{X}_2, $ where $\tilde{X}_i= X_{(T^\ast
Q_i,\omega_i, H_i, F_i, u_i)}, \; i=1,2.$ From the following
commutative Diagram-3:
\[
\begin{CD}
 Q_1 @> \gamma_1 >> T^* Q_1 @> \tilde{X}_1 >> TT^* Q_1 @> T\pi_{Q_1} >> TQ_1 \\
@V \varphi VV @A \varphi^* AA @A T\varphi^* AA @V T\varphi VV \\
 Q_2 @> \gamma_2 >> T^* Q_2 @> \tilde{X}_2 >> TT^* Q_2 @> T\pi_{Q_2} >> TQ_2
\end{CD}
\]
$$\mbox{Diagram-3}$$
we have that $\gamma_1= \varphi^* \cdot \gamma_2\cdot \varphi, $
$\mathbf{d}\gamma_1 = \varphi^* \cdot \mathbf{d}\gamma_2 \cdot
\varphi, $ and
$T\varphi \cdot T\pi_{Q_1} \cdot T\varphi^*= T\pi_{Q_2}. $
For $x\in Q_1, $ and $v, \; w \in TT^*_x Q_1, $ then $\varphi(x) \in Q_2 $
and $T\varphi_*(v), \; T\varphi_*(w) \in TT^*_{\varphi(x)}Q_2. $
Since the one-form $\gamma_2: Q_2 \rightarrow T^* Q_2$ is closed with
respect to $T\pi_{Q_2}: TT^* Q_2 \rightarrow TQ_2, $ then
$$\mathbf{d}\gamma_2 (T\pi_{Q_2}\cdot T\varphi_*(v), \; T\pi_{Q_2}\cdot T\varphi_*(w) )(\varphi(x))=0.$$
Thus,
\begin{align*}
\mathbf{d}\gamma_1 (T\pi_{Q_1}(v), \; T\pi_{Q_1}(w)) (x) &= \varphi^* \cdot \mathbf{d}\gamma_2 \cdot
\varphi (T\pi_{Q_1}(v), \; T\pi_{Q_1}(w))(x)\\
&=  \mathbf{d}\gamma_2 (T\varphi \cdot T\pi_{Q_1}(v),\; T\varphi \cdot T\pi_{Q_1}(w))(\varphi(x))\\
&= \mathbf{d}\gamma_2 (T\varphi \cdot T\pi_{Q_1} \cdot T\varphi^*\cdot T(\varphi^{-1})^*(v), \;
T\varphi \cdot T\pi_{Q_1} \cdot T\varphi^*\cdot T(\varphi^{-1})^*(w)) (\varphi(x))\\
&= \mathbf{d}\gamma_2 (T\pi_{Q_2}\cdot T\varphi_*(v), \; T\pi_{Q_2}\cdot T\varphi_*(w) )(\varphi(x))=0,
\end{align*}
that is, the one-form $\gamma_1=
\varphi^* \cdot \gamma_2\cdot \varphi: Q_1 \rightarrow T^* Q_1 $ is closed with respect to $T\pi_{Q_1}:
TT^* Q_1 \rightarrow TQ_1. $ Moreover, from Theorem 2.6 we know that,
the one-form $\gamma_2$ is
a solution of the Type I of Hamilton-Jacobi equation for the RCH system
$(T^*Q_2,\omega_2, H_2, F_2, W_2), $ that is,
$T\gamma_2\cdot \tilde{X}_2^{\gamma_2}= X_{H_2}\cdot \gamma_2 ,$
where $\tilde{X}_i^{\gamma_i}= T\pi_{Q_i}\cdot \tilde{X}_i \cdot \gamma_i, \; i=1,2.$
Hence,
\begin{align*}
T\gamma_1\cdot \tilde{X}_1^{\gamma_1}& =T(\varphi^* \cdot \gamma_2\cdot \varphi)\cdot
T\pi_{Q_1}\cdot \tilde{X}_1 \cdot \gamma_1 \\
& =T(\varphi^*) \cdot T\gamma_2\cdot T\varphi \cdot T\pi_{Q_1}\cdot \tilde{X}_1 \cdot (\varphi^* \cdot \gamma_2\cdot \varphi) \\
& =T(\varphi^*) \cdot T\gamma_2\cdot T\varphi \cdot T\pi_{Q_1}\cdot (T(\varphi^\ast) \cdot \tilde{X}_2) \cdot \gamma_2\cdot \varphi \\
& =T(\varphi^*) \cdot T\gamma_2\cdot (T\pi_{Q_2}\cdot \tilde{X}_2 \cdot \gamma_2)\cdot \varphi
=T(\varphi^*) \cdot T\gamma_2\cdot \tilde{X}_2^{\gamma_2}\cdot \varphi \\
& =T(\varphi^*) \cdot X_{H_2}\cdot \gamma_2 \cdot \varphi =X_{H_1}\cdot \varphi^*\cdot \gamma_2 \cdot \varphi
=X_{H_1}\cdot \gamma_1,
\end{align*}
where we have used that $T(\varphi^*) \cdot X_{H_2}= X_{H_1}\cdot \varphi^*$,
because $\varphi^*: T^* Q_2 \rightarrow T^* Q_1 $ is symplectic.
Thus, the one-form $\gamma_1=
\varphi^* \cdot \gamma_2\cdot \varphi $ is
a solution of the Type I of Hamilton-Jacobi equation for the RCH system
$(T^*Q_1,\omega_1,H_1,F_1,W_1). $
Note that the map $\varphi: Q_1 \rightarrow Q_2 $ is a diffeomorphism,
and $\varphi^*: T^* Q_2 \rightarrow T^* Q_1 $ is a symplectic isomorphisms,
vice versa. It follows that the assertion $(\mathrm{i})$ of Theorem 2.11 holds.\\

Next, we prove the assertion $(\mathrm{ii})$.
From the following commutative Diagram-4:
\[
\begin{CD}
 Q_1 @> \gamma_1 >> T^* Q_1 @> \varepsilon_1 >> T^* Q_1 @> \tilde{X}_1 >> TT^* Q_1 @> T\pi_{Q_1} >> TQ_1 \\
@V \varphi VV @V \varphi_* VV @ A \varphi^* AA @A T\varphi^* AA @V T\varphi VV \\
 Q_2 @> \gamma_2 >> T^* Q_2 @> \varepsilon_2 >> T^* Q_2 @> \tilde{X}_2 >> TT^* Q_2 @> T\pi_{Q_2} >> TQ_2
\end{CD}
\]
$$\mbox{Diagram-4}$$
we have that $\varepsilon_1=
\varphi^* \cdot \varepsilon_2\cdot \varphi_*: T^*Q_1 \rightarrow
T^* Q_1 .$ Since
$\varepsilon_2: T^*Q_2\rightarrow T^* Q_2$ is symplectic with respect to $\omega_2$, then
for $x\in Q_1 $, $v, \; w \in TT^*_x Q_1, $ and $\varphi(x) \in Q_2 $,
$T\varphi_*(v), \; T\varphi_*(w) \in TT^*_{\varphi(x)}Q_2, $ we have that
$\varepsilon_2^*\cdot \omega_2(T\varphi_*(v), \; T\varphi_*(w))(\varphi(x))=\omega_2 (T\varphi_*(v), \; T\varphi_*(w))(\varphi(x)).$
Note that $\varphi^\ast: T^\ast Q_2\rightarrow T^\ast Q_1$ is symplectic,
$(\varphi^*)^*\omega_1(v,w)(x)= \omega_2 (T\varphi_*(v), \; T\varphi_*(w))(\varphi(x))$,
then we have that
\begin{align*}
\varepsilon_1^*\cdot \omega_1(v,w)(x)& = (\varphi^* \cdot \varepsilon_2\cdot \varphi_*)^*\omega_1(v,w)(x)
= (\varphi_*)^* \cdot \varepsilon_2^*\cdot (\varphi^*)^*\omega_1(v,w)(x)\\
& =(\varphi_*)^* \cdot \varepsilon_2^*\cdot \omega_2(T\varphi_*(v), \; T\varphi_*(w))(\varphi(x))
=((\varphi^{-1})^*)^*\cdot \omega_2(T\varphi_*(v), \; T\varphi_*(w))(\varphi(x))\\
& =\omega_1(T(\varphi^{-1})_*\cdot T\varphi_*(v), \; T(\varphi^{-1})_*\cdot T\varphi_*(w))(\varphi^{-1}\cdot \varphi(x))
= \omega_1(v,w)(x),
\end{align*}
that is, the map $\varepsilon_1=
\varphi^* \cdot \varepsilon_2\cdot \varphi_*: T^*Q_1 \rightarrow
T^* Q_1 $ is symplectic map with respect to $\omega_1$. Moreover, because
the symplectic map $\varepsilon_2: T^*Q_2\rightarrow T^* Q_2$ is
a solution of the Type II of Hamilton-Jacobi equation for the RCH system
$(T^*Q_2,\omega_2,H_2, F_2,W_2)$, that is,
$T\gamma_2\cdot \tilde{X}_2^{\varepsilon_2}= X_{H_2}\cdot
\varepsilon_2 ,$ where $\tilde{X}_i^{\varepsilon_i}= T\pi_{Q_i}\cdot \tilde{X}_i \cdot \varepsilon_i, \; i=1,2.$
Hence, we have that
\begin{align*}
T\gamma_1\cdot \tilde{X}_1^{\varepsilon_1}& =T(\varphi^* \cdot \gamma_2\cdot \varphi)\cdot
T\pi_{Q_1}\cdot \tilde{X}_1 \cdot \varepsilon_1 \\
& =T(\varphi^*) \cdot T\gamma_2\cdot T\varphi \cdot T\pi_{Q_1}\cdot \tilde{X}_1 \cdot (\varphi^* \cdot \varepsilon_2\cdot \varphi_*) \\
& =T(\varphi^*) \cdot T\gamma_2\cdot T\varphi \cdot T\pi_{Q_1}\cdot (T(\varphi^\ast) \cdot \tilde{X}_2) \cdot \varepsilon_2\cdot \varphi_* \\
& =T(\varphi^*) \cdot T\gamma_2\cdot (T\pi_{Q_2}\cdot \tilde{X}_2 \cdot \varepsilon_2)\cdot \varphi_*
=T(\varphi^*) \cdot T\gamma_2\cdot \tilde{X}_2^{\varepsilon_2}\cdot \varphi_* \\
& =T(\varphi^*) \cdot X_{H_2}\cdot \varepsilon_2 \cdot \varphi_* =X_{H_1}\cdot \varphi^*\cdot \varepsilon_2 \cdot \varphi_*
=X_{H_1}\cdot \varepsilon_1,
\end{align*}
that is, the symplectic map $\varepsilon_1=
\varphi^* \cdot \varepsilon_2\cdot \varphi_*$
is a solution of the Type II of Hamilton-Jacobi equation for the RCH
system $(T^*Q_1,\omega_1,H_1,F_1,W_1). $
In the same way, because the map $\varphi: Q_1 \rightarrow Q_2 $ is a diffeomorphism,
and $\varphi^\ast: T^\ast Q_2\rightarrow T^\ast Q_1$ is a symplectic isomorphisms,
vice versa. Hence we prove the assertion $(\mathrm{ii})$ of Theorem 2.11.
\hskip 0.3cm  $\blacksquare$\\

In the following $\S4$ and $\S5$, we shall generalize the above results to regular
point and regular orbit reducible RCH systems with symmetries and momentum maps, and
give precisely a variety of Hamilton-Jacobi equations for the regular reduced RCH
systems.

\section{Magnetic Symplectic Form and Hamilton-Jacobi Theorem}

From the proofs of Theorem 2.6 and Theorem 2.7 we know that Lemma 2.5
is a very important tool. But, if $\omega$ is not the
canonical symplectic form on $T^*Q$, then we cannot prove Lemma 2.5,
and hence we cannot get Theorem 2.6 and Theorem 2.7, that is, the two types of
Hamilton-Jacobi equation for an RCH system.
In order to describe the impact of different geometric structures
for the Hamilton-Jacobi theorem, we shall consider the
magnetic symplectic form on $T^*Q$ as follows:
Assume that $T^*Q$ with the canonical symplectic form $\omega$,
and $B$ is a closed two-form on $Q$,
then $\omega_B= \omega- \pi_Q^*B$ is a symplectic form on $T^*Q$,
where $\pi_Q^*: T^*Q \rightarrow T^*T^*Q $. The
$\omega_B$ is called a magnetic symplectic form, and $\pi_Q^*B$ is called
a magnetic term on $T^*Q$, see Marsden et al. \cite{mamiorpera07}.
Using Lemma 2.5 we can prove the following lemma.
\begin{lemm}
Assume that $\gamma: Q \rightarrow T^*Q$ is an one-form on $Q$, and
$\lambda=\gamma \cdot \pi_{Q}: T^* Q \rightarrow T^* Q .$
For the magnetic symplectic form $\omega_B= \omega- \pi_Q^*B $ on $T^*Q$,
where $\omega$ is the canonical symplectic form on $T^*Q$,
then we have that the following two assertions hold.\\
\noindent $(\mathrm{i})$ For any $v, w \in
TT^* Q, \; \lambda^*\omega_B(v,w)= -(\mathbf{d}\gamma+B)(T\pi_{Q}(v), \;
T\pi_{Q}(w))$; \\
\noindent $(\mathrm{ii})$ For any $v, w \in TT^* Q, \;
\omega_B(T\lambda \cdot v,w)= \omega_B(v, w-T\lambda \cdot
w)-(\mathbf{d}\gamma+B)(T\pi_{Q}(v), \; T\pi_{Q}(w)). $
\end{lemm}
\noindent{\bf Proof:}
We first prove the assertion $(\mathrm{i})$.
Since $\omega$ is the canonical symplectic form on $T^*Q$,
from Lemma 2.5(i) we have that for any $v, w \in
TT^* Q $,
\begin{align*}
\lambda^*\omega_B(v,w)& =\lambda^*\omega(v,w)-\lambda^*\cdot \pi_Q^*B(v,w)\\
& =-\mathbf{d}\gamma(T\pi_{Q}(v), \; T\pi_{Q}(w))
-(\pi_Q\cdot \gamma \cdot \pi_{Q})^*B(v,w)\\
& =-\mathbf{d}\gamma(T\pi_{Q}(v), \; T\pi_{Q}(w))- \pi_Q^*B(v,w)\\
& =-(\mathbf{d}\gamma+B)(T\pi_{Q}(v), \; T\pi_{Q}(w)),
\end{align*}
where we have used the relation $\pi_Q\cdot \gamma\cdot \pi_Q= \pi_Q. $
It follows that the assertion $(\mathrm{i})$ holds.\\

Next, we prove the assertion $(\mathrm{ii})$. For any $v, w \in TT^*
Q,$ from Lemma 2.5(ii) we have that
\begin{align*}
& \omega_B(T\lambda \cdot v,w)= \omega(T\lambda \cdot v,w)-\pi_Q^*B(T\lambda \cdot v,w)\\
& =\omega(v, w-T\lambda \cdot w)-\mathbf{d}\gamma(T\pi_{Q}(v), \; T\pi_{Q}(w))
-B(T\pi_Q\cdot T\lambda \cdot v, \; T\pi_{Q}(w))\\
& =\omega_B(v, w-T\lambda \cdot w)+\pi_Q^*B(v, w-T\lambda \cdot w)\\
& \;\;\;\;\; -\mathbf{d}\gamma(T\pi_{Q}(v), \; T\pi_{Q}(w))
-B(T(\pi_Q\cdot \lambda) \cdot v, \; T\pi_{Q}(w))\\
& =\omega_B(v, w-T\lambda \cdot w)+\pi_Q^*B(v, w)-B(T\pi_{Q}(v), \; T\pi_Q\cdot T\lambda \cdot w)\\
& \;\;\;\;\; -\mathbf{d}\gamma(T\pi_{Q}(v), \; T\pi_{Q}(w))
-B(T(\pi_Q\cdot \gamma \cdot \pi_{Q}) \cdot v, \; T\pi_{Q}(w))\\
& =\omega_B(v, w-T\lambda \cdot w)+\pi_Q^*B(v, w)-B(T\pi_{Q}(v), \; T(\pi_Q\cdot \lambda) \cdot w)\\
& \;\;\;\;\; -\mathbf{d}\gamma(T\pi_{Q}(v), \; T\pi_{Q}(w))
-B(T\pi_Q (v), \; T\pi_{Q}(w))\\
& =\omega_B(v, w-T\lambda \cdot w)+\pi_Q^*B(v, w)-B(T\pi_{Q}(v), \; T\pi_{Q}(w))
-(\mathbf{d}\gamma+B)(T\pi_{Q}(v), \; T\pi_{Q}(w))\\
& =\omega_B(v, w-T\lambda \cdot
w)-(\mathbf{d}\gamma+B)(T\pi_{Q}(v), \; T\pi_{Q}(w)).
\end{align*}
Thus, the assertion $(\mathrm{ii})$ holds.
\hskip 0.3cm $\blacksquare$\\

A magnetic Hamiltonian system is a 3-tuple $(T^\ast Q,\omega_B,H)$,
which is Hamiltonian system defined by the
magnetic symplectic form $\omega_B$. For a given Hamiltonian $H$,
the dynamical vector field $X^B_H$, which is called
the magnetic Hamiltonian vector field,
satisfies the magnetic Hamilton's equation, that is,
$\mathbf{i}_{X^B_{H} }\omega_B= \mathbf{d}H $, see Marsden et al. \cite{mamiorpera07}.
A controlled magnetic Hamiltonian (CMH) system on $T^*Q$ is
a 5-tuple $(T^\ast Q,\omega_B,H,F,W)$, which is a
magnetic Hamiltonian system $(T^\ast Q,\omega_B,H)$
with external force $F$ and control $W$, where $F: T^*Q\rightarrow T^*Q$ is
the fiber-preserving map, and $W\subset T^*£Ñ$ is a fiber submanifold,
see Wang \cite{wa15a}.
Thus, a CMH system is also an RCH system, but its symplectic structure
is given by a magnetic symplectic form, and the set of
the CMH systems is a subset of the set of the RCH systems.
For the given Hamiltonian $H$,
the external force map $F$, and the control subset $W$, assume that
the dynamical vector field of the CMH system $(T^*Q,\omega_B,H,F,W)$ with a control law $u\in W$
can be expressed by
\begin{equation}
X^B_{(T^\ast Q,\omega_B,H,F,u)}= X^B_{H}+\textnormal{vlift}(F)^B
+\textnormal{vlift}(u)^B, \; \label{3.1}
\end{equation}
where $X^B_H$ is the magnetic Hamiltonian vector field.
and $\textnormal{vlift}(F)^B=\textnormal{vlift}(F)X^B_H$,
$\textnormal{vlift}(u)^B=\textnormal{vlift}(u)X^B_H.$ are the
changes of $X^B_H$ under the actions of $F$ and $u$.\\

In the following we shall give precisely the geometric constraint conditions of
magnetic symplectic form for the
dynamical vector field of the CMH system, that is, Type I and Type II of
Hamilton-Jacobi equation for the CMH system.
Assume that $\gamma: Q
\rightarrow T^*Q$ is an one-form on $Q$,
we say that the $\gamma$ satisfies the condition that
$\mathbf{d}\gamma=-B$ with respect to $T\pi_{Q}:
TT^* Q \rightarrow TQ, $ if for any $v, w \in TT^* Q, $ we have
$(\mathbf{d}\gamma +B)(T\pi_{Q}(v),T\pi_{Q}(w))=0. $
For a given CMH system $(T^*Q,\omega_B,H,F,W)$ on $T^*Q$, by using
the above Lemma 3.1, magnetic symplectic form and the dynamical vector field $X^B_{(T^\ast Q,\omega_B,H,F,u)}$,
we can prove the following two types of geometric
Hamilton-Jacobi theorem for the CMH system.
\begin{theo}
(Type I of Hamilton-Jacobi Theorem for a CMH System)
For the CMH system $(T^*Q,\omega_B,H,F,W)$ with the
magnetic symplectic form $\omega_B= \omega- \pi_Q^*B $ on $T^*Q$, assume that $\gamma: Q
\rightarrow T^*Q$ is an one-form on $Q$, and
$\tilde{X}^\gamma = T\pi_{Q}\cdot \tilde{X} \cdot \gamma$,
where $\tilde{X}=X^B_{(T^\ast Q,\omega_B,H,F,u)}$ is the dynamical vector field
of the CMH system $(T^*Q,\omega_B,H,F,W)$ with a control law $u$.
If the one-form $\gamma: Q \rightarrow T^*Q $ satisfies the condition that
$\mathbf{d}\gamma=-B $ with respect to $T\pi_{Q}:
TT^* Q \rightarrow TQ, $ then $\gamma$ is a solution of the equation
$T\gamma\cdot \tilde{X}^\gamma= X^B_H\cdot \gamma ,$ where $X^B_H$
is the magnetic Hamiltonian vector field
of the corresponding magnetic Hamiltonian system $(T^*Q,\omega_B,H),$
and the equation is called the Type I of
Hamilton-Jacobi equation for the CMH system
$(T^*Q,\omega_B,H,F,W)$ with a control law $u$.
Here the maps involved in the theorem are shown
in the following Diagram-5.
\begin{center}
\hskip 0cm \xymatrix{ & T^* Q \ar[d]_{X^B_H}\ar[r]^{\pi_Q}
 & Q \ar[d]_{\tilde{X}^\gamma} \ar[r]^{\gamma} & T^*Q \ar[d]^{\tilde{X}} \\
 & T(T^*Q) & TQ \ar[l]_{T\gamma} & T(T^* Q)\ar[l]_{T\pi_Q}}
\end{center}
$$\mbox{Diagram-5}$$
\end{theo}
\noindent{\bf Proof: }
Since $\tilde{X}=\tilde{X}^B_{(T^\ast Q,\omega,H,F,u)}=X^B_H
+\textnormal{vlift}(F)^B+\textnormal{vlift}(u)^B, $ and
$T\pi_{Q}\cdot \textnormal{vlift}(F)^B=T\pi_{Q}\cdot \textnormal{vlift}(u)^B=0, $
then we have that $T\pi_{Q}\cdot \tilde{X}\cdot \gamma=T\pi_{Q}\cdot X^B_H\cdot \gamma. $
If we take that $v= X^B_H\cdot \gamma \in TT^* Q, $ and for
any $w \in TT^* Q, \; T\pi_{Q}(w)\neq 0, $ from Lemma 3.1(ii) and
$\mathbf{d}\gamma=-B $ with respect to $T\pi_{Q}:
TT^* Q \rightarrow TQ, $ we have that
\begin{align*}
\omega_B(T\gamma \cdot \tilde{X}^\gamma, \; w)&
=\omega_B(T\gamma \cdot T\pi_{Q} \cdot \tilde{X}\cdot\gamma, \; w)
=\omega_B(T\gamma \cdot T\pi_{Q} \cdot X^B_H\cdot\gamma, \; w)\\
& = \omega_B(T(\gamma \cdot \pi_Q)\cdot X^B_H\cdot \gamma, \; w)
= \omega_B(X^B_H\cdot \gamma, \; w-T(\gamma \cdot \pi_Q)\cdot w)\\
& = \omega_B(X^B_H\cdot \gamma, \; w) - \omega_B(X^B_H\cdot \gamma, \;
T\lambda \cdot w).
\end{align*}
Hence, we have that
\begin{equation}
\omega_B(T\gamma \cdot \tilde{X}^\gamma, \; w)- \omega_B(X^B_H\cdot \gamma, \; w)
= -\omega_B(X^B_H\cdot \gamma, \; T\lambda \cdot w). \; \label{3.2}
\end{equation}
If $\gamma$ satisfies the equation $T\gamma\cdot \tilde{X}^\gamma= X^B_H\cdot \gamma ,$
from Lemma 3.1(i) we can obtain that
\begin{align*}
\omega_B(X^B_H\cdot \gamma, \; T\lambda \cdot w) &
= \omega_B(T\gamma \cdot \tilde{X}^\gamma, \; T\lambda \cdot w)\\
&= \omega_B(T\gamma \cdot T\pi_{Q} \cdot \tilde{X}\cdot\gamma, \; T\lambda \cdot w)\\
&= \omega_B(T\gamma \cdot T\pi_{Q} \cdot X^B_H\cdot\gamma, \; T\lambda \cdot w)\\
&= \omega_B(T\lambda \cdot X^B_{H}\cdot\gamma, \; T\lambda \cdot w)\\
&= \lambda^*\omega_B(X^B_{H}\cdot\gamma, \; w)\\
&=-(\mathbf{d}\gamma+B)(T\pi_{Q}\cdot X^B_{H}\cdot\gamma, \; T\pi_{Q}\cdot w)=0,
\end{align*}
since $\gamma: Q \rightarrow T^*Q $ satisfies the condition that
$\mathbf{d}\gamma=-B $ with respect to $T\pi_{Q}:
TT^* Q \rightarrow TQ. $
But, because the magnetic symplectic form $\omega_B$ is non-degenerate,
the left side of (3.2) equals zero, only when
$\gamma$ satisfies the equation $T\gamma\cdot \tilde{X}^\gamma= X^B_H\cdot \gamma .$ Thus,
if the one-form $\gamma: Q \rightarrow T^*Q $ satisfies the condition that
$\mathbf{d}\gamma=-B$ with respect to $T\pi_{Q}:
TT^* Q \rightarrow TQ, $ then $\gamma$ must be a solution of
the Type I of Hamilton-Jacobi equation
$T\gamma\cdot \tilde{X}^\gamma= X^B_H\cdot \gamma ,$ for
the CMH system $(T^*Q,\omega_B,H,F,W)$ with a control law $u$.
\hskip 0.3cm $\blacksquare$\\

It is a natural problem what and how we could do,
if an one-form $\gamma: Q \rightarrow T^*Q $ is not
closed on $Q$ with respect to $T\pi_Q: TT^* Q \rightarrow TQ $ in Theorem 2.6.
Our idea is that we hope to look for a new RCH system,
such that $\gamma$ is a solution of the Type I
of Hamilton-Jacobi equation for the new RCH system.
It is worthy of noting that, if
$\gamma: Q \rightarrow T^*Q $ is not closed on $Q$ with respect to
$T\pi_Q: TT^* Q \rightarrow TQ, $ that is, there exist $v,w \in TT^*Q,$
such that $\mathbf{d}\gamma(T\pi_Q(v), \; T\pi_Q(w))\neq 0,$
and hence $\gamma$ is not yet closed on $Q$.
In this case, we note that
$\mathbf{d}\cdot \mathbf{d}\gamma= \mathbf{d}^2 \gamma =0, $
and hence the $\mathbf{d}\gamma$ is a closed two-form on $Q$.
Thus, we can construct a magnetic symplectic form on $T^*Q$,
that is, $\omega_B= \omega+ \pi_Q^*(\mathbf{d}\gamma), $
where $\omega$ is the canonical symplectic form on $T^*Q$,
and $\pi_Q^*: T^*Q \rightarrow T^*T^*Q $. Moreover,
in this case, for any $x, y \in TQ, $ we have that $(\mathbf{d}\gamma +B)(x, y)=0$,
and hence for any $v, w \in TT^* Q, $ we have
$(\mathbf{d}\gamma +B)(T\pi_{Q}(v),T\pi_{Q}(w))=0, $
that is, the one-form $\gamma: Q \rightarrow T^*Q $ satisfies the condition that
$\mathbf{d}\gamma=-B$ with respect to $T\pi_{Q}:
TT^* Q \rightarrow TQ. $ Thus, we can construct a CMH system,
such that $\gamma$ is just a solution of the Type I
of Hamilton-Jacobi equation for the new CMH system.\\

For a given RCH system $(T^*Q, \omega, H,F,W )$ with
the canonical symplectic form $\omega$ on $T^*Q$, and
$\gamma: Q \rightarrow T^*Q $ is an
one-form on $Q$, and it is not closed with respect to
$T\pi_Q: TT^* Q \rightarrow TQ. $ Then we can construct
a magnetic symplectic form on $T^*Q$,
$\omega_B= \omega+ \pi_Q^*(\mathbf{d}\gamma), $ that is, $B=- \mathbf{d}\gamma,$
and a CMH system $(T^*Q, \omega_B, H,F,W )$,
its dynamical vector field with a control law $u$ is given by
$X^B_{(T^\ast Q,\omega_B,H,F,u)}
=X^B_{H}+\textnormal{vlift}(F)^B +\textnormal{vlift}(u)^B$,
where $X^B_{H}$ satisfies the magnetic Hamiltonian equation, that is,
$\mathbf{i}_{X^B_{H} }\omega_B= \mathbf{d}H $.
Since in this case, the one-form $\gamma: Q \rightarrow T^*Q $ satisfies the condition that
$\mathbf{d}\gamma=-B$ with respect to $T\pi_{Q}:
TT^* Q \rightarrow TQ, $
by using Lemma 3.1 and the dynamical vector field
$X^B_{(T^\ast Q,\omega_B,H,F,u)}$, from Theorem 3.2
we can obtain the following theorem.
\begin{theo}
For a given RCH system $(T^*Q,\omega,H,F,W)$ with
the canonical symplectic form $\omega$ on $T^*Q$,
assume that the one-form $\gamma: Q
\rightarrow T^*Q$ is not closed with respect to
$T\pi_Q: TT^* Q \rightarrow TQ. $ Construct
a magnetic symplectic form on $T^*Q$,
$\omega_B= \omega+ \pi_Q^*(\mathbf{d}\gamma), $
and a CMH system $(T^*Q, \omega_B, H, F, W )$.
Denote $\tilde{X}^\gamma = T\pi_{Q}\cdot \tilde{X} \cdot \gamma$,
where $\tilde{X}= X^B_{(T^\ast Q,\omega_B,H,F,u)}$ is the dynamical vector field
of the CMH system $(T^*Q,\omega_B,H,F,W)$ with a control law $u$.
Then $\gamma$ is a solution of the Type I of
Hamilton-Jacobi equation $T\gamma\cdot \tilde{X}^\gamma= X^B_H\cdot \gamma ,$
for the CMH system $(T^*Q,\omega_B,H,F,W)$ with a control law $u$.
\end{theo}
Next, for any symplectic map $\varepsilon: T^* Q \rightarrow T^* Q $
with respect to the magnetic symplectic form $\omega_B$,
we can prove the following Type II of geometric
Hamilton-Jacobi theorem for the CMH system $(T^*Q,\omega_B,H,F,W)$. For convenience,
the maps involved in the following theorem and its proof are shown
in Diagram-6.
\begin{center}
\hskip 0cm \xymatrix{ & T^* Q \ar[r]^{\varepsilon}
& T^*Q \ar[d]_{X^B_{H\cdot \varepsilon}}
\ar[dr]^{\tilde{X}^\varepsilon} \ar[r]^{\pi_Q}
& Q \ar[r]^{\gamma} & T^*Q \ar[d]^{\tilde{X}} \\
&  & T(T^*Q) & TQ \ar[l]_{T\gamma} & T(T^* Q)\ar[l]_{T\pi_Q}}
\end{center}
$$\mbox{Diagram-6}$$
\begin{theo}
(Type II of Hamilton-Jacobi Theorem for a CMH System)
For the CMH system $(T^*Q,\omega_B,H,F,W)$ with the
magnetic symplectic form $\omega_B= \omega- \pi_Q^*B $ on $T^*Q$,
assume that $\gamma: Q \rightarrow T^*Q$ is an one-form on $Q$, and
$\lambda=\gamma\cdot\pi_{Q}: T^* Q \rightarrow T^* Q $, and for any
symplectic map $\varepsilon: T^* Q \rightarrow T^* Q $ with respect to $\omega_B$,
denote by $\tilde{X}^\varepsilon = T\pi_{Q}\cdot \tilde{X} \cdot \varepsilon$,
where $\tilde{X}=X^B_{(T^\ast Q,\omega_B,H,F,u)}$
is the dynamical vector field of the CMH system
$(T^*Q,\omega_B,H,F,W)$ with a control law $u$.
Then $\varepsilon$ is a solution of the equation
$T\varepsilon\cdot X^B_{H\cdot\varepsilon}= T\lambda \cdot \tilde{X} \cdot \varepsilon,$
if and only if it is a solution of the equation $T\gamma\cdot \tilde{X}^\varepsilon= X^B_H\cdot
\varepsilon, $ where $X^B_H$ and $ X^B_{H\cdot\varepsilon} \in
TT^*Q $ are the magnetic Hamiltonian vector fields of the functions $H$ and $H\cdot\varepsilon:
T^*Q\rightarrow \mathbb{R}, $ respectively.
The equation $T\gamma\cdot \tilde{X}^\varepsilon= X^B_H\cdot
\varepsilon ,$ is called the Type II of Hamilton-Jacobi equation
for the CMH system $(T^*Q,\omega_B,H,F,W)$ with a control law $u$.
\end{theo}
\noindent{\bf Proof: }
Since $\tilde{X}=X^B_{(T^\ast Q,\omega_B,H,F,u)}
=X^B_H +\textnormal{vlift}(F)^B+\textnormal{vlift}(u)^B, $ and
$T\pi_{Q}\cdot \textnormal{vlift}(F)^B=T\pi_{Q}\cdot \textnormal{vlift}(u)^B=0, $
then we have that $T\pi_{Q}\cdot \tilde{X}\cdot \varepsilon
=T\pi_{Q}\cdot X^B_H\cdot \varepsilon. $
If we take that $v= X^B_H\cdot \varepsilon \in TT^* Q, $ and for
any $w \in TT^* Q, \; T\lambda(w)\neq 0, $ from Lemma 3.1(ii) we have that
\begin{align*}
&\omega_B(T\gamma \cdot \tilde{X}^\varepsilon, \; w)
= \omega_B(T\gamma \cdot T\pi_Q\cdot \tilde{X}\cdot \varepsilon, \; w)
= \omega_B(T\gamma \cdot T\pi_Q\cdot X^B_H\cdot \varepsilon, \; w)\\ &
= \omega_B(T(\gamma \cdot \pi_Q)\cdot X^B_H\cdot \varepsilon, \; w)
= \omega_B(X^B_H\cdot \varepsilon, \; w-T(\gamma \cdot \pi_Q)\cdot w)
-(\mathbf{d}\gamma+B)(T\pi_{Q}(X^B_H\cdot \varepsilon), \; T\pi_{Q}(w))\\
& =\omega_B(X^B_H\cdot \varepsilon, \; w) - \omega_B(X^B_H\cdot \varepsilon, \;
T\lambda \cdot w)+\lambda^*\omega_B(X^B_H\cdot \varepsilon, \; w)\\
& =\omega_B(X^B_H\cdot \varepsilon, \; w) - \omega_B(X^B_H\cdot \varepsilon, \;
T\lambda \cdot w)+ \omega_B(T\lambda \cdot X^B_H\cdot \varepsilon, \; T\lambda \cdot w).
\end{align*}
Because $\varepsilon: T^* Q
\rightarrow T^* Q $ is symplectic with respect to $\omega_B$,
and hence $ X^B_H\cdot \varepsilon= T\varepsilon \cdot X^B_{H\cdot\varepsilon}, $
along $\varepsilon$. Note that
$T\lambda \cdot X^B_H\cdot \varepsilon=T\gamma \cdot
T\pi_Q\cdot X^B_H\cdot \varepsilon=T\gamma \cdot
T\pi_Q\cdot \tilde{X}\cdot \varepsilon=T\lambda\cdot \tilde{X}\cdot \varepsilon.$
From the above arguments, we can obtain that
\begin{align*}
&\omega_B(T\gamma \cdot \tilde{X}^\varepsilon, \; w)- \omega_B(X^B_H\cdot \varepsilon, \; w)\\
& =- \omega_B(X^B_H\cdot \varepsilon, \; T\lambda \cdot w)
+ \omega_B(T\lambda \cdot X^B_H\cdot \varepsilon, \; T\lambda \cdot w)\\
& =-\omega_B(T\varepsilon \cdot X^B_{H\cdot\varepsilon}, \; T\lambda \cdot w)
+ \omega_B(T\lambda \cdot \tilde{X}\cdot \varepsilon, \; T\lambda \cdot w)\\
& = \omega_B(T\lambda \cdot \tilde{X}\cdot \varepsilon
-T\varepsilon \cdot X^B_{H\cdot\varepsilon}, \; T\lambda \cdot w).
\end{align*}
Because the magnetic symplectic form $\omega_B$ is non-degenerate,
it follows that $T\gamma\cdot \tilde{X}^\varepsilon= X^B_H\cdot
\varepsilon ,$ is equivalent to $T\varepsilon \cdot X^B_{H\cdot\varepsilon}
= T\lambda\cdot \tilde{X}\cdot \varepsilon $.
Thus, $\varepsilon$ is a solution of the equation
$T\varepsilon\cdot X^B_{H\cdot\varepsilon}= T\lambda \cdot \tilde{X} \cdot\varepsilon,$
if and only if it is a solution of the Type II of Hamilton-Jacobi equation
$T\gamma\cdot \tilde{X}^\varepsilon= X^B_H\cdot \varepsilon .$
\hskip 0.3cm $\blacksquare$

\begin{rema}
When $B=0$, in this case the magnetic symplectic form $\omega_B$
is just the canonical symplectic form $\omega$ on $T^*Q$, and the
condition that the one-form $\gamma: Q \rightarrow T^*Q $ satisfies the condition that
$\mathbf{d}\gamma=-B $ with respect to $T\pi_{Q}:
TT^* Q \rightarrow TQ, $ becomes that  $\gamma $ is closed with respect to
$T\pi_Q: TT^* Q \rightarrow TQ.$ Thus, from above Theorem 3.2 and Theorem 3.4,
we can obtain Theorem 2.6 and Theorem 2.7.
\end{rema}

\section{Hamilton-Jacobi Theorem of the $R_p$-reduced RCH System }

The reduction theory for the mechanical system with symmetry is an
important subject and it is widely studied in the theory of
mathematics and mechanics, as well as applications. The main goal of
reduction theory in mechanics is to use conservation laws and the
associated symmetries to reduce the number of dimensions of a
mechanical system required to be described. So, such reduction
theory is regarded as a useful tool for simplifying and studying
concrete mechanical systems. Over forty years ago, the regular
symplectic reduction for the Hamiltonian system with symmetry and
coadjoint equivariant momentum map was set up by famous professors
Jerrold E. Marsden and Alan Weinstein, which is called
Marsden-Weinstein reduction, and great developments have been
obtained around the work in the theoretical study and applications
of mathematics, mechanics and physics; see
Abraham and Marsden \cite{abma78}, Abraham et al.
\cite{abmara88}, Arnold \cite{ar89}, de Le\'{o}n and Rodrigues \cite{lero89},
Libermann and Marle \cite{lima87}, Marsden \cite{ma92}, Marsden et al.
\cite{mamiorpera07, mamora90}, Marsden and Perlmutter \cite{mape00},
Marsden and Ratiu \cite{mara99}, Marsden and
Weinstein \cite{mawe74}, Meyer \cite{me73},
Nijmeijer and Van der Schaft \cite {nivds90} and Ortega and Ratiu \cite{orra04}.\\

It is worthy of noting that the authors in Marsden et al. \cite{mawazh10}
set up the regular reduction theory for the RCH systems
with symplectic structures and symmetries on a symplectic fiber
bundle, as an extension of the  Marsden-Weinstein reduction theory of
Hamiltonian systems under regular controlled Hamiltonian equivalence
conditions, and from the viewpoint of completeness of regular symplectic
reduction and by analyzing carefully the geometrical and topological
structures of the phase space and the reduced phase space of the
corresponding Hamiltonian system.
Some developments around the work are given in
Wang and Zhang \cite{wazh12}, Ratiu and Wang \cite{rawa12},
Wang \cite{wa15a}.
In this section, we first give the regular point reducible
RCH system with symmetry and momentum map. Then we
give precisely the geometric constraint conditions of the $R_p$-reduced symplectic form for the
dynamical vector field of the regular point reducible RCH system,
and prove the Type I and Type II of
Hamilton-Jacobi theorem for the $R_p$-reduced RCH system.
Moreover, we state the relationship between the solutions of
Type II of Hamilton-Jacobi equations and regular point reduction. Finally,
we consider the RpCH-equivalence, and state that the solutions of
two types of Hamilton-Jacobi equations for the RCH systems with
symmetries leave invariant under the conditions of RpCH-equivalence.
We shall follow the notations and conventions introduced in Marsden
et al \cite{mawazh10}, Wang \cite{wa17}, Wang \cite{wa18}\\

At first, we consider the regular point reducible RCH system,
which is given by Marsden et al \cite{mawazh10}.
Let $Q$ be a smooth manifold and $T^\ast Q$ its cotangent bundle
with the symplectic form $\omega$. Let $\Phi:G\times Q\rightarrow Q$
be a smooth left action of the Lie group $G$ on $Q$, which is free
and proper. Assume that the cotangent lifted left action
$\Phi^{T^\ast}:G\times T^\ast Q\rightarrow T^\ast Q$ is symplectic,
free and proper, and admits an $\operatorname{Ad}^\ast$-equivariant
momentum map $\mathbf{J}:T^\ast Q\rightarrow \mathfrak{g}^\ast$,
where $\mathfrak{g}$ is a Lie algebra of $G$ and $\mathfrak{g}^\ast$
is the dual of $\mathfrak{g}$. Let $\mu\in\mathfrak{g}^\ast$ be a
regular value of $\mathbf{J}$ and denote by $G_\mu$ the isotropy
subgroup of the coadjoint $G$-action at the point
$\mu\in\mathfrak{g}^\ast$, which is defined by $G_\mu=\{g\in
G|\operatorname{Ad}_g^\ast \mu=\mu \}$. Since $G_\mu (\subset G)$
acts freely and properly on $Q$ and on $T^\ast Q$, then
$Q_\mu=Q/G_\mu$ is a smooth manifold and that the canonical
projection $\rho_\mu:Q\rightarrow Q_\mu$ is a surjective submersion.
It follows that $G_\mu$ acts also freely and properly on
$\mathbf{J}^{-1}(\mu)$, so that the space $(T^\ast
Q)_\mu=\mathbf{J}^{-1}(\mu)/G_\mu$ is a symplectic manifold with
the $R_p$-reduced symplectic form $\omega_\mu$ uniquely characterized by the relation
\begin{equation}\pi_\mu^\ast \omega_\mu=i_\mu^\ast
\omega. \label{4.1}
\end{equation} The map
$i_\mu:\mathbf{J}^{-1}(\mu)\rightarrow T^\ast Q$ is the inclusion
and $\pi_\mu:\mathbf{J}^{-1}(\mu)\rightarrow (T^\ast Q)_\mu$ is the
projection. The pair $((T^\ast Q)_\mu,\omega_\mu)$ is called
the $R_p$-reduced space of $(T^\ast Q,\omega)$ at $\mu$,
which is symplectic diffeomorphic to a symplectic fiber bundle. Let
$H: T^\ast Q\rightarrow \mathbb{R}$ be a $G$-invariant Hamiltonian,
the flow $F_t$ of the Hamiltonian vector field $X_H$ leaves the
connected components of $\mathbf{J}^{-1}(\mu)$ invariant and
commutes with the $G$-action, so it induces a flow $f_t^\mu$ on
$(T^\ast Q)_\mu$, defined by $f_t^\mu\cdot \pi_\mu=\pi_\mu \cdot
F_t\cdot i_\mu$, and the vector field $X_{h_\mu}$ generated by the
flow $f_t^\mu$ on $((T^\ast Q)_\mu,\omega_\mu)$ is Hamiltonian with
the associated $R_p$-reduced Hamiltonian function
$h_\mu:(T^\ast Q)_\mu\rightarrow \mathbb{R}$ defined by
$h_\mu\cdot\pi_\mu=H\cdot i_\mu$, and the Hamiltonian vector fields
$X_H$ and $X_{h_\mu}$ are $\pi_\mu$-related.
Moreover, assume that the fiber-preserving map $F:T^\ast Q\rightarrow T^\ast
Q$ and the control subset $W$ of\; $T^\ast Q$ are both $G$-invariant.
In order to get the $R_p$-reduced RCH system, we also assume that
$F(\mathbf{J}^{-1}(\mu))\subset \mathbf{J}^{-1}(\mu)£¬$ and $W \cap
\mathbf{J}^{-1}(\mu)\neq \emptyset $.
Thus, we can introduce a regular point
reducible RCH system as follows, see Marsden et al \cite{mawazh10}
and Wang \cite{wa18}.
\begin{defi}
(Regular Point Reducible RCH System) A 6-tuple $(T^\ast Q, G,
\omega, H, F, W)$ with
the canonical symplectic form $\omega$ on $T^*Q$, where the Hamiltonian $H:T^\ast Q\rightarrow
\mathbb{R}$, the fiber-preserving map $F:T^\ast Q\rightarrow T^\ast
Q$ and the fiber submanifold $W$ of\; $T^\ast Q$ are all
$G$-invariant, is called a regular point reducible RCH system, if
there exists a point $\mu\in\mathfrak{g}^\ast$, which is a regular
value of the momentum map $\mathbf{J}$, such that the regular point
reduced system, that is, the 5-tuple $((T^\ast Q)_\mu,
\omega_\mu,h_\mu,f_\mu,W_\mu)$, where $(T^\ast
Q)_\mu=\mathbf{J}^{-1}(\mu)/G_\mu$, $\pi_\mu^\ast
\omega_\mu=i_\mu^\ast\omega$, $h_\mu\cdot \pi_\mu=H\cdot i_\mu$,
$F(\mathbf{J}^{-1}(\mu))\subset \mathbf{J}^{-1}(\mu) $, $f_\mu\cdot
\pi_\mu=\pi_\mu \cdot F\cdot i_\mu$, $W \cap
\mathbf{J}^{-1}(\mu)\neq \emptyset $, $W_\mu=\pi_\mu(W\cap
\mathbf{J}^{-1}(\mu))$, is an RCH system, which is
simply written as $R_p$-reduced RCH system. Where $((T^\ast
Q)_\mu,\omega_\mu)$ is the $R_p$-reduced space, the function $h_\mu:(T^\ast
Q)_\mu\rightarrow \mathbb{R}$ is called the $R_p$-reduced Hamiltonian, the
fiber-preserving map $f_\mu:(T^\ast Q)_\mu\rightarrow (T^\ast
Q)_\mu$ is called the $R_p$-reduced (external) force map, $W_\mu$ is a
fiber submanifold of \;$(T^\ast Q)_\mu$ and is called the $R_p$-reduced
control subset.
\end{defi}

Denote by $X_{(T^\ast Q,G,\omega,H,F,u)}$ the dynamical vector field of
the regular point reducible RCH system $(T^\ast Q,G,\omega, H,F,W)$
with a control law $u$. Assume that it can be expressed by
\begin{equation}X_{(T^\ast Q,G,\omega,H,F,u)}
=X_H+\textnormal{vlift}(F)+\textnormal{vlift}(u).\label{4.2}
\end{equation}
If an $R_p$-reduced feedback control law $u_\mu:(T^\ast
Q)_\mu\rightarrow W_\mu$ is chosen, the $R_p$-reduced RCH
system $((T^\ast Q)_\mu, \omega_\mu, h_\mu, f_\mu, u_\mu)$ is a
closed-loop regular dynamical system with a control law $u_\mu$.
Assume that its dynamical vector field $X_{((T^\ast Q)_\mu, \omega_\mu, h_\mu,
f_\mu, u_\mu)}$ can be expressed by
\begin{equation}X_{((T^\ast Q)_\mu, \omega_\mu, h_\mu, f_\mu, u_\mu)}
=X_{h_\mu}+\textnormal{vlift}(f_\mu)+\textnormal{vlift}(u_\mu),
\label{4.3}
\end{equation}
where $X_{h_\mu}$ is the Hamiltonian vector field of
the $R_p$-reduced Hamiltonian $h_\mu$, and $\textnormal{vlift}(f_\mu)=
\textnormal{vlift}(f_\mu)X_{h_\mu}$, $\textnormal{vlift}(u_\mu)=
\textnormal{vlift}(u_\mu)X_{h_\mu}$, and satisfies the condition
\begin{equation}X_{((T^\ast Q)_\mu, \omega_\mu, h_\mu, f_\mu,
u_\mu)}\cdot \pi_\mu=T\pi_\mu\cdot X_{(T^\ast
Q,G,\omega,H,F,u)}\cdot i_\mu. \label{4.4}
\end{equation}

From Marsden et al.\cite{mawazh10} and Wang \cite{wa18}, we know that,
the set of Hamiltonian systems with symmetries on the cotangent
bundle is not complete under the Marsden-Weinstein reduction,
then the regular point reduced system of an
RCH system with symmetry defined on the cotangent bundle
$T^*Q$ may not be an RCH system on a cotangent bundle.
On the other hand, from the expression of the dynamical vector
field of an $R_p$-reduced RCH system, we know that under the
actions of the $R_p$-reduced external force $f_\mu$
and the $R_p$-reduced control $u_\mu$, in general, the dynamical vector
field is not Hamiltonian, and the $R_p$-reduced RCH system is not
yet a Hamiltonian system. Thus, we can not describe the Hamilton-Jacobi equation for
the $R_p$-reduced RCH system from the viewpoint of generating
function just like same as Theorem 1.1
given by Abraham and Marsden in \cite{abma78}.
But, for a given regular point reducible RCH system
$(T^*Q,G,\omega,H,F,W)$ with an $R_p$-reduced RCH system $((T^\ast
Q)_\mu, \omega_\mu,h_\mu, f_\mu, u_\mu )$, by using Lemma 2.5,
we can give precisely the following geometric constraint
conditions of the $R_p$-reduced symplectic form for the
dynamical vector field of the regular point reducible RCH system, that is, Type I and Type II of
Hamilton-Jacobi equation for the $R_p$-reduced RCH system $((T^\ast
Q)_\mu, \omega_\mu,h_\mu, f_\mu, u_\mu )$.
At first, by using the fact that the one-form $\gamma: Q
\rightarrow T^*Q $ is closed with respect to
$T\pi_Q: TT^* Q \rightarrow TQ, $ and $\textmd{Im}(\gamma)\subset
\mathbf{J}^{-1}(\mu), $ and it is $G_\mu$-invariant, we can prove the following Type I of
Hamilton-Jacobi theorem for the $R_p$-reduced RCH system $((T^\ast
Q)_\mu, \omega_\mu,h_\mu, f_\mu, u_\mu )$.
For convenience, the maps involved in
the following theorem and its proof are shown in Diagram-7.
\begin{center}
\hskip 0cm \xymatrix{ \mathbf{J}^{-1}(\mu) \ar[r]^{i_\mu} & T^* Q \ar[d]_{X_H} \ar[r]^{\pi_Q}
& Q \ar[d]_{\tilde{X}^\gamma} \ar[r]^{\gamma}
& T^*Q \ar[d]_{\tilde{X}} \ar[r]^{\pi_\mu}
& (T^* Q)_\mu \ar[d]_{X_{h_\mu}} \\
& T(T^*Q)  & TQ \ar[l]^{T\gamma}
& T(T^*Q) \ar[l]^{T\pi_Q} \ar[r]_{T\pi_\mu} & T(T^* Q)_\mu }
\end{center}
$$\mbox{Diagram-7}$$
\begin{theo} (Type I of Hamilton-Jacobi Theorem for an $R_p$-reduced RCH System)
For the regular point reducible RCH system
$(T^*Q,G,\omega,H,F,W)$ with an $R_p$-reduced RCH system \\
$((T^\ast Q)_\mu, \omega_\mu,h_\mu, f_\mu, u_\mu )$,
assume that $\gamma: Q \rightarrow T^*Q$ is an one-form
on $Q$, and $\tilde{X}^\gamma = T\pi_{Q}\cdot \tilde{X} \cdot \gamma$,
where $\tilde{X}=X_{(T^\ast Q,G,\omega,H,F,u)}$ is the dynamical vector
field of the regular point reducible RCH system
$(T^*Q,G,\omega,H,F,W)$ with a control law $u$. Moreover,
assume that $\mu \in \mathfrak{g}^\ast $ is a regular value of the momentum
map $\mathbf{J}$, and $\textmd{Im}(\gamma)\subset
\mathbf{J}^{-1}(\mu), $ and it is $G_\mu$-invariant, and
$\bar{\gamma}=\pi_\mu(\gamma): Q \rightarrow (T^* Q)_\mu. $
If the one-form $\gamma: Q \rightarrow T^*Q $ is closed with respect to
$T\pi_Q: TT^* Q \rightarrow TQ, $
then $\bar{\gamma}$ is a solution of the equation
$T\bar{\gamma}\cdot \tilde{X}^\gamma= X_{h_\mu}\cdot \bar{\gamma}, $
where $X_{h_\mu}$ is the Hamiltonian vector field of the $R_p$-reduced
Hamiltonian function $h_\mu:(T^\ast Q)_\mu\rightarrow \mathbb{R}, $
and the equation is called the Type I of Hamilton-Jacobi equation for
the $R_p$-reduced RCH System $((T^\ast Q)_\mu, \omega_\mu,h_\mu, f_\mu, u_\mu)$.
\end{theo}
\noindent{\bf Proof: } At first, from Theorem 2.6, we know that
$\gamma$ is a solution of the Type I of Hamilton-Jacobi equation
$T\gamma\cdot \tilde{X}^\gamma= X_H\cdot \gamma. $ Next, we note that
$\textmd{Im}(\gamma)\subset \mathbf{J}^{-1}(\mu), $ and it
is $G_\mu$-invariant, in this case $\pi_\mu^*\omega_\mu=
i_\mu^*\omega= \omega, $ along $\textmd{Im}(\gamma)$.
Since $\tilde{X}=X_{(T^\ast Q,G,\omega,H,F,u)}=X_H
+\textnormal{vlift}(F)+\textnormal{vlift}(u), $ and
$T\pi_{Q}\cdot \textnormal{vlift}(F)=T\pi_{Q}\cdot \textnormal{vlift}(u)=0, $
then we have that $T\pi_{Q}\cdot \tilde{X}\cdot \gamma=T\pi_{Q}\cdot X_H\cdot \gamma. $
By using the reduced symplectic form $\omega_\mu$, if we take
that $v= X_H\cdot \gamma \in TT^* Q,$ and for any $w \in TT^* Q, \; T\pi_{Q}(w)\neq 0,$
and $T\pi_{\mu} (w) \neq 0, $ from Lemma 2.5(ii) we have that
\begin{align*}
& \omega_\mu(T\bar{\gamma} \cdot \tilde{X}^\gamma, \; T\pi_\mu \cdot w)
 =\omega_\mu(T(\pi_\mu \cdot \gamma) \cdot T\pi_Q\cdot \tilde{X}\cdot \gamma, \; T\pi_\mu \cdot w )\\
& = \pi_\mu^*\omega_\mu(T\gamma \cdot T\pi_Q\cdot X_H \cdot\gamma, \; w)
= \omega(T(\gamma \cdot \pi_Q)\cdot X_H\cdot \gamma, \; w)\\
& = \omega(X_H\cdot \gamma, \; w-T(\gamma \cdot \pi_Q)\cdot w)
-\mathbf{d}\gamma(T\pi_{Q}(X_H\cdot \gamma), \; T\pi_{Q}(w))\\
& =\omega(X_H\cdot \gamma, \; w) - \omega(X_H\cdot \gamma, \;
T(\gamma \cdot \pi_Q)\cdot w)-\mathbf{d}\gamma(T\pi_{Q}(X_H\cdot \gamma), \; T\pi_{Q}(w))\\
& =\pi_\mu^*\omega_\mu(X_H\cdot
\gamma, \; w) - \pi_\mu^*\omega_\mu(X_H\cdot \gamma, \; T(\gamma \cdot \pi_Q)\cdot w)
-\mathbf{d}\gamma(T\pi_{Q}(X_H\cdot \gamma), \; T\pi_{Q}(w))\\
& = \omega_\mu(T\pi_\mu(X_H\cdot \gamma), \;
T\pi_\mu \cdot w) - \omega_\mu(T\pi_\mu\cdot(X_H\cdot \gamma), \;
T(\pi_\mu \cdot\gamma \cdot \pi_Q) \cdot w) \\
& \;\;\;\; -\mathbf{d}\gamma(T\pi_{Q}(X_H\cdot \gamma), \; T\pi_{Q}(w))\\
& = \omega_\mu(T\pi_\mu(X_H)\cdot \pi_\mu(\gamma), \; T\pi_\mu \cdot w)
- \omega_\mu(T\pi_\mu(X_H)\cdot \pi_\mu(\gamma), \; T\bar{\gamma}\cdot T\pi_{Q}(w))\\
& \;\;\;\; -\mathbf{d}\gamma(T\pi_{Q}(X_H\cdot \gamma), \; T\pi_{Q}(w))\\
& = \omega_\mu(X_{h_\mu} \cdot
\bar{\gamma}, \; T\pi_\mu \cdot w)- \omega_\mu(X_{h_\mu} \cdot
\bar{\gamma}, \; T\bar{\gamma} \cdot T\pi_{Q}(w))-\mathbf{d}\gamma(T\pi_{Q}(X_H\cdot \gamma), \; T\pi_{Q}(w)),
\end{align*}
where we have used that $T\pi_\mu(X_H)= X_{h_\mu}. $
Since the one-form $\gamma: Q \rightarrow T^*Q $ is closed with respect to
$T\pi_Q: TT^* Q \rightarrow TQ, $ then we have that $
\mathbf{d}\gamma(T\pi_{Q}(X_H\cdot \gamma), \; T\pi_{Q}(w))=0, $
and hence
\begin{equation}
 \omega_\mu(T\bar{\gamma} \cdot X_H^\gamma, \; T\pi_\mu \cdot w)-
\omega_\mu(X_{h_\mu} \cdot \bar{\gamma}, \; T\pi_\mu \cdot w)
 = - \omega_\mu(X_{h_\mu} \cdot
\bar{\gamma}, \; T\bar{\gamma} \cdot T\pi_{Q}(w)). \; \label{4.5}
\end{equation}
If $\bar{\gamma}$ satisfies the equation $T\bar{\gamma}\cdot \tilde{X}^\gamma= X_{h_\mu}\cdot \bar{\gamma}, $
from Lemma 2.5(i) we can obtain that
\begin{align*}
- \omega_\mu(X_{h_\mu} \cdot
\bar{\gamma}, \; T\bar{\gamma} \cdot T\pi_{Q}(w))
& = -\omega_\mu (T\bar{\gamma} \cdot \tilde{X}^\gamma, \; T\bar{\gamma} \cdot T\pi_{Q}(w))\\
& = -\bar{\gamma}^*\omega_\mu (T\pi_{Q} \cdot \tilde{X}\cdot\gamma, \; T\pi_Q(w))\\
& = -\gamma^* \cdot \pi_\mu^*\omega_\mu (T\pi_{Q} \cdot X_{H}\cdot\gamma, \; T\pi_{Q}(w))\\
& = -\gamma^*\omega( T\pi_{Q}(X_{H}\cdot\gamma), \; T\pi_{Q}(w))\\
& = \mathbf{d}\gamma(T\pi_{Q}( X_{H}\cdot\gamma ), \; T\pi_{Q}(w))=0.
\end{align*}
Because the reduced symplectic form $\omega_\mu$ is non-degenerate, the left side of (4.5) equals zero, only when
$\bar{\gamma}$ satisfies the equation $T\bar{\gamma}\cdot \tilde{X}^\gamma= X_{h_\mu}\cdot \bar{\gamma}.$ Thus,
if the one-form $\gamma: Q \rightarrow T^*Q $ is closed with respect to
$T\pi_Q: TT^* Q \rightarrow TQ, $ then $\bar{\gamma}$ must be a solution of the Type I of Hamilton-Jacobi equation
$T\bar{\gamma}\cdot \tilde{X}^\gamma= X_{h_\mu}\cdot \bar{\gamma}.$
\hskip 0.3cm $\blacksquare$\\

Next, for any $G_\mu$-invariant symplectic map $\varepsilon: T^* Q \rightarrow T^* Q $,
we can prove the following Type II of
Hamilton-Jacobi theorem for the $R_p$-reduced RCH system $((T^\ast
Q)_\mu, \omega_\mu,h_\mu, f_\mu, u_\mu )$.
For convenience, the maps involved in
the following theorem and its proof are shown in Diagram-8.
\begin{center}
\hskip 0cm \xymatrix{ \mathbf{J}^{-1}(\mu) \ar[r]^{i_\mu} & T^* Q
\ar[d]_{X_{H\cdot \varepsilon}} \ar[dr]^{\tilde{X}^\varepsilon} \ar[r]^{\pi_Q}
& Q \ar[r]^{\gamma} & T^*Q \ar[d]_{\tilde{X}} \ar[dr]^{X_{h_\mu \cdot\bar{\varepsilon}}} \ar[r]^{\pi_\mu} & (T^* Q)_\mu \ar[d]^{X_{h_\mu}} \\
  & T(T^*Q)  & TQ \ar[l]^{T\gamma} & T(T^*Q) \ar[l]^{T\pi_Q} \ar[r]_{T\pi_\mu} & T(T^* Q)_\mu }
\end{center}
$$\mbox{Diagram-8}$$
\begin{theo}
(Type II of Hamilton-Jacobi Theorem for an $R_p$-reduced RCH System)
For the regular point reducible RCH system
$(T^*Q,G,\omega,H,F,W)$ with an $R_p$-reduced RCH system \\
$((T^\ast Q)_\mu, \omega_\mu,h_\mu, f_\mu, u_\mu )$,
assume that $\gamma: Q \rightarrow T^*Q$ is an one-form
on $Q$, and $\lambda=\gamma \cdot \pi_{Q}: T^* Q \rightarrow T^* Q, $
and for any symplectic map $\varepsilon:T^* Q \rightarrow T^* Q, $
denote by $\tilde{X}^\varepsilon = T\pi_{Q}\cdot \tilde{X} \cdot \varepsilon$,
where $\tilde{X}=X_{(T^\ast Q,G,\omega,H,F,u)}$ is the dynamical vector
field of the regular point reducible RCH system
$(T^*Q,G,\omega,H,F,W)$ with a control law $u$. Moreover,
assume that $\mu \in \mathfrak{g}^\ast $ is a regular value of the momentum
map $\mathbf{J}$, and $\textmd{Im}(\gamma)\subset
\mathbf{J}^{-1}(\mu), $ and it is $G_\mu$-invariant,
and $\varepsilon$ is $G_\mu$-invariant and $\varepsilon(\mathbf{J}^{-1}(\mu)) \subset \mathbf{J}^{-1}(\mu). $
Denote by $\bar{\gamma}=\pi_\mu(\gamma): Q \rightarrow (T^* Q)_\mu $,
$\bar{\lambda}=\pi_\mu(\lambda): \mathbf{J}^{-1}(\mu) (\subset T^*Q) \rightarrow (T^* Q)_\mu $,
and $\bar{\varepsilon}=\pi_\mu(\varepsilon): \mathbf{J}^{-1}(\mu) (\subset T^*Q) \rightarrow (T^* Q)_\mu $.
Then $\varepsilon$ and $\bar{\varepsilon}$ satisfy
the equation $T\bar{\varepsilon}\cdot(X_{h_\mu \cdot \bar{\varepsilon}})= T\bar{\lambda}\cdot \tilde{X} \cdot \varepsilon, $
if and only if they satisfy the equation
$T\bar{\gamma}\cdot \tilde{X}^\varepsilon= X_{h_\mu}\cdot \bar{\varepsilon}, $
where $X_{h_\mu}$ and $X_{h_\mu \cdot \bar{\varepsilon}} \in TT^* Q$ are the Hamiltonian
vector fields of the $R_p$-reduced Hamiltonian functions $h_\mu$ and $h_\mu \cdot \bar{\varepsilon}: T^* Q
\rightarrow \mathbb{R}, $ respectively.
The equation $T\bar{\gamma}\cdot \tilde{X}^\varepsilon=
X_{h_\mu}\cdot \bar{\varepsilon}$ is called the Type II of Hamilton-Jacobi equation for
the $R_p$-reduced RCH System $((T^\ast Q)_\mu, \omega_\mu,h_\mu, f_\mu, u_\mu)$.
\end{theo}
\noindent{\bf Proof: }
At first, we note that
$\textmd{Im}(\gamma)\subset \mathbf{J}^{-1}(\mu), $ and it
is $G_\mu$-invariant, in this case, $\pi_\mu^*\omega_\mu=
i_\mu^*\omega= \omega, $ along $\textmd{Im}(\gamma)$.
Since $\tilde{X}=X_{(T^\ast Q,G,\omega,H,F,u)}=X_H
+\textnormal{vlift}(F)+\textnormal{vlift}(u), $ and
$T\pi_{Q}\cdot \textnormal{vlift}(F)=T\pi_{Q}\cdot \textnormal{vlift}(u)=0, $
then we have that $T\pi_{Q}\cdot \tilde{X}\cdot \varepsilon=T\pi_{Q}\cdot X_H\cdot \varepsilon. $
By using the $R_p$-reduced symplectic form $\omega_\mu$, if we take
that $v= X_H\cdot \varepsilon \in TT^* Q,$ and for any $w \in TT^* Q, \; T\bar{\lambda}(w)\neq 0,$
and $T\pi_{\mu} (w) \neq 0, $ from Lemma 2.5(ii) we have that
\begin{align*}
& \omega_\mu(T\bar{\gamma} \cdot \tilde{X}^\varepsilon, \; T\pi_\mu \cdot w) =
\omega_\mu(T(\pi_\mu \cdot \gamma) \cdot T\pi_Q \cdot \tilde{X}\cdot \varepsilon, \; T\pi_\mu
\cdot w )\\
& = \pi_\mu^*\omega_\mu(T\gamma \cdot T\pi_Q \cdot X_H \cdot\varepsilon, \; w)
= \omega(T(\gamma \cdot \pi_Q)\cdot X_H\cdot \varepsilon, \; w)\\
& = \omega(X_H\cdot \varepsilon, \; w-T(\gamma \cdot \pi_Q)\cdot w)
-\mathbf{d}\gamma(T\pi_{Q}(X_H\cdot \varepsilon), \; T\pi_{Q}(w))\\
& =\omega(X_H\cdot \varepsilon, \; w) - \omega(X_H\cdot \varepsilon, \;
T\lambda\cdot w)-\mathbf{d}\gamma(T\pi_{Q}(\tilde{X}\cdot \varepsilon), \; T\pi_{Q}(w))\\
& =\pi_\mu^*\omega_\mu(X_H\cdot
\varepsilon, \; w) - \pi_\mu^*\omega_\mu(X_H\cdot \varepsilon, \; T\lambda\cdot w)+ \lambda^*\omega(\tilde{X}\cdot \varepsilon, \; w)\\
& = \omega_\mu(T\pi_\mu(X_H\cdot \varepsilon), \;
T\pi_\mu \cdot w) - \omega_\mu(T\pi_\mu\cdot(X_H\cdot \varepsilon), \;
T(\pi_\mu \cdot\lambda) \cdot w)+ \lambda^*\cdot \pi_\mu^*\cdot \omega_\mu(\tilde{X}\cdot \varepsilon, \; w)\\
& = \omega_\mu(T\pi_\mu(X_H)\cdot \pi_\mu(\varepsilon), \; T\pi_\mu \cdot w)
- \omega_\mu(T\pi_\mu(X_H)\cdot \pi_\mu(\varepsilon), \; T\bar{\lambda}\cdot w)
+ (\pi_\mu\cdot\lambda)^*\cdot \omega_\mu(\tilde{X}\cdot \varepsilon, \; w)\\
& = \omega_\mu(X_{h_\mu} \cdot
\bar{\varepsilon}, \; T\pi_\mu \cdot w)- \omega_\mu(X_{h_\mu} \cdot
\bar{\varepsilon}, \; T\bar{\lambda} \cdot w)
+ \omega_\mu(T\bar{\lambda}\cdot \tilde{X}\cdot \varepsilon, \; T\bar{\lambda}\cdot w),
\end{align*}
where we have used that $T\pi_\mu(X_H)= X_{h_\mu}. $
Note that $\varepsilon: T^* Q \rightarrow T^* Q $ is
symplectic, and $\pi_\mu^*\omega_\mu= i_\mu^*\omega
= \omega, $ along $\textmd{Im}(\gamma)$, and hence
$\bar{\varepsilon}= \pi_\mu(\varepsilon): \mathbf{J}^{-1}(\mu)(\subset T^* Q) \rightarrow (T^*
Q)_\mu $ is also symplectic along
$\textmd{Im}(\gamma)$, and $X_{h_\mu}\cdot \bar{\varepsilon}=
T\bar{\varepsilon} \cdot X_{h_\mu \cdot \bar{\varepsilon}}, $
along $\textmd{Im}(\gamma)\cap\textmd{Im}(\varepsilon)$.
From the above arguments, we can obtain that
\begin{align*}
& \omega_\mu(T\bar{\gamma} \cdot \tilde{X}^\varepsilon, \; T\pi_\mu \cdot w)-
\omega_\mu(X_{h_\mu} \cdot \bar{\varepsilon}, \; T\pi_\mu \cdot w) \\
& =- \omega_\mu(X_{h_\mu} \cdot
\bar{\varepsilon}, \; T\bar{\lambda} \cdot w)
+ \omega_\mu(T\bar{\lambda}\cdot \tilde{X}\cdot \varepsilon, \; T\bar{\lambda}\cdot w)\\
& = \omega_\mu(T\bar{\lambda}\cdot \tilde{X}\cdot \varepsilon, \; T\bar{\lambda}\cdot w)
-\omega_\mu(T\bar{\varepsilon} \cdot X_{h_\mu \cdot \bar{\varepsilon}}, \; T\bar{\lambda} \cdot w)\\
& = \omega_\mu(T\bar{\lambda}\cdot \tilde{X}\cdot \varepsilon- T\bar{\varepsilon} \cdot X_{h_\mu \cdot \bar{\varepsilon}}, \; T\bar{\lambda}\cdot w).
\end{align*}
Because the $R_p$-reduced symplectic form $\omega_\mu$ is non-degenerate,
it follows that
$T\bar{\gamma}\cdot \tilde{X}^\varepsilon=
X_{h_\mu}\cdot \bar{\varepsilon}, $ is equivalent to
$T\bar{\lambda}\cdot \tilde{X}\cdot \varepsilon
= T\bar{\varepsilon} \cdot X_{h_\mu \cdot \bar{\varepsilon}}. $
Thus, we know that the $\varepsilon$ and $\bar{\varepsilon}$ satisfy
the equation $T\bar{\varepsilon}\cdot(X_{h_\mu \cdot \bar{\varepsilon}})
= T\bar{\lambda}\cdot \tilde{X} \cdot \varepsilon, $
if and only if they satisfy the Type II of Hamilton-Jacobi equation
$T\bar{\gamma}\cdot \tilde{X}^\varepsilon= X_{h_\mu}\cdot \bar{\varepsilon}. $
\hskip 0.3cm $\blacksquare$\\

Moreover, for a given RCH system
$(T^*Q,G,\omega,H,F,W)$ with a $R_p$-reduced RCH system $((T^\ast
Q)_\mu, \omega_\mu,h_\mu, f_\mu, u_\mu )$,
we know that the Hamiltonian vector fields
$X_{H}$ and $X_{h_\mu}$ for the corresponding
Hamiltonian system $(T^*Q,G,\omega,H)$
and its $R_p$-reduced system $((T^\ast Q)_\mu,
\omega_\mu, h_\mu )$, are $\pi_\mu$-related, that is,
$X_{h_\mu}\cdot \pi_\mu=T\pi_\mu\cdot X_{H}\cdot i_\mu.$ Then we can
prove the following Theorem 4.4, which states the relationship
between the solutions of Type II of Hamilton-Jacobi equations and the
regular point reduction.
\begin{theo}
For a given RCH system
$(T^*Q,G,\omega,H,F,W)$ with an $R_p$-reduced RCH system \\
$((T^\ast Q)_\mu, \omega_\mu,h_\mu, f_\mu, u_\mu )$,
assume that $\gamma: Q \rightarrow T^*Q$ is an one-form on $Q$,
and $\lambda=\gamma \cdot \pi_{Q}: T^* Q
\rightarrow T^* Q $, and $\varepsilon: T^* Q \rightarrow T^* Q $ is a symplectic map.
Denote by $\tilde{X}^\varepsilon = T\pi_{Q}\cdot \tilde{X} \cdot \varepsilon$,
where $\tilde{X}=X_{(T^\ast Q,G,\omega,H,F,u)}$ is the dynamical vector
field of the regular point reducible RCH system
$(T^*Q,G,\omega,H,F,W)$ with a control law $u$. Moreover, assume that $\mu \in
\mathfrak{g}^\ast $ is a regular value of the momentum
map $\mathbf{J}$, and $\textmd{Im}(\gamma)\subset
\mathbf{J}^{-1}(\mu), $ and it is $G_\mu$-invariant,
and $\varepsilon$ is $G_\mu$-invariant
and $\varepsilon(\mathbf{J}^{-1}(\mu)) \subset \mathbf{J}^{-1}(\mu). $
Denote by $\bar{\gamma}=\pi_\mu(\gamma): Q \rightarrow (T^* Q)_\mu $,
$\bar{\lambda}=\pi_\mu(\lambda):
\mathbf{J}^{-1}(\mu) (\subset T^*Q) \rightarrow (T^* Q)_\mu $,
and $\bar{\varepsilon}=\pi_\mu(\varepsilon):
\mathbf{J}^{-1}(\mu) (\subset T^*Q) \rightarrow (T^* Q)_\mu $.
Then $\varepsilon$ is a solution of the Type II of Hamilton-Jacobi equation
$T\gamma\cdot \tilde{X}^\varepsilon= X_H\cdot \varepsilon, $ for the
regular point reducible RCH system $(T^*Q,G,\omega,H,F,W), $ if and only if
$\varepsilon$ and $\bar{\varepsilon} $ satisfy the Type II of Hamilton-Jacobi equation
$T\bar{\gamma}\cdot \tilde{X}^\varepsilon
= X_{h_\mu}\cdot \bar{\varepsilon}, $ for the
$R_p$-reduced RCH system $((T^\ast Q)_\mu,
\omega_\mu,h_\mu, f_\mu, u_\mu)$.
\end{theo}
\noindent{\bf Proof: }
Note that $\textmd{Im}(\gamma)\subset \mathbf{J}^{-1}(\mu), $ and it
is $G_\mu$-invariant, in this case, $\pi_\mu^*\omega_\mu=
i_\mu^*\omega= \omega, $ along $\textmd{Im}(\gamma)$.
Since the Hamiltonian vector fields
$X_{H}$ and $X_{h_\mu}$ are $\pi_\mu$-related, that is,
$X_{h_\mu}\cdot \pi_\mu= T\pi_\mu\cdot X_{H}\cdot i_\mu, $ and
by using the $R_p$-reduced symplectic form $\omega_\mu$,
for any $w \in TT^* Q,$ and $T\pi_{\mu} \cdot w
\neq 0, $ we have that
\begin{align*}
& \omega_\mu(T\bar{\gamma} \cdot \tilde{X}^\varepsilon
- X_{h_\mu} \cdot \bar{\varepsilon}, \; T\pi_\mu \cdot w) \\
& = \omega_\mu(T\bar{\gamma} \cdot \tilde{X}^\varepsilon, \; T\pi_\mu \cdot w)-
\omega_\mu(X_{h_\mu} \cdot \bar{\varepsilon}, \; T\pi_\mu \cdot w) \\
& = \omega_\mu(T\pi_\mu \cdot T\gamma \cdot \tilde{X}^\varepsilon, \; T\pi_\mu \cdot w)-
\omega_\mu(X_{h_\mu} \cdot \pi_\mu \cdot \varepsilon, \; T\pi_\mu \cdot w) \\
& = \pi_\mu^*\omega_\mu(T\gamma \cdot \tilde{X}^\varepsilon, \; w)
-\omega_\mu(T\pi_\mu\cdot X_{H}\cdot \varepsilon, \; T\pi_\mu \cdot w) \\
& = \pi_\mu^*\omega_\mu(T\gamma \cdot \tilde{X}^\varepsilon, \; w)
-\pi_\mu^*\omega_\mu(X_{H}\cdot \varepsilon, \; w)\\
& = \omega(T\gamma \cdot \tilde{X}^\varepsilon- X_{H}\cdot \varepsilon, \; w).
\end{align*}
Because both the symplectic form $\omega$ and
the $R_p$-reduced symplectic form $\omega_\mu$ are non-degenerate,
it follows that the equation
$T\bar{\gamma}\cdot \tilde{X}^\varepsilon= X_{h_\mu}\cdot \bar{\varepsilon}, $
is equivalent to the equation
$T\gamma\cdot \tilde{X}^\varepsilon= X_H\cdot \varepsilon$. Thus,
$\varepsilon$ is a solution of the Type II of Hamilton-Jacobi equation
$T\gamma\cdot \tilde{X}^\varepsilon= X_H\cdot \varepsilon, $ for the
regular point reducible RCH system $(T^*Q,G,\omega,H,F,W), $ if and only if
$\varepsilon$ and $\bar{\varepsilon} $ satisfy the Type II of Hamilton-Jacobi equation
$T\bar{\gamma}\cdot \tilde{X}^\varepsilon= X_{h_\mu}\cdot \bar{\varepsilon}, $ for the
$R_p$-reduced RCH system $((T^\ast Q)_\mu,
\omega_\mu,h_\mu, f_\mu, u_\mu)$.  \hskip 0.3cm
$\blacksquare$
\begin{rema}
In particular, if both the external force and control of a regular
point reducible RCH system $(T^*Q,G,\omega,H,F,u)$ are zero, in this
case the RCH system is just a regular point reducible Hamiltonian
system $(T^*Q,G,\omega,H)$. From the proofs of
the above Theorem 4.2, 4.3 and 4.4, we can also get
two types of Hamilton-Jacobi theorem
for the associated Marsden-Weinstein reduced Hamiltonian system, which is given in
Wang \cite{wa17}. Thus, Theorem 4.2, 4.3 and 4.4 can be regarded as an extension
of two types of Hamilton-Jacobi theorem for the $R_p$-reduced Hamiltonian
system to that for the $R_p$-reduced RCH system.
\end{rema}
Moreover, for the regular point reducible RCH system we can also
introduce the regular point reducible controlled Hamiltonian
equivalence (RpCH-equivalence) as follows.
\begin{defi}(RpCH-equivalence)
Suppose that we have two regular point reducible RCH systems
$(T^\ast Q_i, G_i,\omega_i,H_i, F_i, W_i),\; i=1,2$, we say them to
be RpCH-equivalent, or simply,\\ $(T^\ast Q_1,
G_1,\omega_1,H_1,F_1,W_1)\stackrel{RpCH}{\sim}(T^\ast
Q_2,G_2,\omega_2,H_2,F_2,W_2)$, if there exists a diffeomorphism
$\varphi:Q_1\rightarrow Q_2$ such that the following controlled Hamiltonian
matching conditions hold:\\
\noindent {\bf RpCH-1:} The cotangent lift map $\varphi^\ast:T^\ast
Q_2\rightarrow T^\ast Q_1$ is symplectic.\\
\noindent {\bf RpCH-2:} For $\mu_i\in \mathfrak{g}^\ast_i $, the
regular reducible points of RCH systems $(T^\ast Q_i, G_i,\omega_i,
H_i, F_i, W_i),\\ i=1,2$, the map
$\varphi_\mu^\ast=i_{\mu_1}^{-1}\cdot\varphi^\ast\cdot i_{\mu_2}:
\mathbf{J}_2^{-1}(\mu_2)\rightarrow \mathbf{J}_1^{-1}(\mu_1)$ is
$(G_{2\mu_2},G_{1\mu_1})$-equivariant and $W_1=\varphi_\mu^\ast
(W_2)$, where $\mu=(\mu_1, \mu_2)$, and denote by
$i_{\mu_1}^{-1}(S)$ the preimage of a subset $S\subset T^\ast Q_1$
for the map $i_{\mu_1}:\mathbf{J}_1^{-1}(\mu_1)\rightarrow T^\ast
Q_1$.\\
\noindent {\bf RpCH-3:}
$Im[X_{H_1}+ \textnormal{vlift}(F_1)-
T\varphi^\ast (X_{H_2})-\textnormal{vlift}(\varphi^\ast
F_2\varphi_\ast)]\subset\textnormal{vlift}(W_1)$.
\end{defi}

Then we can obtain the following regular point reduction theorem for an RCH system,
which explains the relationship between the RpCH-equivalence for the
regular point reducible RCH systems with symmetries and the
RCH-equivalence for the associated $R_p$-reduced RCH systems, its proof
is given in Marsden et al \cite{mawazh10}. This theorem can be
regarded as an extension of regular point reduction theorem of
Hamiltonian systems under regular controlled Hamiltonian equivalence
conditions.
\begin{theo}
Two regular point reducible RCH systems $(T^\ast Q_i, G_i, \omega_i,
H_i, F_i,W_i)$, $i=1,2,$ are RpCH-equivalent if and only if the
associated $R_p$-reduced RCH systems $((T^\ast
Q_i)_{\mu_i},\omega_{i\mu_i},h_{i\mu_i},f_{i\mu_i},\\ W_{i\mu_i}),$
$i=1,2,$ are RCH-equivalent.
\end{theo}

Moreover, if considering the RpCH-equivalence of the regular point
reducible RCH systems and using the above
Theorem 4.7, Theorem 4.4 and Theorem 2.11,
we can obtain the following Theorem 4.8,
which states that the solutions of two types of Hamilton-Jacobi equations for
the regular point reducible RCH systems leave invariant
under the conditions of RpCH-equivalence.
\begin{theo}
Suppose that two regular point reducible RCH systems
$(T^\ast Q_i,G_i,\omega_i,H_i,F_i,W_i)$,
$i=1,2,$ are RpCH-equivalent with an equivalent map $\varphi: Q_1
\rightarrow Q_2 $, and the associated $R_p$-reduced RCH systems
$((T^\ast Q_i)_{\mu_i},\omega_{i\mu_i},h_{i\mu_i},f_{i\mu_i}, u_{i\mu_i}),$
$i=1,2.$ Under the hypotheses and notations of Theorem 4.2,
Theorem 4.3 and Theorem 4.4, then we have that\\

\noindent $(\mathrm{i})$ If the one-form
$\gamma_2: Q_2 \rightarrow T^* Q_2$ is closed with
respect to $T\pi_{Q_2}: TT^* Q_2 \rightarrow TQ_2, $ and
$\bar{\gamma}_2=\pi_{\mu_2}(\gamma_2): Q_2 \rightarrow (T^* Q_2)_{\mu_2} $
is a solution of the Type I of Hamilton-Jacobi equation for
the $R_p$-reduced RCH System $((T^\ast Q_2)_{\mu_2},
\omega_{2\mu_2},h_{2\mu_2}, f_{2\mu_2}, u_{2\mu_2})$.
Then $\gamma_1=\varphi^* \cdot \gamma_2\cdot \varphi: Q_1 \rightarrow T^* Q_1, $ is
a solution of the Type I of Hamilton-Jacobi equation for the
RCH System $(T^*Q_1,G_1,\omega_1,H_1,F_1,W_1) $,
and $\bar{\gamma}_1=\pi_{\mu_1}(\gamma_1): Q_1 \rightarrow (T^* Q_1)_{\mu_1}$ is
a solution of the Type I of Hamilton-Jacobi equation for the $R_p$-reduced
RCH System $((T^\ast Q_1)_{\mu_1}, \omega_{1\mu_1},
h_{1\mu_1}, f_{1\mu_1}, u_{1\mu_1})$.
Vice versa.\\

\noindent $(\mathrm{ii})$ If the $G_{2\mu_2}$-invariant
symplectic map $\varepsilon_2: T^*Q_2\rightarrow T^* Q_2$ and
$\bar{\varepsilon}_2=\pi_{\mu_2}(\varepsilon_2):
\mathbf{J_2}^{-1}(\mu_2) (\subset T^*Q_2) \rightarrow (T^* Q_2)_{\mu_2} $
satisfy the Type II of Hamilton-Jacobi equation for the $R_p$-reduced RCH System
$((T^\ast Q_2)_{\mu_2}, \omega_{2\mu_2},h_{2\mu_2}, f_{2\mu_2}, u_{2\mu_2})$.
Then $\varepsilon_1=
\varphi^* \cdot \varepsilon_2\cdot \varphi_*: T^*Q_1 \rightarrow T^* Q_1 $ and
$\bar{\varepsilon}_1=\pi_{\mu_1}(\varepsilon_1):
\mathbf{J_1}^{-1}(\mu_1) (\subset T^*Q_1) \rightarrow (T^* Q_1)_{\mu_1} $
satisfy the Type II of Hamilton-Jacobi equation for the $R_p$-reduced
RCH System $((T^\ast Q_1)_{\mu_1}, \omega_{1\mu_1},h_{1\mu_1}, f_{1\mu_1}, u_{1\mu_1}). $
Vice versa.
\end{theo}
\noindent{\bf Proof: }
We first prove the assertion $(\mathrm{i})$.
If two regular point reducible RCH systems $(T^\ast Q_i,G_i,\omega_i,\\ H_i,F_i,W_i)$,
$i=1,2,$ are RpCH-equivalent, from Definition 4.6 we know that the two RCH
systems are also RCH-equivalent. If the one-form
$\gamma_2: Q_2 \rightarrow T^* Q_2$ is closed with
respect to $T\pi_{Q_2}: TT^* Q_2 \rightarrow TQ_2, $ from Theorem 2.6, we know that
$\gamma_2$ is a solution of the Type I of Hamilton-Jacobi equation
for the RCH system $(T^*Q_2,G_2,\omega_2,H_2,F_2,W_2). $
Moreover, from Theorem 2.11, we know that
$\gamma_1=\varphi^* \cdot \gamma_2\cdot \varphi: Q_1 \rightarrow T^* Q_1, $ is
a solution of the Type I of Hamilton-Jacobi equation for the
RCH System $(T^*Q_1,G_1,\omega_1,H_1,F_1,W_1) $.
On the other hand,
if $(T^\ast Q_1, G_1, \omega_1, H_1, F_1,
W_1)\stackrel{RpCH}{\sim}(T^\ast Q_2, G_2, \omega_2, H_2,F_2,
W_2)$, then from Definition 4.6 there exists a diffeomorphism
$\varphi: Q_1 \rightarrow Q_2, $ such that $\varphi^\ast: T^\ast Q_2
\rightarrow T^\ast Q_1$ is symplectic, and for $\mu_i\in
\mathfrak{g}^\ast_i, \; i=1,2, $ $\varphi_\mu^\ast=
i_{\mu_1}^{-1}\cdot \varphi^\ast \cdot
i_{\mu_2}:\mathbf{J}_2^{-1}(\mu_2)\rightarrow
\mathbf{J}_1^{-1}(\mu_1)$ is $(G_{2\mu_2},G_{1\mu_1})$-equivariant.
From the following commutative Diagram-9:
\[
\begin{CD}
Q_2 @> \gamma_2 >> T^\ast Q_2 @<i_{\mu_2}
<< \mathbf{J}_2^{-1}(\mu_2) @>\pi_{\mu_2}>> (T^\ast Q_2)_{\mu_2}\\
@A \varphi AA @V\varphi^\ast VV @V\varphi^\ast_\mu
VV @V\varphi^\ast_{\mu/G}VV\\
Q_1 @> \gamma_1 >> T^\ast Q_1 @<i_{\mu_1}
<< \mathbf{J}_1^{-1}(\mu_1)
@>\pi_{\mu_1}>>(T^\ast Q_1)_{\mu_1}
\end{CD}
\]
$$\mbox{Diagram-9}$$
we have a well-defined symplectic map $\varphi_{\mu/G}^\ast:(T^\ast
Q_2)_{\mu_2}\rightarrow (T^\ast Q_1)_{\mu_1}, $ such that
$\varphi_{\mu/G}^\ast \cdot
\pi_{\mu_2}=\pi_{\mu_1}\cdot\varphi^\ast_\mu$, see Marsden et al
\cite{mawazh10}. Then from the regular point reduction
Theorem 4.7 we know that the associated
$R_p$-reduced RCH systems $((T^\ast
Q_i)_{\mu_i},\omega_{i\mu_i},h_{i\mu_i},f_{i\mu_i}, W_{i\mu_i}),$
$i=1,2,$ are RCH-equivalent with an equivalent map
$\varphi_{\mu/G}^\ast:(T^\ast Q_2)_{\mu_2}\rightarrow (T^\ast
Q_1)_{\mu_1}. $ If $\bar{\gamma}_2=\pi_{\mu_2}(\gamma_2): Q_2 \rightarrow (T^*
Q_2)_{\mu_2}$ is a solution of the Type I of Hamilton-Jacobi equation for
$R_p$-reduced RCH system $((T^\ast Q_2)_{\mu_2},
\omega_{2\mu_2},h_{2\mu_2},f_{2\mu_2},u_{2\mu_2})$.
and note that $\bar{\gamma}_1=\pi_{\mu_1}(\gamma_1)=
\pi_{\mu_1}\cdot\varphi^\ast (\gamma_2)\cdot \varphi=
\varphi_{\mu/G}^\ast \cdot \pi_{\mu_2}(\gamma_2)\cdot \varphi=
\varphi_{\mu/G}^\ast \cdot \bar{\gamma}_2 \cdot \varphi $,
from Theorem 2.11, we know that
$\bar{\gamma}_1=\pi_{\mu_1}(\gamma_1): Q_1 \rightarrow (T^* Q_1)_{\mu_1}$ is
a solution of the Type I of Hamilton-Jacobi equation for the $R_p$-reduced
RCH System $((T^\ast Q_1)_{\mu_1}, \omega_{1\mu_1},h_{1\mu_1}, f_{1\mu_1}, u_{1\mu_1})$.
Because the map $\varphi: Q_1 \rightarrow Q_2 $ is a diffeomorphism,
and $\varphi^\ast: T^\ast Q_2\rightarrow T^\ast Q_1$ is a symplectic isomorphisms,
vice versa. It follows that the
assertion $(\mathrm{i})$ of Theorem 4.8 holds.\\

Next, we prove the assertion $(\mathrm{ii})$.
If the $G_{2\mu_2}$-invariant symplectic map
$\varepsilon_2: T^*Q_2\rightarrow T^* Q_2$ and
$\bar{\varepsilon}_2=\pi_{\mu_2}(\varepsilon_2):
\mathbf{J_2}^{-1}(\mu_2) (\subset T^*Q_2) \rightarrow (T^* Q_2)_{\mu_2} $
satisfy the Type II of Hamilton-Jacobi equation for the $R_p$-reduced RCH System
$((T^\ast Q_2)_{\mu_2}, \omega_{2\mu_2},h_{2\mu_2}, f_{2\mu_2}, u_{2\mu_2})$,
from Theorem 4.4 we know that
$\varepsilon_2$ is a solution of the Type II of Hamilton-Jacobi equation
for the RCH system $(T^*Q_2,G_2,\omega_2,H_2,F_2,W_2). $ Since
the two regular point reducible RCH systems $(T^\ast Q_i,G_i,\omega_i,H_i,F_i,W_i)$,
$i=1,2,$ are RpCH-equivalent, and hence are also RCH-equivalent,
and from Theorem 2.11 we know that
$\varepsilon_1= \varphi^* \cdot \varepsilon_2\cdot \varphi_*: T^*Q_1 \rightarrow T^* Q_1 $
is a solution of the Type II of Hamilton-Jacobi equation for the
RCH System $(T^*Q_1,G_1,\omega_1,H_1,F_1,W_1) $.
Moreover, from Theorem 4.4 we know that
$\varepsilon_1$ and $\bar{\varepsilon}_1=\pi_{\mu_1}(\varepsilon_1) $
satisfy the Type II of Hamilton-Jacobi equation for the $R_p$-reduced
RCH System $((T^\ast Q_1)_{\mu_1}, \omega_{1\mu_1},
h_{1\mu_1}, f_{1\mu_1}, u_{1\mu_1}). $
In the same way, because the map $\varphi: Q_1 \rightarrow Q_2 $ is a diffeomorphism,
and $\varphi^\ast: T^\ast Q_2\rightarrow T^\ast Q_1$ is a symplectic isomorphisms,
vice versa. We prove the assertion $(\mathrm{ii})$ of Theorem 4.8. \hskip 0.3cm
$\blacksquare$

\begin{rema}
If $(T^\ast Q, \omega)$ is a connected symplectic manifold, and
$\mathbf{J}:T^\ast Q\rightarrow \mathfrak{g}^\ast$ is a
non-equivariant momentum map with a non-equivariance group
one-cocycle $\sigma: G\rightarrow \mathfrak{g}^\ast$, which is
defined by $\sigma(g):=\mathbf{J}(g\cdot
z)-\operatorname{Ad}^\ast_{g^{-1}}\mathbf{J}(z)$, where $g\in G$ and
$z\in T^\ast Q$. Then we know that $\sigma$ produces a new affine
action $\Theta: G\times \mathfrak{g}^\ast \rightarrow
\mathfrak{g}^\ast $ defined by
$\Theta(g,\mu):=\operatorname{Ad}^\ast_{g^{-1}}\mu + \sigma(g)$,
where $\mu \in \mathfrak{g}^\ast$, with respect to which the given
momentum map $\mathbf{J}$ is equivariant. Assume that $G$ acts
freely and properly on $T^\ast Q$, and $\tilde{G}_\mu$ denotes the
isotropy subgroup of $\mu \in \mathfrak{g}^\ast$ relative to this
affine action $\Theta$ and $\mu$ is a regular value of $\mathbf{J}$.
Then the quotient space $(T^\ast
Q)_\mu=\mathbf{J}^{-1}(\mu)/\tilde{G}_\mu$ is also a symplectic
manifold with the $R_p$-reduced symplectic form
$\omega_\mu$ uniquely characterized by
$(4.1)$, see Ortega and Ratiu \cite{orra04} and Marsden et al
\cite{mamiorpera07}. In this case, we can also define the regular
point reducible RCH system $(T^*Q,G,\omega,H,F,W)$ and
RpCH-equivalence, and prove the two types of Hamilton-Jacobi theorem for
the $R_p$-reduced RCH system $((T^\ast Q)_\mu, \omega_\mu,h_\mu,
f_\mu,u_\mu )$ by using the above similar way, and state that the
solutions of two types of Hamilton-Jacobi equations
for the regular point reducible RCH
systems with symmetries leave invariant under the conditions of
RpCH-equivalence, where the $R_p$-reduced space $((T^\ast Q)_\mu,
\omega_\mu )$ is determined by the affine action.
\end{rema}

\section{Hamilton-Jacobi Theorem of the $R_o$-reduced RCH System }

We know that the orbit reduction for a Hamiltonian system
is an alternative approach to symplectic reduction
given by Marle \cite{ma76}
and Kazhdan, Kostant and Sternberg \cite{kakost78}, which is different from the
Marsden-Weinstein reduction. Thus, the regular orbit reduction for an RCH system
is different from the regular point reduction.
In this section, we first give the regular orbit reducible
RCH system with symmetry and momentum map. Then we
give precisely the geometric constraint conditions of the $R_o$-reduced symplectic form for the
dynamical vector field of the regular orbit reducible RCH system,
and prove the Type I and Type II of
Hamilton-Jacobi theorem for the $R_o$-reduced RCH system. ¡¡
Moreover, we state the relationship between the solutions of
Type II of Hamilton-Jacobi equations and regular orbit reduction. Finally,
we consider the RoCH-equivalence, and state that the solutions of two types of
Hamilton-Jacobi equations for the RCH systems with
symmetries leave invariant under the conditions of RoCH-equivalence.
We shall follow the notations and conventions introduced in Marsden
et al \cite{mawazh10}, Wang \cite{wa17}, Wang \cite{wa18}.\\

At first, we consider the regular orbit reducible RCH system,
which is given by Marsden et al \cite{mawazh10}.
Assume that the cotangent lifted
left action $\Phi^{T^\ast}:G\times T^\ast Q\rightarrow T^\ast Q$ is
symplectic, free and proper, and admits an
$\operatorname{Ad}^\ast$-equivariant momentum map $\mathbf{J}:T^\ast
Q\rightarrow \mathfrak{g}^\ast$. Let $\mu\in \mathfrak{g}^\ast$ be a
regular value of the momentum map $\mathbf{J}$ and
$\mathcal{O}_\mu=G\cdot \mu\subset \mathfrak{g}^\ast$ be the
$G$-orbit of the coadjoint $G$-action through the point $\mu$. Since
$G$ acts freely, properly and symplectically on $T^\ast Q$, then the
quotient space $(T^\ast Q)_{\mathcal{O}_\mu}=
\mathbf{J}^{-1}(\mathcal{O}_\mu)/G$ is a regular quotient symplectic
manifold with the $R_o$-reduced symplectic form $\omega_{\mathcal{O}_\mu}$
uniquely characterized by the relation
\begin{equation}i_{\mathcal{O}_\mu}^\ast \omega=\pi_{\mathcal{O}_{\mu}}^\ast
\omega_{\mathcal{O}
_\mu}+\mathbf{J}_{\mathcal{O}_\mu}^\ast\omega_{\mathcal{O}_\mu}^+,
\label{5.1}
\end{equation} where $\mathbf{J}_{\mathcal{O}_\mu}$ is
the restriction of the momentum map $\mathbf{J}$ to
$\mathbf{J}^{-1}(\mathcal{O}_\mu)$, that is,
$\mathbf{J}_{\mathcal{O}_\mu}=\mathbf{J}\cdot i_{\mathcal{O}_\mu}$
and $\omega_{\mathcal{O}_\mu}^+$ is the $+$-symplectic structure on
the orbit $\mathcal{O}_\mu$ given by
\begin{equation}\omega_{\mathcal{O}_\mu}^
+(\nu)(\xi_{\mathfrak{g}^\ast}(\nu),\eta_{\mathfrak{g}^\ast}(\nu))
=<\nu,[\xi,\eta]>,\;\; \forall\;\nu\in\mathcal{O}_\mu, \;
\xi,\eta\in \mathfrak{g}. \label{5.2}
\end{equation} The maps
$i_{\mathcal{O}_\mu}:\mathbf{J}^{-1}(\mathcal{O}_\mu)\rightarrow
T^\ast Q$ and
$\pi_{\mathcal{O}_\mu}:\mathbf{J}^{-1}(\mathcal{O}_\mu)\rightarrow
(T^\ast Q)_{\mathcal{O}_\mu}$ are natural injection and the
projection, respectively. The pair $((T^\ast
Q)_{\mathcal{O}_\mu},\omega_{\mathcal{O}_\mu})$ is called
the $R_o$-reduced symplectic space of $(T^\ast Q,\omega)$ at $\mu$.
In general case, we
maybe thought that the structure of the $R_o$-reduced symplectic
space $((T^\ast Q)_{\mathcal{O}_\mu},\omega_{\mathcal{O}_\mu})$ is
more complex than that of the $R_p$-reduced symplectic space
$((T^\ast Q)_\mu,\omega_\mu)$, but, from the regular reduction
diagram, see Ortega and Ratiu \cite{orra04},
we know that the $R_o$-reduced space $((T^\ast
Q)_{\mathcal{O}_\mu},\omega_{\mathcal{O}_\mu})$ is symplectic
diffeomorphic to the $R_p$-reduced space $((T^*Q)_\mu,
\omega_\mu)$, and hence is also symplectic diffeomorphic to a
symplectic fiber bundle.\\

Let $H:T^\ast Q\rightarrow \mathbb{R}$ be a $G$-invariant
Hamiltonian, the flow $F_t$ of the Hamiltonian vector field $X_H$
leaves the connected components of
$\mathbf{J}^{-1}(\mathcal{O}_\mu)$ invariant and commutes with the
$G$-action, so it induces a flow $f_t^{\mathcal{O}_\mu}$ on $(T^\ast
Q)_{\mathcal{O}_\mu}$, defined by $f_t^{\mathcal{O}_\mu}\cdot
\pi_{\mathcal{O}_\mu}=\pi_{\mathcal{O}_\mu} \cdot F_t\cdot
i_{\mathcal{O}_\mu}$, and the vector field $X_{h_{\mathcal{O}_\mu}}$
generated by the flow $f_t^{\mathcal{O}_\mu}$ on $((T^\ast
Q)_{\mathcal{O}_\mu},\omega_{\mathcal{O}_\mu})$ is Hamiltonian with
the associated $R_o$-reduced Hamiltonian function
$h_{\mathcal{O}_\mu}:(T^\ast Q)_{\mathcal{O}_\mu}\rightarrow
\mathbb{R}$ defined by $h_{\mathcal{O}_\mu}\cdot
\pi_{\mathcal{O}_\mu}= H\cdot i_{\mathcal{O}_\mu}$, and the
Hamiltonian vector fields $X_H$ and $X_{h_{\mathcal{O}_\mu}}$ are
$\pi_{\mathcal{O}_\mu}$-related.
Moreover, assume that the fiber-preserving map $F:T^\ast Q\rightarrow T^\ast
Q$ and the control subset $W$ of\; $T^\ast Q$ are both $G$-invariant.
In order to get the $R_o$-reduced RCH system, we also assume that
$F(\mathbf{J}^{-1}(\mathcal{O}_\mu))\subset
\mathbf{J}^{-1}(\mathcal{O}_\mu)$, and $W \cap
\mathbf{J}^{-1}(\mathcal{O}_\mu)\neq \emptyset $.
Thus, we can introduce a regular orbit
reducible RCH systems as follows, see Marsden et al \cite{mawazh10}
and Wang \cite{wa18}.

\begin{defi}
(Regular Orbit Reducible RCH System) A 6-tuple $(T^\ast Q, G,
\omega,H,F,W)$ with
the canonical symplectic form $\omega$ on $T^*Q$, where the Hamiltonian $H: T^\ast Q\rightarrow
\mathbb{R}$, the fiber-preserving map $F: T^\ast Q\rightarrow T^\ast
Q$ and the fiber submanifold $W$ of $T^\ast Q$ are all
$G$-invariant, is called a regular orbit reducible RCH system, if
there exists an orbit $\mathcal{O}_\mu, \; \mu\in\mathfrak{g}^\ast$,
where $\mu$ is a regular value of the momentum map $\mathbf{J}$,
such that the regular orbit reduced system, that is, the 5-tuple
$((T^\ast
Q)_{\mathcal{O}_\mu},\omega_{\mathcal{O}_\mu},
h_{\mathcal{O}_\mu},f_{\mathcal{O}_\mu},
W_{\mathcal{O}_\mu})$, where $(T^\ast
Q)_{\mathcal{O}_\mu}=\mathbf{J}^{-1}(\mathcal{O}_\mu)/G$,
$\pi_{\mathcal{O}_\mu}^\ast \omega_{\mathcal{O}_\mu}
=i_{\mathcal{O}_\mu}^\ast\omega
-\mathbf{J}_{\mathcal{O}_\mu}^\ast\omega_{\mathcal{O}_\mu}^+$,
$h_{\mathcal{O}_\mu}\cdot \pi_{\mathcal{O}_\mu} =H\cdot
i_{\mathcal{O}_\mu}$, $F(\mathbf{J}^{-1}(\mathcal{O}_\mu))\subset
\mathbf{J}^{-1}(\mathcal{O}_\mu)$, $f_{\mathcal{O}_\mu}\cdot
\pi_{\mathcal{O}_\mu}=\pi_{\mathcal{O}_\mu}\cdot F\cdot
i_{\mathcal{O}_\mu}$, and $W \cap
\mathbf{J}^{-1}(\mathcal{O}_\mu)\neq \emptyset $,
$W_{\mathcal{O}_\mu}=\pi_{\mathcal{O}_\mu}(W \cap
\mathbf{J}^{-1}(\mathcal{O}_\mu))$, is an RCH system,
which is simply written as $R_o$-reduced RCH system. Where $((T^\ast
Q)_{\mathcal{O}_\mu},\omega_{\mathcal{O}_\mu})$ is the $R_o$-reduced
space, the function $h_{\mathcal{O}_\mu}:(T^\ast
Q)_{\mathcal{O}_\mu}\rightarrow \mathbb{R}$ is called the $R_o$-reduced
Hamiltonian, the fiber-preserving map $f_{\mathcal{O}_\mu}:(T^\ast
Q)_{\mathcal{O}_\mu} \rightarrow (T^\ast Q)_{\mathcal{O}_\mu}$ is
called the $R_o$-reduced (external) force map, $W_{\mathcal{O}_\mu}$ is a
fiber submanifold of $(T^\ast Q)_{\mathcal{O}_\mu}$, and is called
the $R_o$-reduced control subset.
\end{defi}

Denote by $X_{(T^\ast Q,G,\omega,H,F,u)}$ the dynamical vector field of the
regular orbit reducible RCH system $(T^\ast Q, G,\omega, H,F,W)$
with a control law $u$. Assume that it can be expressed by
\begin{equation}X_{(T^\ast Q,G,\omega,H,F,u)}
=X_H+\textnormal{vlift}(F)+\textnormal{vlift}(u).\label{5.3}
\end{equation}
If an $R_o$-reduced feedback control law
$u_{\mathcal{O}_\mu}:(T^\ast Q)_{\mathcal{O}_\mu}\rightarrow
W_{\mathcal{O}_\mu}$ is chosen, the $R_o$-reduced RCH system
$((T^\ast Q)_{\mathcal{O}_\mu},
\omega_{\mathcal{O}_\mu},h_{\mathcal{O}_\mu},f_{\mathcal{O}_\mu},u_{\mathcal{O}_\mu})$
is a closed-loop regular dynamical system with a control law
$u_{\mathcal{O}_\mu}$. Assume that its dynamical vector field $X_{((T^\ast
Q)_{\mathcal{O}_\mu}, \omega_{\mathcal{O}_\mu},
h_{\mathcal{O}_\mu},f_{\mathcal{O}_\mu},u_{\mathcal{O}_\mu})}$ can
be expressed by
\begin{equation}X_{((T^\ast Q)_{\mathcal{O}_\mu},
\omega_{\mathcal{O}_\mu},h_{\mathcal{O}_\mu},f_{\mathcal{O}_\mu},u_{\mathcal{O}_\mu})}=
X_{h_{\mathcal{O}_\mu}}+\textnormal{vlift}(f_{\mathcal{O}_\mu})
+\textnormal{vlift}(u_{\mathcal{O}_\mu}),
\label{5.4}\end{equation}
where $X_{h_{\mathcal{O}_\mu}}$ is the Hamiltonian vector field of
the $R_o$-reduced Hamiltonian $h_{\mathcal{O}_\mu}$, $\textnormal{vlift}(f_{\mathcal{O}_\mu})=
\textnormal{vlift}(f_{\mathcal{O}_\mu})X_{h_{\mathcal{O}_\mu}}$,
$\textnormal{vlift}(u_{\mathcal{O}_\mu})=
\textnormal{vlift}(u_{\mathcal{O}_\mu})X_{h_{\mathcal{O}_\mu}}$, and
satisfies the condition
\begin{equation}X_{((T^\ast Q)_{\mathcal{O}_\mu},
\omega_{\mathcal{O}_\mu},h_{\mathcal{O}_\mu},f_{\mathcal{O}_\mu},u_{\mathcal{O}_\mu})}\cdot
\pi_{\mathcal{O}_\mu} =T\pi_{\mathcal{O}_\mu} \cdot X_{(T^\ast
Q,G,\omega,H,F,u)}\cdot
i_{\mathcal{O}_\mu}.\label{5.5}
\end{equation}

Since the set of Hamiltonian systems with symmetries on the cotangent
bundle is not complete under the Marsden-Weinstein reduction,
then it is not yet complete under the regular orbit reduction,
and the $R_o$-reduced system of an
RCH system with symmetry defined on the cotangent bundle
$T^*Q$ may not be an RCH system on a cotangent bundle.
On the other hand, from the expression of the dynamical vector
field of an $R_o$-reduced RCH system, we know that under the
actions of the $R_o$-reduced external force $f_{\mathcal{O}_\mu}$
and the $R_o$-reduced control $u_{\mathcal{O}_\mu}$, in general, the dynamical vector
field is not Hamiltonian, and the $R_o$-reduced RCH system is not
yet a Hamiltonian system. Thus, we can not describe the Hamilton-Jacobi equation for
the $R_o$-reduced RCH system from the viewpoint of generating
function just like same as Theorem 1.1
given by Abraham and Marsden in \cite{abma78}.
But, for a given regular orbit reducible RCH system
$(T^*Q,G,\omega,H,F,W)$ with an $R_o$-reduced RCH system $((T^\ast
Q)_{\mathcal{O}_\mu}, \omega_{\mathcal{O}_\mu},h_{\mathcal{O}_\mu},
f_{\mathcal{O}_\mu},u_{\mathcal{O}_\mu})$,
by using Lemma 2.5, we can give precisely the following geometric constraint
conditions of the $R_o$-reduced symplectic form for the
dynamical vector field of the regular orbit reducible RCH system,
that is, Type I and Type II of
Hamilton-Jacobi equation for the $R_o$-reduced RCH system
$((T^\ast Q)_{\mathcal{O}_\mu}, \omega_{\mathcal{O}_\mu},h_{\mathcal{O}_\mu},
f_{\mathcal{O}_\mu},u_{\mathcal{O}_\mu})$.
At first, by using Lemma 2.5 and the $R_o$-reduced symplectic
form $\omega_{\mathcal{O}_\mu}$,
and the fact that the one-form $\gamma: Q
\rightarrow T^*Q $ is closed with respect to
$T\pi_Q: TT^* Q \rightarrow TQ, $ and a stronger
assumption condition $\textmd{Im}(\gamma)\subset
\mathbf{J}^{-1}(\mu), $ and it is $G$-invariant,
we can also prove the following Type I of Hamilton-Jacobi
theorem for the $R_o$-reduced RCH system $((T^\ast
Q)_{\mathcal{O}_\mu}, \omega_{\mathcal{O}_\mu},h_{\mathcal{O}_\mu},
f_{\mathcal{O}_\mu},u_{\mathcal{O}_\mu})$.

\begin{theo} (Type I of Hamilton-Jacobi Theorem for an $R_o$-reduced RCH system)
For a given regular orbit reducible RCH system
$(T^*Q,G,\omega,H,F,W)$ with an $R_o$-reduced RCH system \\
$((T^\ast Q)_{\mathcal{O}_\mu}, \omega_{\mathcal{O}_\mu},h_{\mathcal{O}_\mu},
f_{\mathcal{O}_\mu},u_{\mathcal{O}_\mu})$,
assume that $\gamma: Q \rightarrow T^*Q$ is an
one-form on $Q$, and $\tilde{X}^\gamma = T\pi_{Q}\cdot \tilde{X} \cdot \gamma$,
where $\tilde{X}=X_{(T^\ast Q,G,\omega,H,F,u)}$ is the dynamical vector
field of the regular orbit reducible RCH system
$(T^*Q,G,\omega,H,F,W)$ with a control law $u$.
Moreover, assume that $\mu \in \mathfrak{g}^\ast $ is a regular
value of the momentum map $\mathbf{J}$, and $\mathcal{O}_\mu, \;
(\mu \in \mathfrak{g}^\ast), $ is the regular reducible orbit of the
corresponding Hamiltonian system $(T^*Q,G,\omega,H)$,
and $\textmd{Im}(\gamma)\subset
\mathbf{J}^{-1}(\mu), $ and it is $G$-invariant,
$\bar{\gamma}=\pi_{\mathcal{O}_\mu}(\gamma): Q \rightarrow (T^*
Q)_{\mathcal{O}_\mu}. $ If the one-form $\gamma: Q \rightarrow T^*Q
$ is closed with respect to $T\pi_Q: TT^* Q \rightarrow TQ, $ then
$\bar{\gamma}$ is a solution of the equation $T\bar{\gamma}\cdot
\tilde{X}^\gamma= X_{h_{\mathcal{O}_\mu}}\cdot \bar{\gamma}, $
where $X_{h_{\mathcal{O}_\mu}}$
is the Hamiltonian vector field of the $R_o$-reduced
Hamiltonian function $h_{\mathcal{O}_\mu}:(T^\ast
Q)_{\mathcal{O}_\mu}\rightarrow \mathbb{R}, $
and the equation is
called the Type I of Hamilton-Jacobi equation for
the $R_o$-reduced RCH system $((T^\ast Q)_{\mathcal{O}_\mu},
\omega_{\mathcal{O}_\mu},h_{\mathcal{O}_\mu},
f_{\mathcal{O}_\mu},u_{\mathcal{O}_\mu})$. Here the maps involved in
the theorem are shown in the following Diagram-10.
\begin{center}
\hskip 0cm \xymatrix{ \mathbf{J}^{-1}(\mathcal{O}_\mu)
\ar[r]^{i_{\mathcal{O}_\mu}} & T^* Q \ar[d]_{X_H} \ar[r]^{\pi_Q}
& Q \ar[d]_{\tilde{X}^\gamma} \ar[r]^{\gamma} & T^*Q \ar[d]_{\tilde{X}} \ar[r]^{\pi_{\mathcal{O}_\mu}}
& (T^* Q)_{\mathcal{O}_\mu} \ar[d]_{X_{h_{\mathcal{O}_\mu}}} \\
  & T(T^*Q)  & TQ \ar[l]^{T\gamma} & T(T^*Q) \ar[l]^{T\pi_Q} \ar[r]_{T\pi_{\mathcal{O}_\mu}}
  & T(T^* Q)_{\mathcal{O}_\mu} }
\end{center}
$$\mbox{Diagram-10}$$
\end{theo}

\noindent{\bf Proof: } At first, from Theorem 2.6, we know that
$\gamma$ is a solution of the Type I of Hamilton-Jacobi equation
$T\gamma\cdot \tilde{X}^\gamma= X_H\cdot \gamma .$
Next, we note that the $R_o$-reduced symplectic space
$(T^\ast Q)_{\mathcal{O}_\mu}= \mathbf{J}^{-1}(\mathcal{O}_\mu)/G
\cong \mathbf{J}^{-1}(\mu)/G \times \mathcal{O}_\mu, $ with the $R_o$-reduced
symplectic form $\omega_{\mathcal{O}_\mu}$ uniquely characterized by
the relation $i_{\mathcal{O}_\mu}^\ast
\omega=\pi_{\mathcal{O}_{\mu}}^\ast \omega_{\mathcal{O}
_\mu}+\mathbf{J}_{\mathcal{O}_\mu}^\ast\omega_{\mathcal{O}_\mu}^+. $
Since $\textmd{Im}(\gamma)\subset \mathbf{J}^{-1}(\mu), $ and it is
$G$-invariant, in this case for any $V\in TQ, $ and $w\in TT^*Q, $
we have that
$\mathbf{J}_{\mathcal{O}_\mu}^\ast\omega_{\mathcal{O}_\mu}^+(T\gamma
\cdot V, \; w)=0, $ and hence
$\pi_{\mathcal{O}_\mu}^*\omega_{\mathcal{O}_\mu}=
i_{\mathcal{O}_\mu}^*\omega= \omega, $ along $\textmd{Im}(\gamma)$.
Note that $\tilde{X}=X_{(T^\ast Q,G,\omega,H,F,u)}=X_H
+\textnormal{vlift}(F)+\textnormal{vlift}(u), $ and
$T\pi_{Q}\cdot \textnormal{vlift}(F)=T\pi_{Q}\cdot \textnormal{vlift}(u)=0, $
then we have that $T\pi_{Q}\cdot \tilde{X}\cdot \gamma
=T\pi_{Q}\cdot X_H\cdot \gamma. $
By using the $R_o$-reduced symplectic form
$\omega_{\mathcal{O}_\mu}$, and if
we take that $v= X_H\cdot \gamma \in TT^* Q,$ and
for any $w \in TT^* Q, \; T\pi_{Q}(w)\neq 0,$
and $T\pi_{\mathcal{O}_\mu} (w) \neq 0,
$ from Lemma 2.5(ii) we have that
\begin{align*}
& \omega_{\mathcal{O}_\mu}(T\bar{\gamma} \cdot \tilde{X}^\gamma, \;
T\pi_{\mathcal{O}_\mu} \cdot w) =
\omega_{\mathcal{O}_\mu}(T(\pi_{\mathcal{O}_\mu} \cdot \gamma) \cdot
T\pi_Q \cdot \tilde{X}\cdot \gamma, \; T\pi_{\mathcal{O}_\mu} \cdot w )\\
& = \pi_{\mathcal{O}_\mu}^*\omega_{\mathcal{O}_\mu}(T\gamma \cdot
T\pi_Q \cdot X_H \cdot\gamma, \; w)
= \omega(T(\gamma \cdot \pi_Q)\cdot X_H\cdot \gamma, \; w)\\
& = \omega(X_H\cdot \gamma, \; w-T(\gamma \cdot
\pi_Q)\cdot w)-\mathbf{d}\gamma(T\pi_{Q}(X_H\cdot \gamma), \; T\pi_{Q}(w))\\
& = \omega_{\mathcal{O}_\mu}(X_{h_{\mathcal{O}_\mu}} \cdot
\bar{\gamma}, \; T\pi_{\mathcal{O}_\mu} \cdot w)-
\omega_{\mathcal{O}_\mu}(X_{h_{\mathcal{O}_\mu}} \cdot
\bar{\gamma}, \; T\bar{\gamma} \cdot T\pi_{Q}(w))
-\mathbf{d}\gamma(T\pi_{Q}(X_H\cdot \gamma), \; T\pi_{Q}(w)),
\end{align*}
in which we have used that $ T\pi_{\mathcal{O}_\mu} \cdot X_H= X_{
h_{\mathcal{O}_\mu}}. $
Since the one-form $\gamma: Q \rightarrow T^*Q $ is closed with respect to
$T\pi_Q: TT^* Q \rightarrow TQ, $ then we have that $
\mathbf{d}\gamma(T\pi_{Q}(X_H\cdot \gamma), \; T\pi_{Q}(w))=0, $
and hence
\begin{equation}
 \omega_{\mathcal{O}_\mu}(T\bar{\gamma} \cdot \tilde{X}^\gamma, \;
 T\pi_{\mathcal{O}_\mu} \cdot w)-
\omega_{\mathcal{O}_\mu}(X_{h_{\mathcal{O}_\mu}} \cdot \bar{\gamma}, \;
T\pi_{\mathcal{O}_\mu} \cdot w)
 = - \omega_{\mathcal{O}_\mu}(X_{h_{\mathcal{O}_\mu}} \cdot
\bar{\gamma}, \; T\bar{\gamma} \cdot T\pi_{Q}(w)). \; \label{5.6}
\end{equation}
If $\bar{\gamma}$ satisfies the equation
$T\bar{\gamma}\cdot \tilde{X}^\gamma
= X_{h_{\mathcal{O}_\mu}}\cdot \bar{\gamma}, $
from Lemma 2.5(i) we can obtain that
\begin{align*}
- \omega_{\mathcal{O}_\mu}(X_{h_{\mathcal{O}_\mu}} \cdot
\bar{\gamma}, \; T\bar{\gamma} \cdot T\pi_{Q}(w))
& = -\omega_{\mathcal{O}_\mu} (T\bar{\gamma}
\cdot \tilde{X}^\gamma, \; T\bar{\gamma} \cdot T\pi_{Q}(w))\\
& = -\bar{\gamma}^*\omega_{\mathcal{O}_\mu}
(T\pi_{Q} \cdot \tilde{X}\cdot\gamma, \; T\pi_Q(w))\\
& = -\gamma^* \cdot \pi_{\mathcal{O}_\mu}^*
\omega_{\mathcal{O}_\mu}(T\pi_{Q} \cdot X_{H}\cdot\gamma, \; T\pi_{Q}(w))\\
& = -\gamma^*\omega( T\pi_{Q}(X_{H}\cdot\gamma), \; T\pi_{Q}(w))\\
& = \textbf{d}\gamma(T\pi_{Q}( X_{H}\cdot\gamma ), \; T\pi_{Q}(w))=0.
\end{align*}
Because the $R_o$-reduced symplectic form $\omega_{\mathcal{O}_\mu}$
is non-degenerate, the left side of (5.6) equals zero, only when
$\bar{\gamma}$ satisfies the equation
$T\bar{\gamma}\cdot \tilde{X}^\gamma
= X_{h_{\mathcal{O}_\mu}}\cdot \bar{\gamma}.$
Thus, if the one-form $\gamma: Q \rightarrow T^*Q $ is closed with respect to
$T\pi_Q: TT^* Q \rightarrow TQ, $ then $\bar{\gamma}$ must be
a solution of the Type I of Hamilton-Jacobi equation
$T\bar{\gamma}\cdot \tilde{X}^\gamma
= X_{h_{\mathcal{O}_\mu}}\cdot \bar{\gamma}.$
\hskip 0.3cm $\blacksquare$\\

Next, for any $G$-invariant symplectic map
$\varepsilon: T^* Q \rightarrow T^* Q $,
we can also prove the following Type II of
Hamilton-Jacobi theorem for the $R_o$-reduced
RCH system $((T^\ast Q)_{\mathcal{O}_\mu}, \omega_{\mathcal{O}_\mu},h_{\mathcal{O}_\mu},
f_{\mathcal{O}_\mu},u_{\mathcal{O}_\mu})$.

\begin{theo} (Type II of Hamilton-Jacobi Theorem for an $R_o$-reduced RCH system)
For a given regular orbit reducible RCH system
$(T^*Q,G,\omega,H,F,W)$ with an $R_o$-reduced RCH system \\
$((T^\ast Q)_{\mathcal{O}_\mu}, \omega_{\mathcal{O}_\mu},h_{\mathcal{O}_\mu},
f_{\mathcal{O}_\mu},u_{\mathcal{O}_\mu})$,
assume that $\gamma: Q \rightarrow T^*Q$ is an one-form on $Q$,
and $\lambda=\gamma \cdot \pi_{Q}: T^* Q
\rightarrow T^* Q $, and for any $G$-invariant symplectic map
$\varepsilon: T^* Q \rightarrow T^* Q $,
denote by $\tilde{X}^\varepsilon = T\pi_{Q}\cdot
\tilde{X} \cdot \varepsilon$, where
$\tilde{X}=X_{(T^\ast Q,G,\omega,H,F,u)}$ is the dynamical vector
field of the regular orbit reducible RCH system
$(T^*Q,G,\omega,H,F,W)$ with a control law $u$.
Moreover, assume that $\mu \in \mathfrak{g}^\ast $ is a regular value of the momentum
map $\mathbf{J}$, and $\mathcal{O}_\mu, \; (\mu \in \mathfrak{g}^\ast), $ is the
regular reducible orbit of the corresponding Hamiltonian system $(T^*Q,G,\omega,H)$, and
$\textmd{Im}(\gamma)\subset \mathbf{J}^{-1}(\mu), $ and it is
$G$-invariant, and $\varepsilon(\mathbf{J}^{-1}(\mathcal{O}_\mu))
\subset \mathbf{J}^{-1}(\mathcal{O}_\mu). $
Denote by $\bar{\gamma}=\pi_{\mathcal{O}_\mu}(\gamma): Q
\rightarrow (T^* Q)_{\mathcal{O}_\mu} $,
$\bar{\lambda}=\pi_{\mathcal{O}_\mu}(\lambda):
\mathbf{J}^{-1}(\mathcal{O}_\mu) (\subset T^*Q) \rightarrow (T^*
Q)_{\mathcal{O}_\mu} $, and $\bar{\varepsilon}=\pi_{\mathcal{O}_\mu}(\varepsilon):
\mathbf{J}^{-1}(\mathcal{O}_\mu) (\subset T^*Q) \rightarrow (T^*
Q)_{\mathcal{O}_\mu} $. Then $\varepsilon$ and $\bar{\varepsilon}$ satisfy the equation
$T\bar{\varepsilon}\cdot X_{h_{\mathcal{O}_\mu}\cdot\bar{\varepsilon}}
= T\bar{\lambda} \cdot \tilde{X}\cdot\varepsilon, $
if and only if they satisfy the equation $T\bar{\gamma}\cdot \tilde{X}^\varepsilon=
X_{h_{\mathcal{O}_\mu}}\cdot \bar{\varepsilon}$,
where $X_{h_{\mathcal{O}_\mu}}$ and
$ X_{h_{\mathcal{O}_\mu} \cdot \bar{\varepsilon}} \in TT^*Q $
are the Hamiltonian vector fields of the $R_o$-reduced Hamiltonian
functions $h_{\mathcal{O}_\mu}$ and
$h_{\mathcal{O}_\mu} \cdot \bar{\varepsilon}: T^*Q\rightarrow
\mathbb{R}, $ respectively.
The equation $T\bar{\gamma}\cdot \tilde{X}^\varepsilon=
X_{h_{\mathcal{O}_\mu}}\cdot \bar{\varepsilon}$
is called the Type II of Hamilton-Jacobi
equation for the $R_o$-reduced RCH system $((T^\ast
Q)_{\mathcal{O}_\mu}, \omega_{\mathcal{O}_\mu},h_{\mathcal{O}_\mu},
f_{\mathcal{O}_\mu},u_{\mathcal{O}_\mu})$.
Here the maps involved in the theorem are shown in the following Diagram-11.

\begin{center}
\hskip 0cm \xymatrix{ \mathbf{J}^{-1}(\mathcal{O}_\mu)
\ar[r]^{i_{\mathcal{O}_\mu}} & T^* Q \ar[d]_{X_{H\cdot \varepsilon}}
\ar[dr]^{\tilde{X}^\varepsilon} \ar[r]^{\pi_Q}
& Q \ar[r]^{\gamma} & T^*Q \ar[d]_{\tilde{X}}
\ar[dr]^{X_{h_{\mathcal{O}_\mu} \cdot\bar{\varepsilon}}} \ar[r]^{\pi_{\mathcal{O}_\mu}}
& (T^* Q)_{\mathcal{O}_\mu} \ar[d]^{X_{h_{\mathcal{O}_\mu}}} \\
  & T(T^*Q)  & TQ \ar[l]^{T\gamma} & T(T^*Q) \ar[l]^{T\pi_Q}
   \ar[r]_{T\pi_{\mathcal{O}_\mu}} & T(T^* Q)_{\mathcal{O}_\mu} }
\end{center}
$$\mbox{Diagram-11}$$
\end{theo}

\noindent{\bf Proof: } At first, we note that the $R_o$-reduced symplectic space
$(T^\ast Q)_{\mathcal{O}_\mu}= \mathbf{J}^{-1}(\mathcal{O}_\mu)/G
\cong (\mathbf{J}^{-1}(\mu)/G) \times \mathcal{O}_\mu, $ with the $R_o$-reduced
symplectic form $\omega_{\mathcal{O}_\mu}$ uniquely characterized by
the relation $i_{\mathcal{O}_\mu}^\ast
\omega=\pi_{\mathcal{O}_{\mu}}^\ast \omega_{\mathcal{O}
_\mu}+\mathbf{J}_{\mathcal{O}_\mu}^\ast\omega_{\mathcal{O}_\mu}^+. $
Since $\textmd{Im}(\gamma)\subset \mathbf{J}^{-1}(\mu), $ and it is
$G$-invariant, in this case for any $V\in TQ, $ and $w\in TT^*Q, $
we have that
$\mathbf{J}_{\mathcal{O}_\mu}^\ast\omega_{\mathcal{O}_\mu}^+(T\gamma
\cdot V, \; w)=0, $ and hence
$\pi_{\mathcal{O}_\mu}^*\omega_{\mathcal{O}_\mu}=
i_{\mathcal{O}_\mu}^*\omega= \omega, $ along $\textmd{Im}(\gamma)$.
Note that $\tilde{X}=X_{(T^\ast Q,G,\omega,H,F,u)}=X_H
+\textnormal{vlift}(F)+\textnormal{vlift}(u), $ and
$T\pi_{Q}\cdot \textnormal{vlift}(F)=T\pi_{Q}\cdot \textnormal{vlift}(u)=0, $
then we have that $T\pi_{Q}\cdot \tilde{X}\cdot \varepsilon
=T\pi_{Q}\cdot X_H\cdot \varepsilon. $
By using the reduced symplectic form $\omega_{\mathcal{O}_\mu}$, and if
we take that $v= X_H\cdot \varepsilon \in TT^* Q, $ and
for any $w \in TT^* Q, \; T\bar{\lambda}(w)\neq 0, $
and $T\pi_{\mathcal{O}_\mu} (w) \neq 0,
$ from Lemma 2.5 we have that
\begin{align*}
& \omega_{\mathcal{O}_\mu}(T\bar{\gamma} \cdot \tilde{X}^\varepsilon, \;
T\pi_{\mathcal{O}_\mu} \cdot w) =
\omega_{\mathcal{O}_\mu}(T(\pi_{\mathcal{O}_\mu} \cdot \gamma) \cdot
T\pi_{Q}\cdot \tilde{X}\cdot \varepsilon, \; T\pi_{\mathcal{O}_\mu} \cdot w )\\
& = \pi_{\mathcal{O}_\mu}^*\omega_{\mathcal{O}_\mu}(T\gamma \cdot
T\pi_{Q}\cdot X_H\cdot \varepsilon, \; w)
= \omega(T(\gamma \cdot \pi_Q)\cdot X_H\cdot \varepsilon, \; w)\\
& = \omega(X_H\cdot \varepsilon, \; w-T(\gamma \cdot
\pi_Q)\cdot w)-\mathbf{d}\gamma(T\pi_{Q}(X_H\cdot \varepsilon), \; T\pi_{Q}(w))\\
& =\omega(X_H\cdot \varepsilon, \; w) - \omega(X_H\cdot \varepsilon, \;
T\lambda\cdot w)-\mathbf{d}\gamma(T\pi_{Q}(\tilde{X}\cdot \varepsilon), \; T\pi_{Q}(w))\\
& =\omega(X_H\cdot \varepsilon, \; w) - \omega(X_H\cdot \varepsilon, \;
T\lambda\cdot w)+ \lambda^*\omega(\tilde{X}\cdot \varepsilon, \; w)\\
& = \omega_{\mathcal{O}_\mu}(X_{h_{\mathcal{O}_\mu}} \cdot
\bar{\varepsilon}, \; T\pi_{\mathcal{O}_\mu} \cdot w)-
\omega_{\mathcal{O}_\mu}(X_{h_{\mathcal{O}_\mu}} \cdot
\bar{\varepsilon}, \; T\bar{\lambda} \cdot w)
+ \omega_{\mathcal{O}_\mu}(T\bar{\lambda}
\cdot \tilde{X}\cdot \varepsilon, \; T\bar{\lambda} \cdot w),
\end{align*}
in which we have used that $ T\pi_{\mathcal{O}_\mu} \cdot X_H= X_{
h_{\mathcal{O}_\mu}}. $
Note that $\varepsilon: T^* Q
\rightarrow T^* Q $ is symplectic, and
$\pi_{\mathcal{O}_\mu}^*\omega_{\mathcal{O}_\mu}=
i_{\mathcal{O}_\mu}^*\omega= \omega, $ along $\textmd{Im}(\gamma)$,
and hence $\bar{\varepsilon}=
\pi_{\mathcal{O}_\mu}\cdot \varepsilon:
\mathbf{J}^{-1}(\mathcal{O}_\mu) (\subset T^* Q)
\rightarrow (T^* Q)_{\mathcal{O}_\mu}$ is also symplectic along
$\textmd{Im}(\gamma)$, and $X_{h_{\mathcal{O}_\mu}}\cdot
\bar{\varepsilon} = T\bar{\varepsilon} \cdot X_{h_{\mathcal{O}_\mu} \cdot
\bar{\varepsilon}}, $ along
$\textmd{Im}(\gamma)\cap\textmd{Im}(\varepsilon). $
From the above arguments, we can obtain that
\begin{align*}
& \omega_{\mathcal{O}_\mu}(T\bar{\gamma} \cdot X_H^\varepsilon, \;
T\pi_{\mathcal{O}_\mu} \cdot w)-
\omega_{\mathcal{O}_\mu}(X_{h_{\mathcal{O}_\mu}} \cdot
\bar{\varepsilon}, \; T\pi_{\mathcal{O}_\mu} \cdot w)\\
& = -\omega_{\mathcal{O}_\mu}(X_{h_{\mathcal{O}_\mu}} \cdot
\bar{\varepsilon}, \; T\bar{\lambda} \cdot w)
+ \omega_{\mathcal{O}_\mu}(T\bar{\lambda}
\cdot \tilde{X}\cdot \varepsilon, \; T\bar{\lambda} \cdot w)\\
& = \omega_{\mathcal{O}_\mu}(T\bar{\lambda}
\cdot \tilde{X}\cdot \varepsilon, \; T\bar{\lambda}\cdot w)
-\omega_{\mathcal{O}_\mu}(T\bar{\varepsilon}
\cdot X_{h_{\mathcal{O}_\mu} \cdot \bar{\varepsilon}}, \; T\bar{\lambda} \cdot w)\\
& = \omega_{\mathcal{O}_\mu}(T\bar{\lambda}\cdot \tilde{X}\cdot \varepsilon
- T\bar{\varepsilon} \cdot X_{h_{\mathcal{O}_\mu}
\cdot \bar{\varepsilon}}, \; T\bar{\lambda}\cdot w).
\end{align*}
Because the $R_o$-reduced symplectic form
$\omega_{\mathcal{O}_\mu}$ is non-degenerate, it follows that
$T\bar{\gamma}\cdot \tilde{X}^\varepsilon=
X_{h_{\mathcal{O}_\mu}}\cdot \bar{\varepsilon}$ is equivalent to
$T\bar{\varepsilon}\cdot X_{h_{\mathcal{O}_\mu}\cdot\bar{\varepsilon}}
= T\bar{\lambda} \cdot \tilde{X}\cdot\varepsilon. $
Thus, we know that the $\varepsilon$ and $\bar{\varepsilon}$ satisfy the equation
$T\bar{\varepsilon}\cdot X_{h_{\mathcal{O}_\mu}\cdot\bar{\varepsilon}}
= T\bar{\lambda} \cdot \tilde{X}\cdot\varepsilon, $
if and only if they satisfy the Type II of Hamilton-Jacobi equation
$T\bar{\gamma}\cdot \tilde{X}^\varepsilon=
X_{h_{\mathcal{O}_\mu}}\cdot \bar{\varepsilon}. $
\hskip 0.3cm $\blacksquare$\\

It is worthy of noting that the different symplectic forms determine
the different regular reduced RCH systems. From (5.1) we know that,
for the regular orbit reduced symplectic space
$(T^\ast Q)_{\mathcal{O}_\mu}= \mathbf{J}^{-1}(\mathcal{O}_\mu)/G
\cong \mathbf{J}^{-1}(\mu)/G \times \mathcal{O}_\mu, $ if we give a stronger
assumption condition, that is, for the one-form $\gamma: Q \rightarrow T^*Q$ on $Q,$
assume that $\textmd{Im}(\gamma)\subset
\mathbf{J}^{-1}(\mu), $ (note that it is not $\textmd{Im}(\gamma)\subset
\mathbf{J}^{-1}(\mathcal{O}_\mu) $),
and it is $G$-invariant, then in this case for any $V\in TQ, $ and $w\in TT^*Q, $
we have that
$\mathbf{J}_{\mathcal{O}_\mu}^\ast\omega_{\mathcal{O}_\mu}^+(T\gamma
\cdot V, \; w)=0, $ and hence from (5.1), $i_{\mathcal{O}_\mu}^\ast
\omega=\pi_{\mathcal{O}_{\mu}}^\ast \omega_{\mathcal{O}
_\mu}+\mathbf{J}_{\mathcal{O}_\mu}^\ast\omega_{\mathcal{O}_\mu}^+, $
we have that $\pi_{\mathcal{O}_\mu}^*\omega_{\mathcal{O}_\mu}=
i_{\mathcal{O}_\mu}^*\omega= \omega, $ along $\textmd{Im}(\gamma)$.
Thus, we can use Lemma 2.5 for the regular orbit reduced symplectic form $\omega_{\mathcal{O}_\mu}$
in the proofs of Theorem 5.2 and Theorem 5.3. It is easy to give the wrong results
without the precise analysis for the regular orbit reduction case.\\

Moreover, for the regular orbit reducible RCH system
$(T^*Q,G,\omega,H,F,W)$ with an $R_o$-reduced RCH system $((T^\ast
Q)_{\mathcal{O}_\mu}, \omega_{\mathcal{O}_\mu},h_{\mathcal{O}_\mu},f_{\mathcal{O}_\mu},u_{\mathcal{O}_\mu})$,
we know that the Hamiltonian vector fields
$X_{H}$ and $X_{h_{\mathcal{O}_\mu}}$ for
the corresponding Hamiltonian system $(T^*Q,G,\omega,H)$
and its $R_o$-reduced system $((T^\ast Q)_{\mathcal{O}_\mu},
\omega_{\mathcal{O}_\mu},h_{\mathcal{O}_\mu})$, are
$\pi_{\mathcal{O}_\mu}$-related, that is, $X_{
h_{\mathcal{O}_\mu}}\cdot \pi_{\mathcal{O}_\mu} =
T\pi_{\mathcal{O}_\mu}\cdot X_{H}\cdot i_{\mathcal{O}_\mu}. $ Then
we can also prove the following Theorem 5.4, which states the
relationship between the solutions of Type II of Hamilton-Jacobi
equations and the regular orbit reduction.

\begin{theo}
For the regular orbit reducible RCH system
$(T^*Q,G,\omega,H,F,W)$ with an $R_o$-reduced RCH system $((T^\ast
Q)_{\mathcal{O}_\mu}, \omega_{\mathcal{O}_\mu},
h_{\mathcal{O}_\mu},f_{\mathcal{O}_\mu},u_{\mathcal{O}_\mu})$,
assume that $\gamma: Q \rightarrow T^*Q$ is an
one-form on $Q$, and $\lambda=\gamma \cdot \pi_{Q}: T^* Q
\rightarrow T^* Q $, and $\varepsilon: T^* Q \rightarrow T^* Q $ is
a $G$-invariant symplectic map. Denote
$\tilde{X}^\varepsilon = T\pi_{Q}\cdot
\tilde{X} \cdot \varepsilon$, where
$\tilde{X}=X_{(T^\ast Q,G,\omega,H,F,u)}$ is the dynamical vector
field of the regular orbit reducible RCH system
$(T^*Q,G,\omega,H,F,W)$ with a control law $u$. Moreover, assume that $\mu \in
\mathfrak{g}^\ast $ is a regular value of the momentum map
$\mathbf{J}$, and $\mathcal{O}_\mu, \; (\mu \in \mathfrak{g}^\ast),
$ is the regular reducible orbit of the corresponding Hamiltonian system, and
$\textmd{Im}(\gamma)\subset \mathbf{J}^{-1}(\mu), $ and it is
$G$-invariant, and $\varepsilon(\mathbf{J}^{-1}(\mathcal{O}_\mu))
\subset \mathbf{J}^{-1}(\mathcal{O}_\mu). $ Denote by
$\bar{\gamma}=\pi_{\mathcal{O}_\mu}(\gamma): Q \rightarrow (T^*
Q)_{\mathcal{O}_\mu} $,
$\bar{\lambda}=\pi_{\mathcal{O}_\mu}(\lambda):
\mathbf{J}^{-1}(\mathcal{O}_\mu) (\subset T^*Q) \rightarrow (T^*
Q)_{\mathcal{O}_\mu} $, and
$\bar{\varepsilon}=\pi_{\mathcal{O}_\mu}(\varepsilon):
\mathbf{J}^{-1}(\mathcal{O}_\mu) (\subset T^*Q) \rightarrow (T^*
Q)_{\mathcal{O}_\mu} $. Then $\varepsilon$ is a solution of the Type
II of Hamilton-Jacobi equation $T\gamma\cdot \tilde{X}^\varepsilon=
X_H\cdot \varepsilon, $ for the regular orbit reducible RCH
system $(T^*Q,G,\omega,H,F,W), $ if and only if $\varepsilon$ and
$\bar{\varepsilon}$ satisfy the Type II of Hamilton-Jacobi equation
$T\bar{\gamma}\cdot \tilde{X}^\varepsilon=
X_{h_{\mathcal{O}_\mu}\cdot\bar{\varepsilon}}, $ for the
$R_o$-reduced RCH system $((T^\ast Q)_{\mathcal{O}_\mu},
\omega_{\mathcal{O}_\mu},h_{\mathcal{O}_\mu},
f_{\mathcal{O}_\mu},u_{\mathcal{O}_\mu})$.
\end{theo}

\noindent{\bf Proof: } Note that the $R_o$-reduced symplectic space
$(T^\ast Q)_{\mathcal{O}_\mu}= \mathbf{J}^{-1}(\mathcal{O}_\mu)/G
\cong (\mathbf{J}^{-1}(\mu)/G) \times \mathcal{O}_\mu, $ with the $R_o$-reduced
symplectic form $\omega_{\mathcal{O}_\mu}$ uniquely characterized by
the relation $i_{\mathcal{O}_\mu}^\ast
\omega=\pi_{\mathcal{O}_{\mu}}^\ast \omega_{\mathcal{O}
_\mu}+\mathbf{J}_{\mathcal{O}_\mu}^\ast\omega_{\mathcal{O}_\mu}^+. $
Since $\textmd{Im}(\gamma)\subset \mathbf{J}^{-1}(\mu), $ and it is
$G$-invariant, in this case for any $V\in TQ, $ and $w\in TT^*Q, $
we have that
$\mathbf{J}_{\mathcal{O}_\mu}^\ast\omega_{\mathcal{O}_\mu}^+(T\gamma
\cdot V, \; w)=0, $ and hence
$\pi_{\mathcal{O}_\mu}^*\omega_{\mathcal{O}_\mu}=
i_{\mathcal{O}_\mu}^*\omega= \omega, $ along $\textmd{Im}(\gamma)$.
Since the Hamiltonian vector fields
$X_{H}$ and $X_{h_{\mathcal{O}_\mu}}$ are
$\pi_{\mathcal{O}_\mu}$-related, that is,
$X_{h_{\mathcal{O}_\mu}}\cdot \pi_{\mathcal{O}_\mu}
= T\pi_{\mathcal{O}_\mu}\cdot X_{H}\cdot i_{\mathcal{O}_\mu}, $
and by using the $R_o$-reduced symplectic form
$\omega_{\mathcal{O}_\mu}$, for any $w \in TT^* Q,$ and
$T\pi_{\mathcal{O}_\mu} \cdot w \neq 0, $ we have that
\begin{align*}
& \omega_{\mathcal{O}_\mu}(T\bar{\gamma}\cdot \tilde{X}^\varepsilon
-X_{h_{\mathcal{O}_\mu}}\cdot
\bar{\varepsilon}, \; T\pi_{\mathcal{O}_\mu} \cdot w)\\
& = \omega_{\mathcal{O}_\mu}(T\bar{\gamma} \cdot \tilde{X}^\varepsilon, \;
T\pi_{\mathcal{O}_\mu} \cdot w)-
\omega_{\mathcal{O}_\mu}(X_{h_{\mathcal{O}_\mu}} \cdot
\bar{\varepsilon}, \; T\pi_{\mathcal{O}_\mu} \cdot w)\\
& = \omega_{\mathcal{O}_\mu}(T\pi_{\mathcal{O}_\mu}
\cdot T\gamma \cdot \tilde{X}^\varepsilon, \; T\pi_{\mathcal{O}_\mu} \cdot w)
-\omega_{\mathcal{O}_\mu}(X_{h_{\mathcal{O}_\mu}}
\cdot \pi_{\mathcal{O}_\mu} \cdot \varepsilon, \; T\pi_{\mathcal{O}_\mu} \cdot w)\\
& = \pi_{\mathcal{O}_\mu}^*\omega_{\mathcal{O}_\mu}(T\gamma\cdot \tilde{X}^\varepsilon, \; w)
- \omega_{\mathcal{O}_\mu}(T\pi_{\mathcal{O}_\mu}
\cdot X_{H}\cdot \varepsilon, \; T\pi_{\mathcal{O}_\mu} \cdot w)\\
& = \pi_{\mathcal{O}_\mu}^*\omega_{\mathcal{O}_\mu}
(T\gamma\cdot \tilde{X}^\varepsilon, \; w)
-\pi_{\mathcal{O}_\mu}^*\omega_{\mathcal{O}_\mu}(X_{H}\cdot \varepsilon, \; w)\\
& = \omega(T\gamma\cdot \tilde{X}^\varepsilon, \; w)
-\omega(X_{H}\cdot \varepsilon, \; w)\\
& = \omega(T\gamma\cdot \tilde{X}^\varepsilon- X_{H}\cdot \varepsilon, \; w).
\end{align*}
Because the symplectic form $\omega$ and the $R_o$-reduced symplectic
form $\omega_{\mathcal{O}_\mu}$ are non-degenerate,
it follows that the equation
$T\bar{\gamma}\cdot \tilde{X}^\varepsilon
= X_{h_{\mathcal{O}_\mu}}\cdot \bar{\varepsilon}, $
is equivalent to the equation
$T\gamma\cdot \tilde{X}^\varepsilon= X_H\cdot \varepsilon. $ Thus,
$\varepsilon$ is a solution of the Type II of Hamilton-Jacobi equation
$T\gamma\cdot \tilde{X}^\varepsilon= X_H\cdot \varepsilon, $ for the
regular orbit reducible RCH system $(T^*Q,G,\omega,H,F,W), $ if and only if
$\varepsilon$ and $\bar{\varepsilon} $ satisfy the Type II of Hamilton-Jacobi equation
$T\bar{\gamma}\cdot \tilde{X}^\varepsilon
= X_{h_{\mathcal{O}_\mu}}\cdot \bar{\varepsilon}, $ for the
$R_o$-reduced RCH system $((T^\ast Q)_{\mathcal{O}_\mu},
\omega_{\mathcal{O}_\mu}, h_{\mathcal{O}_\mu},
f_{\mathcal{O}_\mu},u_{\mathcal{O}_\mu}). $
\hskip 0.3cm $\blacksquare$

\begin{rema}
In particular, if both the external force and control of a regular
orbit reducible RCH system $(T^*Q,G,\omega,H,F,u)$ are zero, in this
case the RCH system is just a regular orbit reducible Hamiltonian
system $(T^*Q,G,\omega,H)$. From the proofs of
the above Theorem 5.2, 5.3 and 5.4, we can also get the two type of Hamilton-Jacobi theorem
for the associated $R_o$-reduced Hamiltonian system, which is given in Wang
\cite{wa17}. Thus, Theorem 5.2, 5.3 and 5.4 can be regarded as an extension of
two types of Hamilton-Jacobi theorem for the $R_o$-reduced Hamiltonian
system to that for the $R_o$-reduced RCH system.
\end{rema}

Moreover, for the regular orbit reducible RCH system we can also
introduce the regular orbit reducible controlled Hamiltonian
equivalence (RoCH-equivalence) as follows.

\begin{defi}
(RoCH-equivalence) Suppose that we have two regular orbit reducible
RCH systems $(T^\ast Q_i, G_i, \omega_i, H_i, F_i, W_i)$, $i=1,2$,
we say them to be RoCH-equivalent, or simply,\\ $(T^\ast Q_1, G_1,
\omega_1, H_1, F_1, W_1)\stackrel{RoCH}{\sim}(T^\ast Q_2, G_2,
\omega_2, H_2, F_2, W_2)$, if there exists a diffeomorphism
$\varphi:Q_1\rightarrow Q_2$ such that the following controlled Hamiltonian
matching conditions hold:

\noindent {\bf RoCH-1:} The cotangent lift map $\varphi^\ast: T^\ast
Q_2\rightarrow T^\ast Q_1$ is symplectic.

\noindent {\bf RoCH-2:} For $\mathcal{O}_{\mu_i},\; \mu_i\in
\mathfrak{g}^\ast_i$, the regular reducible orbits of RCH systems
$(T^\ast Q_i, G_i, \omega_i, H_i, F_i,\\ W_i)$, $i=1,2$, the map
$\varphi^\ast_{\mathcal{O}_\mu}=i_{\mathcal{O}_{\mu_1}}^{-1}\cdot\varphi^\ast\cdot
i_{\mathcal{O}_{\mu_2}}:\mathbf{J}_2^{-1}(\mathcal{O}_{\mu_2})\rightarrow
\mathbf{J}_1^{-1}(\mathcal{O}_{\mu_1})$ is $(G_2,G_1)$-equivariant,
$W_1=\varphi_{\mathcal{O}_\mu}^\ast (W_2)$, and
$\mathbf{J}_{2\mathcal{O}_{\mu_2}}^\ast
\omega_{2\mathcal{O}_{\mu_2}}^{+}=(\varphi_{\mathcal{O}_\mu}^\ast)^\ast
\cdot\mathbf{J}_{1\mathcal{O}_{\mu_1}}^\ast\omega_{1\mathcal{O}_{\mu_1}}^{+},$
where $\mu=(\mu_1, \mu_2)$, and denote by
$i_{\mathcal{O}_{\mu_1}}^{-1}(S)$ the preimage of a subset $S\subset
T^\ast Q_1$ for the map
$i_{\mathcal{O}_{\mu_1}}:\mathbf{J}_1^{-1}(\mathcal{O}_{\mu_1})\rightarrow
T^\ast Q_1$.

\noindent {\bf RoCH-3:}
$Im[X_{H_1}+\textnormal{vlift}(F_1)-T\varphi^\ast (X_{H_2})
-\textnormal{vlift}(\varphi^\ast F_2\varphi_\ast)]\subset
\textnormal{vlift}(W_1).$
\end{defi}

It is worthy of noting that for the regular orbit reducible RCH
system, the induced equivalent map $\varphi^*$ not only keeps the
symplectic structure and the restriction of the $+$-symplectic
structure on the regular orbit to
$\mathbf{J}^{-1}(\mathcal{O}_\mu)$, but also keeps the equivariance
of $G$-action on the regular orbit.\\

Moreover, we can obtain the
following regular orbit reduction theorem for an RCH system, which
explains the relationship between the RoCH-equivalence for the
regular orbit reducible RCH systems with symmetries and the
RCH-equivalence for the associated $R_o$-reduced RCH systems, its proof
is given in Marsden et al \cite{mawazh10}. This theorem can be
regarded as an extension of regular orbit reduction theorem of
Hamiltonian systems under regular controlled Hamiltonian equivalence
conditions.

\begin{theo}
If two regular orbit reducible RCH systems $(T^\ast Q_i, G_i,
\omega_i, H_i, F_i,W_i)$, $i=1,2,$ are RoCH-equivalent, then their
associated $R_o$-reduced RCH systems $((T^\ast
Q)_{\mathcal{O}_{\mu_i}}, \omega_{i\mathcal{O}_{\mu_i}},
h_{i\mathcal{O}_{\mu_i}}, f_{i\mathcal{O}_{\mu_i}},\\
W_{i\mathcal{O}_{\mu_i}})$, $i=1,2,$ must be RCH-equivalent.
Conversely, if the $R_o$-reduced RCH systems $((T^\ast
Q)_{\mathcal{O}_{\mu_i}},\\ \omega_{i\mathcal{O}_{\mu_i}},
h_{i\mathcal{O}_{\mu_i}}, f_{i\mathcal{O}_{\mu_i}},
W_{i\mathcal{O}_{\mu_i}})$, $i=1,2,$ are RCH-equivalent and the
induced map
$\varphi^\ast_{\mathcal{O}_\mu}:\mathbf{J}_2^{-1}(\mathcal{O}_{\mu_2})\rightarrow
\mathbf{J}_1^{-1}(\mathcal{O}_{\mu_1})$, such that
$\mathbf{J}_{2\mathcal{O}_{\mu_2}}^\ast
\omega_{2\mathcal{O}_{\mu_2}}^{+}=(\varphi_{\mathcal{O}_\mu}^\ast)^\ast
\cdot\mathbf{J}_{1\mathcal{O}_{\mu_1}}^\ast\omega_{1\mathcal{O}_{\mu_1}}^{+},$
then the regular orbit reducible RCH systems $(T^\ast Q_i, G_i,
\omega_i, H_i, F_i, W_i)$, $i=1,2,$ are RoCH-equivalent.
\end{theo}

Comparing with Theorem 4.7 and Theorem 5.7 we know that,
if the $R_o$-reduced RCH systems $((T^\ast
Q)_{\mathcal{O}_{\mu_i}}, \omega_{i\mathcal{O}_{\mu_i}},
h_{i\mathcal{O}_{\mu_i}}, f_{i\mathcal{O}_{\mu_i}},
W_{i\mathcal{O}_{\mu_i}})$, $i=1,2,$ are RCH-equivalent,
then the regular orbit reducible RCH systems $(T^\ast Q_i, G_i,
\omega_i, H_i, F_i, W_i)$, $i=1,2,$ may not be RoCH-equivalent.
The condition that $\varphi^\ast_{\mathcal{O}_\mu}:
\mathbf{J}_2^{-1}(\mathcal{O}_{\mu_2})\rightarrow
\mathbf{J}_1^{-1}(\mathcal{O}_{\mu_1})$, such that
$\mathbf{J}_{2\mathcal{O}_{\mu_2}}^\ast
\omega_{2\mathcal{O}_{\mu_2}}^{+}=(\varphi_{\mathcal{O}_\mu}^\ast)^\ast
\cdot\mathbf{J}_{1\mathcal{O}_{\mu_1}}^\ast\omega_{1\mathcal{O}_{\mu_1}}^{+},$
is determined by (5.1), which guarantee that
the regular orbit reduction theorem for RCH system holds.
Thus, we emphasize explicitly the impact of different symplectic structures
in the study of RCH-equivalence and regular reduction for RCH systems.\\

Moreover, if considering the RoCH-equivalence of the regular orbit
reducible RCH systems, and using the above Theorem 5.7, Theorem 5.4 and Theorem 2.11,
we can obtain the following Theorem 5.8,
which states that the solutions of two types of Hamilton-Jacobi equations for
the regular orbit reducible RCH systems leave invariant
under the conditions of RoCH-equivalence.

\begin{theo} Suppose that two regular orbit reducible RCH systems
$(T^\ast Q_i,G_i,\omega_i,H_i,F_i,W_i)$,
$i=1,2,$ are RoCH-equivalent with an equivalent map $\varphi: Q_1
\rightarrow Q_2 $ and the associated $R_o$-reduced RCH
systems $((T^\ast Q)_{\mathcal{O}_{\mu_i}},
\omega_{i\mathcal{O}_{\mu_i}}, h_{i\mathcal{O}_{\mu_i}},
f_{i\mathcal{O}_{\mu_i}}, W_{i\mathcal{O}_{\mu_i}}), $ $i=1,2 $.
Under the hypotheses and notations of Theorem 5.2, Theorem 5.3
and Theorem 5.4, then we have that\\

\noindent $(\mathrm{i})$ If the one-form
$\gamma_2: Q_2 \rightarrow T^* Q_2$ is closed with
respect to $T\pi_{Q_2}: TT^* Q_2 \rightarrow TQ_2, $ and
$\bar{\gamma}_2=\pi_{\mathcal{O}_{\mu_2}}(\gamma_2):
Q_2 \rightarrow (T^* Q_2)_{\mathcal{O}_{\mu_2}} $
is a solution of the Type I of Hamilton-Jacobi equation for
the $R_o$-reduced RCH System
$((T^\ast Q_2)_{\mathcal{O}_{\mu_2}}, \omega_{2\mathcal{O}_{\mu_2}},
h_{2\mathcal{O}_{\mu_2}}, f_{2\mathcal{O}_{\mu_2}}, u_{2\mathcal{O}_{\mu_2}})$.
Then $\gamma_1=\varphi^* \cdot \gamma_2\cdot \varphi: Q_1 \rightarrow T^* Q_1, $
is a solution of the Type I of Hamilton-Jacobi equation for the
RCH System $(T^*Q_1,G_1,\omega_1,H_1,F_1,W_1) $,
and $\bar{\gamma}_1=\pi_{\mathcal{O}_{\mu_1}}(\gamma_1):
Q_1 \rightarrow (T^* Q_1)_{\mathcal{O}_{\mu_1}}$ is
a solution of the Type I of Hamilton-Jacobi equation for the $R_o$-reduced
RCH System $((T^\ast Q_1)_{\mathcal{O}_{\mu_1}}, \omega_{1\mathcal{O}_{\mu_1}},
h_{1\mathcal{O}_{\mu_1}}, f_{1\mathcal{O}_{\mu_1}}, u_{1\mathcal{O}_{\mu_1}})$.
Vice versa.\\

\noindent $(\mathrm{ii})$ If the $G_2$-invariant symplectic map
$\varepsilon_2: T^*Q_2\rightarrow T^* Q_2$ and
$\bar{\varepsilon}_2=\pi_{\mathcal{O}_{\mu_2}}(\varepsilon_2):
\mathbf{J_2}^{-1}(\mathcal{O}_{\mu_2})
(\subset T^*Q_2) \rightarrow (T^* Q_2)_{\mathcal{O}_{\mu_2}} $
satisfy the Type II of Hamilton-Jacobi equation for the $R_o$-reduced RCH System
$((T^\ast Q_2)_{\mathcal{O}_{\mu_2}},
\omega_{2\mathcal{O}_{\mu_2}},h_{2\mathcal{O}_{\mu_2}}, f_{2\mathcal{O}_{\mu_2}}, u_{2\mathcal{O}_{\mu_2}})$.
Then $\varepsilon_1= \varphi^* \cdot \varepsilon_2\cdot \varphi_*: T^*Q_1 \rightarrow T^* Q_1 $
and $\bar{\varepsilon}_1=\pi_{\mathcal{O}_{\mu_1}}(\varepsilon_1):
\mathbf{J_1}^{-1}(\mathcal{O}_{\mu_1})
(\subset T^*Q_1) \rightarrow (T^* Q_1)_{\mathcal{O}_{\mu_1}} $
satisfy the Type II of Hamilton-Jacobi equation for the $R_o$-reduced
RCH System $((T^\ast Q_1)_{\mathcal{O}_{\mu_1}},
\omega_{1\mathcal{O}_{\mu_1}},h_{1\mathcal{O}_{\mu_1}}, f_{1\mathcal{O}_{\mu_1}}, u_{1\mathcal{O}_{\mu_1}}). $ Vice versa.
\end{theo}

\noindent{\bf Proof: } We first prove the assertion $(\mathrm{i})$.
If two regular orbit reducible RCH systems $(T^\ast Q_i,G_i,\omega_i,\\ H_i,F_i,W_i)$,
$i=1,2,$ are RoCH-equivalent, from Definition 5.6 we know that the two RCH
systems are also RCH-equivalent. If the one-form
$\gamma_2: Q_2 \rightarrow T^* Q_2$ is closed with
respect to $T\pi_{Q_2}: TT^* Q_2 \rightarrow TQ_2, $ from Theorem 2.6, we know that
$\gamma_2$ is a solution of the Type I of Hamilton-Jacobi equation
for the RCH system $(T^*Q_2,G_2,\omega_2,H_2,F_2,W_2). $
Moreover, from Theorem 2.11, we know that
$\gamma_1=\varphi^* \cdot \gamma_2\cdot \varphi: Q_1 \rightarrow T^* Q_1, $ is
a solution of the Type I of Hamilton-Jacobi equation for the
RCH System $(T^*Q_1,G_1,\omega_1,H_1,F_1,W_1) $.
On the other hand, If $(T^\ast
Q_1,G_1,\omega_1,H_1,F_1,W_1)\stackrel{RoCH}{\sim}(T^\ast
Q_2,G_2,\omega_2,H_2, F_2,W_2)$, then from Definition 5.6 there
exists a diffeomorphism $\varphi: Q_1 \rightarrow Q_2$, such that
$\varphi^\ast: T^\ast Q_2 \rightarrow T^\ast Q_1$ is symplectic, and
for $\mu_i\in \mathfrak{g}_i^\ast, \; i=1,2, $
$\varphi_{\mathcal{O}_\mu}^\ast=i_{\mathcal{O}_{\mu_1}}^{-1}\cdot
\varphi^\ast\cdot
i_{\mathcal{O}_{\mu_2}}:\mathbf{J}_2^{-1}(\mathcal{O}_{\mu_2})\rightarrow
\mathbf{J}_1^{-1}(\mathcal{O}_{\mu_1})$ is $(G_2,G_1)$-equivariant,
$\mathbf{J}_{2\mathcal{O}_{\mu_2}}^\ast\omega_{2\mathcal{O}_{\mu_2}}^+
=(\varphi_{\mathcal{O}_\mu}^\ast)^\ast
\cdot\mathbf{J}_{1\mathcal{O}_{\mu_1}}^\ast\omega_{1\mathcal{O}_{\mu_1}}^+.
$ From the following commutative Diagram-12,
\[
\begin{CD}
Q_2 @> \gamma_2 >> T^\ast Q_2 @<i_{\mathcal{O}_{\mu_2}}<<
\mathbf{J}_2^{-1}(\mathcal{O}_{\mu_2}) @>\pi_{\mathcal{O}_{\mu_2}}
>> (T^\ast Q_2)_{\mathcal{O}_{\mu_2}}\\
@A \varphi AA @V\varphi^\ast VV @V\varphi^\ast_{\mathcal{O}_\mu}
VV @V\varphi^\ast_{\mathcal{O}_{\mu/G}}VV\\
Q_1 @> \gamma_1 >> T^\ast Q_1 @<i_{\mathcal{O}_{\mu_1}}<<
\mathbf{J}_1^{-1}(\mathcal{O}_{\mu_1})
@>\pi_{\mathcal{O}_{\mu_1}}>>(T^\ast Q_1)_{\mathcal{O}_{\mu_1}}
\end{CD}
\]
$$\mbox{Diagram-12}$$
we have a well-defined symplectic map
$\varphi_{\mathcal{O}_\mu/G}^\ast:(T^\ast
Q_2)_{\mathcal{O}_{\mu_2}}\rightarrow (T^\ast
Q_1)_{\mathcal{O}_{\mu_1}}$, such that
$\varphi_{\mathcal{O}_\mu/G}^\ast \cdot
\pi_{\mathcal{O}_{\mu_2}}=\pi_{\mathcal{O}_{\mu_1}}\cdot
\varphi_{\mathcal{O}_\mu}^\ast$, see Marsden et al \cite{mawazh10}.
Then from the regular orbit reduction Theorem 5.7,
we know that the associated $R_o$-reduced RCH
systems $((T^\ast Q)_{\mathcal{O}_{\mu_i}},
\omega_{i\mathcal{O}_{\mu_i}}, h_{i\mathcal{O}_{\mu_i}},
f_{i\mathcal{O}_{\mu_i}}, W_{i\mathcal{O}_{\mu_i}}), $ $i=1,2, $ are
RCH-equivalent with an equivalent map
$\varphi_{\mathcal{O}_\mu/G}^\ast: (T^\ast
Q_2)_{\mathcal{O}_{\mu_2}} \rightarrow (T^\ast
Q_1)_{\mathcal{O}_{\mu_1}}. $ If
$\bar{\gamma}_2=\pi_{\mathcal{O}_{\mu_2}}(\gamma_2): Q_2 \rightarrow
(T^* Q_2)_{\mathcal{O}_{\mu_2}}$ is a solution of the type I of Hamilton-Jacobi
equation for $R_o$-reduced RCH system $((T^\ast Q_2)_{\mathcal{O}_{\mu_2}},
\omega_{2\mathcal{O}_{\mu_2}},h_{2\mathcal{O}_{\mu_2}},
f_{2\mathcal{O}_{\mu_2}},u_{2\mathcal{O}_{\mu_2}})$,
and note that $\bar{\gamma}_1=\pi_{\mathcal{O}_{\mu_1}}(\gamma_1)=
\pi_{\mathcal{O}_{\mu_1}} \cdot \varphi^* \cdot \gamma_2\cdot
\varphi= \varphi_{\mathcal{O}_\mu/G}^\ast \cdot
\pi_{\mathcal{O}_{\mu_2}} \cdot \gamma_2\cdot \varphi=
\varphi_{\mathcal{O}_\mu/G}^\ast \cdot \bar{\gamma}_2\cdot \varphi
$, from Theorem 2.11 we know that $\bar{\gamma}_1=
\varphi_{\mathcal{O}_\mu/G}^\ast \cdot \bar{\gamma}_2\cdot \varphi $
is a solution of the type I of Hamilton-Jacobi equation for the $R_o$-reduced RCH system
$((T^\ast Q_1)_{\mathcal{O}_{\mu_1}},
\omega_{1\mathcal{O}_{\mu_1}},h_{1\mathcal{O}_{\mu_1}},
f_{1\mathcal{O}_{\mu_1}},u_{1\mathcal{O}_{\mu_1}})
$. Because the map $\varphi: Q_1 \rightarrow Q_2 $ is a diffeomorphism,
and $\varphi^\ast: T^\ast Q_2\rightarrow T^\ast Q_1$ is a symplectic isomorphisms,
vice versa. It follows that the assertion $(\mathrm{i})$ of Theorem 5.8 holds.\\

Next, we prove the assertion $(\mathrm{ii})$.
If the $G_2$-invariant symplectic map $\varepsilon_2: T^*Q_2\rightarrow T^* Q_2$ and
$\bar{\varepsilon}_2=\pi_{\mathcal{O}_{\mu_2}}(\varepsilon_2) $
satisfy the Type II of Hamilton-Jacobi equation for
the $R_o$-reduced RCH System
$((T^\ast Q_2)_{\mathcal{O}_{\mu_2}},
\omega_{2\mathcal{O}_{\mu_2}},h_{2\mathcal{O}_{\mu_2}}, f_{2\mathcal{O}_{\mu_2}}, u_{2\mathcal{O}_{\mu_2}})$, from Theorem 5.4 we know that
$\varepsilon_2$ is a solution of the Type
II of Hamilton-Jacobi equation for the regular orbit reducible RCH
system $(T^*Q_2,G_2,\omega_2,H_2,F_2,W_2). $ Since
the two regular orbit reducible RCH systems $(T^\ast Q_i,G_i,\omega_i,H_i,F_i,W_i)$,
$i=1,2,$ are RoCH-equivalent, and hence are also RCH-equivalent, from Theorem 2.11,
then $\varepsilon_1=\varphi^*\cdot \varepsilon_2  \cdot\varphi_*$
is a solution of the Type II of Hamilton-Jacobi equation for the
RCH System $(T^*Q_1,G_1,\omega_1,H_1,F_1,W_1) $. Moreover,
from Theorem 5.4 we know that
$\varepsilon_1$ and
$\bar{\varepsilon}_1=\pi_{\mathcal{O}_{\mu_1}}(\varepsilon_1) $
satisfy the Type II of Hamilton-Jacobi equation for the $R_o$-reduced
RCH System \\ $((T^\ast Q_1)_{\mathcal{O}_{\mu_1}},
\omega_{1\mathcal{O}_{\mu_1}},h_{1\mathcal{O}_{\mu_1}}, f_{1\mathcal{O}_{\mu_1}}, u_{1\mathcal{O}_{\mu_1}}). $ In the same way,
because the map $\varphi: Q_1 \rightarrow Q_2 $ is a diffeomorphism,
and $\varphi^\ast: T^\ast Q_2\rightarrow T^\ast Q_1$ is a symplectic isomorphisms,
vice versa.
We prove the assertion $(\mathrm{ii})$ of Theorem 5.8. \hskip 0.3cm
$\blacksquare$

\begin{rema}
If $(T^\ast Q, \omega)$ is a connected symplectic manifold, and
$\mathbf{J}:T^\ast Q\rightarrow \mathfrak{g}^\ast$ is a
non-equivariant momentum map with a non-equivariance group
one-cocycle $\sigma: G\rightarrow \mathfrak{g}^\ast$, which is
defined by $\sigma(g):=\mathbf{J}(g\cdot
z)-\operatorname{Ad}^\ast_{g^{-1}}\mathbf{J}(z)$, where $g\in G$ and
$z\in T^\ast Q$. Then we know that $\sigma$ produces a new affine
action $\Theta: G\times \mathfrak{g}^\ast \rightarrow
\mathfrak{g}^\ast $ defined by
$\Theta(g,\mu):=\operatorname{Ad}^\ast_{g^{-1}}\mu + \sigma(g)$,
where $\mu \in \mathfrak{g}^\ast$, with respect to which the given
momentum map $\mathbf{J}$ is equivariant. Assume that $G$ acts
freely and properly on $T^\ast Q$, and $\mathcal{O}_\mu= G\cdot \mu
\subset \mathfrak{g}^\ast$ denotes the G-orbit of the point $\mu \in
\mathfrak{g}^\ast$ with respect to this affine action $\Theta$, and
$\mu$ is a regular value of $\mathbf{J}$. Then the quotient space
$(T^\ast Q)_{\mathcal{O}_\mu}=\mathbf{J}^{-1}(\mathcal{O}_\mu)/ G $
is also a symplectic manifold with the $R_o$-reduced symplectic form
$\omega_{\mathcal{O}_\mu}$ uniquely characterized by $(5.1)$, see
Ortega and Ratiu \cite{orra04} and Marsden et al
\cite{mamiorpera07}. In this case, we can also define the regular
orbit reducible RCH system $(T^*Q,G,\omega,H,F,W)$ and
RoCH-equivalence, and prove the two types of Hamilton-Jacobi theorem for the
$R_o$-reduced RCH system $((T^\ast Q)_{\mathcal{O}_\mu},
\omega_{\mathcal{O}_\mu}, h_{\mathcal{O}_\mu},
f_{\mathcal{O}_\mu},u_{\mathcal{O}_\mu})$ by using the above similar
way, and state that the solutions of two types of Hamilton-Jacobi equations for
the regular orbit reducible RCH systems with symmetries leave invariant
under the conditions of RoCH-equivalence, where the $R_o$-reduced
space $((T^\ast Q)_{\mathcal{O}_\mu}, \omega_{\mathcal{O}_\mu})$ is
determined by the affine action.
\end{rema}

\section{Applications}

Now, it is a natural problem
if there is a practical RCH system and how to show the effect on controls
in regular point reduction and Hamilton-Jacobi theory of the system.
In this section, as the applications of the above theoretical results, we
consider first the regular point reducible RCH system on the
generalization of a Lie group, and give precisely the geometric constraint conditions
of the $R_p$-reduced symplectic form for the
dynamical vector field of the regular point reducible RCH system,
that is, the two types of Hamilton-Jacobi equation for
the $R_p$- reduced RCH system. Then
we show the Type I and Type II of Hamilton-Jacobi
equations for rigid body and heavy top with internal rotors on the
generalization of rotation group $\textmd{SO}(3)$ and on the
generalization of Euclidean group $\textmd{SE}(3)$,
respectively. Note that these given equations are more
complex than that of Hamiltonian systems without control, which
describe explicitly the effect on controls in the regular point
reduction and Hamilton-Jacobi theory for RCH systems. We shall follow the notations
and conventions introduced in Marsden et al. \cite{mamora90}, Marsden
and Ratiu \cite{mara99}, Marsden et al. \cite{mawazh10}, and Wang
\cite{wa17}.

\subsection{RCH System on the Generalization of a Lie Group}

In order to describe the two types of Hamilton-Jacobi equations of rigid body and
heavy top with internal rotors, we need to first consider the
regular point reducible RCH system $(T^*Q,G,\omega_Q,H,F,W)$
on the generalization of a Lie
group $Q=G\times V$, where $G$ is a Lie group with Lie algebra
$\mathfrak{g}$ and $V$ is a $k$-dimensional vector space. Defined
the left $G$-action $\Phi: G\times Q \rightarrow Q, \;
\Phi(g,(h,\theta)):=(gh,\theta)$, for any $g,h \in G,\; \theta \in
V$, that is, the $G$-action on $Q$ is by the left translation on the
first factor $G$, and the trivial action on the second factor $V$.
Because locally, $T^\ast Q \cong T^\ast G \times T^\ast V$, and $T^\ast V\cong
V\times V^\ast$, and by using the local left trivialization of $T^\ast G$,
that is, $T^\ast G \cong G \times \mathfrak{g}^\ast, $ where
$\mathfrak{g}^\ast$ is the dual of $\mathfrak{g}$, and hence we have
that locally, $T^\ast Q \cong G \times \mathfrak{g}^\ast \times V \times V^\ast$.
If the left $G$-action $\Phi: G\times Q \rightarrow Q $ is free and
proper, then the cotangent lift of the action to its cotangent
bundle $T^\ast Q$, given by $\Phi^{T^*}: G \times T^*Q \rightarrow
T^*Q, \; \Phi^{T^*}(g,(h,\mu,\theta,\lambda)):=(gh,\mu,\theta,\lambda)$, for
any $g,h \in G,\; \mu \in \mathfrak{g}^\ast, \; \theta \in V, \;
\lambda \in V^\ast$, is also a free and proper action, and the orbit
space $(T^\ast Q)/ G $ is a smooth manifold and $\pi: T^*Q
\rightarrow (T^*Q )/G $ is a smooth submersion. Since $G$ acts
trivially on $\mathfrak{g}^\ast$, $V$ and $V^\ast$, it follows that
$(T^\ast Q)/ G$ is diffeomorphic locally to $\mathfrak{g}^\ast \times V
\times V^\ast$.\\

We know that $\mathfrak{g}^\ast$ is a Poisson manifold with respect
to the $(\pm)$-Lie-Poisson bracket $\{\cdot,\cdot\}_\pm$ defined by
\begin{equation}
\{f,g\}_\pm(\mu):=\pm<\mu,[\frac{\delta f}{\delta \mu},\frac{\delta
g}{\delta\mu}]>,\;\; \forall f,g\in C^\infty(\mathfrak{g}^\ast),\;\;
\mu\in \mathfrak{g}^\ast,\label{6.1}
\end{equation}
where the element $\frac{\delta f}{\delta \mu}\in\mathfrak{g}$ is
defined by the equality $<v,\frac{\delta f}{\delta
\mu}>:=Df(\mu)\cdot v$, for any $v\in \mathfrak{g}^\ast$, see
Marsden and Ratiu \cite{mara99}. Thus, for
the coadjoint orbit $\mathcal{O}_\mu \subset \mathfrak{g}^\ast$, $\mu \in \mathfrak{g}^\ast$,
the induced orbit symplectic forms $\omega^\pm_{\mathcal{O}_\mu}$ are
given by
\begin{equation}\omega_{\mathcal{O}_\mu}^\pm(\nu)(\operatorname{ad}_\xi^\ast(\nu),
\operatorname{ad}_\eta^\ast(\nu))=\pm \langle\nu,[\xi,\eta]\rangle,
\qquad \forall\; \xi,\eta\in\mathfrak{g}, \;\;
\nu\in\mathcal{O}_\mu\subset\mathfrak{g}^\ast, \label{6.2}
\end{equation}
which coincide with the restriction of the Lie-Poisson brackets
on $\mathfrak{g}^\ast$ to the coadjoint orbit $\mathcal{O}_\mu$.
From the symplectic stratification theorem we know that
a finite dimensional Poisson manifold is the disjoint union of
its symplectic leaves, and its each symplectic leaf is an
injective immersed Poisson submanifold
whose induced Poisson structure is symplectic.
In consequence, when $\mathfrak{g}^\ast$ is endowed with one of the Lie
Poisson structures $\{\cdot,\cdot\}_\pm$, the symplectic leaves of
the Poisson manifolds $(\mathfrak{g}^\ast,\{\cdot,\cdot\}_\pm)$
coincide with the connected components of the orbits of the elements
in $\mathfrak{g}^\ast$ under the coadjoint action. From Abraham and
Marsden \cite{abma78}, we know that the coadjoint orbit
$(\mathcal{O}_\mu, \omega_{\mathcal{O}_\mu}^{-}), \; \mu\in
\mathfrak{g}^\ast,$ is symplectically diffeomorphic to
a $R_p$-reduced space $((T^\ast G)_\mu,\omega_\mu)$ of $T^*G$.\\

Let $\omega_V$ be the canonical symplectic form on
$T^\ast V \cong V \times V^\ast$ given by
$$\omega_V((\theta_1, \lambda_1),(\theta_2,
\lambda_2))=<\lambda_2,\theta_1> -<\lambda_1,\theta_2>,$$ where
$(\theta_i, \lambda_i)\in V\times V^\ast, \; i=1,2$, $<\cdot,\cdot>$
is the natural pairing between $V^\ast$ and $V$. Thus, we can induce
the symplectic forms $\tilde{\omega}^\pm_{\mathcal{O}_\mu \times V
\times V^\ast}= \pi_{\mathcal{O}_\mu}^\ast
\omega^\pm_{\mathcal{O}_\mu}+ \pi_V^\ast \omega_V$ on the smooth
manifold $\mathcal{O}_\mu \times V \times V^\ast$, where the maps
$\pi_{\mathcal{O}_\mu}: \mathcal{O}_\mu \times V \times V^\ast \to
\mathcal{O}_\mu$ and $\pi_V: \mathcal{O}_\mu \times V \times V^\ast
\to V\times V^\ast$ are canonical projections. On the other hand,
note that for $F,K: T^\ast V \cong V \times V^\ast \to \mathbb{R}$,
by using the canonical symplectic form $\omega_V$ on $T^\ast V \cong
V \times V^\ast$, we can define the Poisson bracket
$\{\cdot,\cdot\}_V$ on $T^\ast V$ as follows
$$\{F,K\}_V(\theta,\lambda)=<\frac{\delta F}{\delta \theta},
\frac{\delta K}{\delta \lambda}>- <\frac{\delta K}{\delta
\theta},\frac{\delta F}{\delta \lambda}>
$$
If $\theta_i, \; i=1,\cdots, k,$ is a base of $V$, and $\lambda_i,
\; i=1,\cdots, k,$ a base of $V^\ast$, then we have that
\begin{equation}\{F,K\}_V(\theta,\lambda)=\sum_{i=1}^k(\frac{\partial F}{\partial \theta_i}
\frac{\partial K}{\partial \lambda_i}- \frac{\partial K}{\partial
\theta_i}\frac{\partial F}{\partial \lambda_i}). \; \label{6.3}
\end{equation}
Thus, by the $(\pm)$-Lie-Poisson brackets on $\mathfrak{g}^\ast$ and
the Poisson bracket $\{\cdot,\cdot\}_V$ on $T^\ast V$, for $F,K:
\mathfrak{g}^\ast \times V \times V^\ast \to \mathbb{R}$, we can
define the Poisson bracket on $\mathfrak{g}^\ast \times V \times
V^\ast$ as follows
\begin{align} \{F,K\}_\pm(\mu,\theta,\lambda) & = \{F,K\}_\pm(\mu)+
\{F,K\}_V(\theta,\lambda) \nonumber \\
& = \pm<\mu,[\frac{\delta F}{\delta
\mu},\frac{\delta K}{\delta\mu}]>+
 \sum_{i=1}^k(\frac{\partial F}{\partial \theta_i} \frac{\partial
K}{\partial \lambda_i}- \frac{\partial K}{\partial
\theta_i}\frac{\partial F}{\partial \lambda_i}). \; \label{6.4}
\end{align}
See Krishnaprasad and Marsden \cite{krma87}. In particular, for
$F_\mu,K_\mu: \mathcal{O}_\mu \times V \times V^\ast \to
\mathbb{R}$, we have that the induced symplectic form
from the above Poisson bracket on $\mathcal{O}_\mu \times V\times V^\ast $
is given by
$\tilde{\omega}_{\mathcal{O}_\mu \times V
\times V^\ast}^{\pm}(X_{F_\mu}, X_{K_\mu})=
\{F_\mu,K_\mu\}_{\pm}|_{\mathcal{O}_\mu \times V \times V^\ast}. $\\

On the other hand, we note that the cotangent bundle $T^\ast Q$ has a canonical
symplectic form $\omega_Q$, from locally, $T^\ast Q \cong T^\ast G \times T^\ast V$,
we have that $\omega_Q= \pi^\ast_1 \omega_0 +\pi^\ast_2 \omega_V$
on $T^\ast Q$, where $\omega_0$ is
the canonical symplectic form on $T^\ast G$ and the maps $\pi_1:
Q=G\times V \to G$ and $\pi_2: Q=G\times V \to V$ are canonical
projections. Assume that the cotangent lift of the left $G$-action
$\Phi^{T^*}: G \times T^\ast Q \to T^\ast Q$ is symplectic, and
admits an associated $\operatorname{Ad}^\ast$-equivariant momentum
map $\mathbf{J}_Q: T^\ast Q \to \mathfrak{g}^\ast$ such that
$\mathbf{J}_Q\cdot \pi^\ast_1=\mathbf{J}_L$, where
$\mathbf{J}_L:T^\ast G \rightarrow \mathfrak{g}^\ast$ is a momentum
map of the left G-translation on $T^\ast G$ and we assume that it exists,
and $\pi^\ast_1: T^\ast G \to T^\ast Q$.
If $\mu\in\mathfrak{g}^\ast$ is a regular value of
$\mathbf{J}_Q$, then $\mu\in\mathfrak{g}^\ast$ is also a regular
value of $\mathbf{J}_L$ and $\mathbf{J}_Q^{-1}(\mu)\cong
\mathbf{J}_L^{-1}(\mu)\times V \times V^\ast$. Denote
$G_\mu=\{g\in G|\operatorname{Ad}_g^\ast \mu=\mu \}$ the isotropy
subgroup of co-adjoint $G$-action at the point
$\mu\in\mathfrak{g}^\ast$. It follows that $G_\mu$ acts also freely
and properly on $\mathbf{J}_Q^{-1}(\mu)$, the $R_p$-reduced
space $(T^\ast Q)_\mu=\mathbf{J}_Q^{-1}(\mu)/G_\mu\cong (T^\ast
G)_\mu \times V \times V^\ast$ of $(T^\ast Q,\omega_Q)$ at $\mu$, is
a symplectic manifold with symplectic form $\omega_\mu$ uniquely
characterized by the relation $\pi_\mu^\ast \omega_\mu=i_\mu^\ast
\omega_Q=i_\mu^\ast \pi^\ast_1 \omega_0 +i_\mu^\ast \pi^\ast_2
\omega_V$, where the map $i_\mu:\mathbf{J}_Q^{-1}(\mu)\rightarrow
T^\ast Q$ is the inclusion and
$\pi_\mu:\mathbf{J}_Q^{-1}(\mu)\rightarrow (T^\ast Q)_\mu$ is the
projection. From Abraham and Marsden \cite{abma78}, we know
that $((T^\ast G)_\mu,\omega_\mu)$ is symplectically diffeomorphic
to $(\mathcal{O}_\mu,\omega_{\mathcal{O}_\mu}^{-})$, and hence we
have that $((T^\ast Q)_\mu,\omega_\mu)$ is symplectically
diffeomorphic to $(\mathcal{O}_\mu \times V\times
V^\ast,\tilde{\omega}_{\mathcal{O}_\mu \times V \times
V^\ast}^{-})$, which is a symplectic leaf of the Poisson manifold
$(\mathfrak{g}^\ast \times V \times V^\ast, \{\cdot,\cdot\}_{-}). $\\

In the following we shall give the regular point reduction of the RCH system
$(T^\ast Q,G,\omega_Q,H,F,W)$ with symmetry and momentum map
on the generalization of a Lie group.
Now we identify $TG$ and $G\times \mathfrak{g}$ locally,
by using the left translation, and $TV\cong V\times V$,
then locally, $TQ\cong G\times \mathfrak{g} \times V\times V$.
In consequence, we consider the Lagrangian $L(g,\xi,\theta,\dot{\theta}):TQ
\cong G\times \mathfrak{g} \times V\times V \to \mathbb{R}$, which is usually
the total kinetic minus potential energy of the system, where
$(g,\xi) \in G\times \mathfrak{g}$, and $\theta \in V$, $\xi^i$ and
$\dot{\theta}^j=\frac{\mathrm{d} \theta^j}{\mathrm{d} t}$,
($i=1,\cdots,n, \; j=1,\cdots,k$, $n=\dim G$, $k=\dim V$), regarded
as the velocity variables of the system. If we introduce the conjugate momentum
$p_i=\frac{\partial L}{\partial \xi^i}, \; l_j=\frac{\partial
L}{\partial \dot{\theta}^j}$, $i=1,\cdots,n, \; j=1,\cdots,k,$ and
by the Legendre transformation $FL:TQ \cong G\times
\mathfrak{g}\times V \times V \to T^\ast Q \cong G\times
\mathfrak{g}^\ast \times V \times V^\ast$,
$(g^i,\xi^i,\theta^j,\dot{\theta}^j)\to (g^i,p_i,\theta^j,l_j)$, we
have the Hamiltonian $H(g,p,\theta,l):T^\ast Q \cong G\times
\mathfrak{g}^\ast \times V \times V^\ast \to \mathbb{R}$ given by
\begin{equation}
H(g^i,p_i,\theta^j,l_j)=\sum_{i=1}^{n}p_i\xi^i+\sum_{j=1}^{k}l_j\dot{\theta}^j
-L(g^i,\xi^i,\theta^j,\dot{\theta}^j). \; \label{6.5}
\end{equation}

If the Hamiltonian $H(g,p,\theta,l): T^\ast Q \cong G\times
\mathfrak{g}^\ast \times V \times V^\ast \to \mathbb{R}$ is the left
cotangent lifted $G$-action $\Phi^{T^*}$ invariant, for
$\mu\in\mathfrak{g}^\ast$ the regular value of
the momentum map $\mathbf{J}_Q: T^\ast Q \rightarrow
\mathfrak{g}^\ast$, we have the associated $R_p$-reduced Hamiltonian
$h_\mu(\nu,\theta,l): (T^\ast Q)_\mu \cong\mathcal{O}_\mu \times V
\times V^\ast \to \mathbb{R}$, defined by $h_\mu\cdot \pi_\mu=H\cdot
i_\mu$, and the $R_p$-reduced Hamiltonian vector field $X_{h_\mu}$ given
by $X_{h_\mu}(K_\mu)=\{K_\mu,h_\mu\}_{-}|_{\mathcal{O}_\mu \times V
\times V^\ast}, $ where
$K_\mu(\nu,\theta,l): \mathcal{O}_\mu
\times V\times V^* \to \mathbb{R}.$\\

Thus, if we consider that the fiber-preserving map $F: T^\ast Q \to
T^\ast Q$ and the fiber submanifold $W$ of $T^\ast Q$ are all left
cotangent lifted $G$-action $\Phi^{T^*}$ invariant, then for $u \in W$, the 6-tuple
$(T^\ast Q,G,\omega_Q,H,F,u)$ is a regular point reducible RCH
system with a control law $u$, and its dynamical vector field can be expressed by
\begin{align}
\tilde{X}= X_{(T^\ast Q,G,\omega_Q,H,F,u)}= X_H
+ \textnormal{vlift}(F)+ \textnormal{vlift}(u), \; \label{6.6}
\end{align}
where $\textnormal{vlift}(F)= \textnormal{vlift}(F)\cdot X_H, \;\;
\textnormal{vlift}(u)= \textnormal{vlift}(u)\cdot X_H $ are
the changes of $X_H$ under the action of the external force $F$ and the control law $u$.\\

For a point $\mu\in\mathfrak{g}^\ast$, the regular value of
the momentum map $\mathbf{J}_Q $, assume that
$F(\mathbf{J}_Q^{-1}(\mu))\subset \mathbf{J}_Q^{-1}(\mu) $,
$f_\mu\cdot \pi_\mu=\pi_\mu\cdot F\cdot i_\mu$,
$u \in (W \cap \mathbf{J}_Q^{-1}(\mu))\neq \emptyset $,
and $u_\mu \in
W_\mu=\pi_\mu(W \cap \mathbf{J}_Q^{-1}(\mu))\subset
\mathcal{O}_\mu \times V \times V^\ast$,
$u_\mu\cdot \pi_\mu=\pi_\mu\cdot u\cdot i_\mu$,
then the $R_p$-reduced system is the 5-tuple
$(\mathcal{O}_\mu \times V\times
V^\ast,\tilde{\omega}_{\mathcal{O}_\mu \times V \times V^\ast
}^{-},h_\mu,f_\mu,W_\mu)$, where $\mathcal{O}_\mu \subset
\mathfrak{g}^\ast$ is the co-adjoint orbit,
$\tilde{\omega}_{\mathcal{O}_\mu \times V \times V^\ast }^{-}$ is
the induced symplectic form on $\mathcal{O}_\mu \times V\times V^\ast $.
Moreover, we can give precisely the geometric constraint
conditions of the $R_p$-reduced symplectic form
$\tilde{\omega}_{\mathcal{O}_\mu \times V\times
V*}$ for the dynamical vector field
of the regular point reducible RCH system,
that is, Type I and Type II of
Hamilton-Jacobi equation for the $R_p$-reduced RCH system.
Assume that $\gamma: Q \rightarrow T^* Q$ is an one-form on $Q$, and
$\lambda=\gamma \cdot \pi_{Q}: T^* Q \rightarrow T^* Q $, and
$\varepsilon: T^* Q \rightarrow T^* Q $ is a $G_\mu$-invariant
symplectic map, such that
$\varepsilon(\mathbf{J}_Q^{-1}(\mu))\subset \mathbf{J}_Q^{-1}(\mu),
$ and $\textmd{Im}(\gamma)\subset \mathbf{J}_Q^{-1}(\mu), $ and it
is $G_\mu$-invariant, where $G_\mu= \{g\in
G|\operatorname{Ad}_g^\ast \mu=\mu \}$ is the isotropy subgroup of
the coadjoint $G$-action at the point $\mu\in\mathfrak{g}^\ast$.
Denote $\bar{\gamma}=\pi_\mu(\gamma): Q \rightarrow
\mathcal{O}_\mu\times V\times
V^\ast, $ and $\bar{\lambda}=\pi_\mu(\lambda):
\mathbf{J}_Q^{-1}(\mu) \rightarrow \mathcal{O}_\mu\times V\times
V^\ast, $ and
$\bar{\varepsilon}=\pi_\mu(\varepsilon):
\mathbf{J}_Q^{-1}(\mu)\rightarrow \mathcal{O}_\mu\times V\times
V^\ast, $ where
$\pi_\mu:\mathbf{J}_Q^{-1}(\mu)\rightarrow \mathcal{O}_\mu\times V\times
V^\ast$ is the projection. By using the similar way in
the proofs of two types of Hamilton-Jacobi theorem for the $R_p$-reduced RCH
system, see Theorem 4.2 and Theorem 4.3, we can obtain the following theorem.

\begin{theo}
If the 6-tuple $(T^\ast Q,G,\omega_Q,H,F,W)$ is a regular point reducible RCH system
on the generalization of a Lie group $Q=G\times
V$, where $G$ is a Lie group and $V$ is a $k$-dimensional vector
space, then for $\mu \in \mathfrak{g}^\ast $, the regular value of
the momentum map $\mathbf{J}_Q: T^\ast Q \rightarrow
\mathfrak{g}^\ast$,
the $R_p$-reduced system is the 5-tuple
$(\mathcal{O}_\mu \times V\times
V^\ast,\tilde{\omega}_{\mathcal{O}_\mu \times V \times V^\ast
}^{-},h_\mu,f_\mu,W_\mu)$, where $\mathcal{O}_\mu \subset
\mathfrak{g}^\ast$ is the co-adjoint orbit,
$\tilde{\omega}_{\mathcal{O}_\mu \times V \times V^\ast }^{-}$ is
the induced symplectic form on $\mathcal{O}_\mu \times V\times V^\ast $,
$h_\mu\cdot \pi_\mu=H\cdot i_\mu$,
$F(\mathbf{J}_Q^{-1}(\mu))\subset \mathbf{J}_Q^{-1}(\mu) $,
$f_\mu\cdot \pi_\mu=\pi_\mu\cdot F\cdot i_\mu$,
$u \in (W \cap \mathbf{J}_Q^{-1}(\mu))\neq \emptyset $,
and $u_\mu \in
W_\mu=\pi_\mu(W \cap \mathbf{J}_Q^{-1}(\mu))\subset
\mathcal{O}_\mu \times V \times V^\ast$,
$u_\mu\cdot \pi_\mu=\pi_\mu\cdot u\cdot i_\mu$.
Assume that $\gamma: Q \rightarrow T^* Q$ is an one-form on
$Q$, and $\lambda=\gamma \cdot \pi_{Q}: T^* Q \rightarrow T^* Q $, and
$\varepsilon: T^* Q \rightarrow T^* Q $ is a
symplectic map. Denote $\tilde{X}^\gamma = T\pi_{Q}\cdot \tilde{X} \cdot \gamma$,
and $\tilde{X}^\varepsilon = T\pi_{Q}\cdot \tilde{X} \cdot \varepsilon$,
where $\tilde{X}=X_{(T^\ast Q,G,\omega_Q,H,F,u)}$ is the dynamical
vector field of the regular point reducible RCH system
$(T^*Q,G,\omega_Q,H,F,W)$ with a control law $u$. Moreover, assume
that $\textmd{Im}(\gamma)\subset
\mathbf{J}_Q^{-1}(\mu), $ and it is $G_\mu$-invariant, and
$\varepsilon$ is $G_\mu$-invariant
and $\varepsilon(\mathbf{J}_Q^{-1}(\mu))\subset \mathbf{J}_Q^{-1}(\mu).$
Denote
$\bar{\gamma}=\pi_\mu(\gamma): Q \rightarrow \mathcal{O}_\mu \times
V\times V^\ast, $ and $\bar{\lambda}=\pi_\mu(\lambda):
\mathbf{J}_Q^{-1}(\mu) \rightarrow \mathcal{O}_\mu\times V\times
V^\ast, $ and
$\bar{\varepsilon}=\pi_\mu(\varepsilon):
\mathbf{J}_Q^{-1}(\mu)\rightarrow \mathcal{O}_\mu\times V\times
V^\ast. $
Then the following two assertions hold:\\
\noindent $(\mathbf{i})$
If the one-form $\gamma: Q \rightarrow T^*Q $ is closed with respect to
$T\pi_Q: TT^* Q \rightarrow TQ, $
then $\bar{\gamma}$ is a solution of the Type I of Hamilton-Jacobi equation
$T\bar{\gamma}\cdot \tilde{X}^\gamma= X_{h_\mu}\cdot \bar{\gamma}; $\\
\noindent $(\mathbf{ii})$
The $\varepsilon$ and $\bar{\varepsilon} $ satisfy the Type II of Hamilton-Jacobi equation
$T\bar{\gamma}\cdot \tilde{X}^\varepsilon= X_{h_\mu}\cdot \bar{\varepsilon}, $
if and only if they satisfy
the equation $T\bar{\varepsilon}\cdot(X_{h_\mu \cdot \bar{\varepsilon}})
= T\bar{\lambda}\cdot \tilde{X} \cdot\varepsilon. $
Here $X_{h_{\mu}}$ and $ X_{h_{\mu} \cdot \bar{\varepsilon}} \in TT^*Q $
are the Hamiltonian vector fields of the $R_p$-reduced Hamiltonian
functions $h_{\mu}$ and $h_{\mu} \cdot \bar{\varepsilon}: T^*Q\rightarrow
\mathbb{R}, $ respectively. \hskip 0.3cm $\blacksquare$
\end{theo}

In particular, when $Q=G, $ we can obtain the two types of Hamilton-Jacobi
equation for the $R_p$-reduced RCH system on Lie group $G$.
In this case, if the RCH system we considered has not any control,
then from the above Theorem 6.1, we can get the Theorem 5.1
in Wang \cite{wa17}, that is, the two types of Hamilton-Jacobi
equation for the Marsden-Weinstein reduced Hamiltonian system on a Lie group.
In addition, note that the symplecitic structure
on the co-adjoint orbit
$\mathcal{O}_\mu$ is induced by the (-)-Lie-Poisson brackets on
$\mathfrak{g}^\ast$, then the Type I and Type II of Hamilton-Jacobi equation
$T\bar{\gamma}\cdot \tilde{X}^\gamma= X_{h_\mu}\cdot \bar{\gamma}$
and
$T\bar{\gamma}\cdot \tilde{X}^\varepsilon= X_{h_\mu}\cdot \bar{\varepsilon}, $
for $R_p$-reduced RCH system $(\mathcal{O}_\mu,
\omega_{\mathcal{O}_\mu}^{-}, h_\mu, f_\mu,u_\mu )$ are also
called the Type I and Type II of
Lie-Poisson Hamilton-Jacobi equation, respectively.
See Wang \cite{wa17}, Marsden
and Ratiu \cite{mara99}, and Ge and Marsden \cite{gema88}.

\subsection {Hamilton-Jacobi Equations of Rigid Body with Internal Rotors}

In the following we regard the rigid body with three symmetric
internal rotors as a regular point reducible RCH system on the
generalization of rotation group $\textmd{SO}(3)\times
\mathbb{R}^3$, and give its two types of Hamilton-Jacobi equation
by calculation in detail. Note that our description of
the motion and the equations of rigid body with internal rotors in
this subsection follows some of the notations and conventions in
Marsden \cite{ma92}, Marsden and Ratiu \cite{mara99},  Marsden et
al. \cite{mawazh10}, Wang \cite{wa17}.\\

We consider a rigid body carrying
three "non-mass" internal rotors which is called a carrier body,
where "non-mass" means that the mass of a rotor is very very small
comparing with the mass of the rigid body.
Denote the system center of mass by $O$ in
the carrier body frame and at $O$ place a set of (orthonormal) body axes,
and assume that the rotor and the body coordinate axes are aligned with
principal axes of the carrier body.
If translations are ignored and only rotations of the rigid body-rotor system
are considered, and the rotor spins under the
influence of a torque $u$ acting on the rotor.
In this case, the configuration
space is $Q=\textmd{SO}(3)\times V$, where $V=S^1\times S^1\times
S^1$, with the first factor being carrier body attitude and the second
factor being the angles of rotors. The corresponding phase space is
the cotangent bundle $T^\ast Q$ and locally,
$T^\ast Q \cong T^\ast \textmd{SO}(3)\times T^\ast V$,
where $T^\ast V = T^\ast (S^1\times S^1\times S^1)\cong T^\ast
\mathbb{R}^3$ locally, with the canonical symplectic form $\omega_Q$.
By using the local left trivialization, locally,
$T^\ast \textmd{SO}(3)\cong \textmd{SO}(3)\times \mathfrak{so}^\ast(3)$
and $T^*\mathbb{R}^3 \cong \mathbb{R}^3 \times \mathbb{R}^{3*}$,
then we have that locally,
$T^*Q \cong \textmd{SO}(3)\times \mathfrak{so}^\ast(3)
\times \mathbb{R}^3 \times \mathbb{R}^{3*}$.
For convenience, in the following we denote uniformly that, locally,
$Q= \textmd{SO}(3)\times \mathbb{R}^3, $ and
$T^* Q= T^*(\textmd{SO}(3)\times \mathbb{R}^3)
\cong \textmd{SO}(3)\times \mathfrak{so}^\ast(3)
\times \mathbb{R}^3 \times \mathbb{R}^{3*}$.
Assume that Lie
group $G=\textmd{SO}(3)$ acts freely and properly on $Q$ by the left
translation on the first factor $\textmd{SO}(3)$ and
the trivial action on the second factor $\mathbb{R}^3$.
Then the action of $\textmd{SO}(3)$
on the phase space $T^\ast Q$ is by the cotangent lift of  the left
$\textmd{SO}(3)$ action on $Q$, that is, $\Phi^{T*}:
\textmd{SO}(3)\times T^\ast Q \cong
\textmd{SO}(3)\times \textmd{SO}(3)\times \mathfrak{so}^\ast(3)
\times \mathbb{R}^3 \times \mathbb{R}^{3*}\to T^* Q \cong \textmd{SO}(3)\times
\mathfrak{so}^\ast(3)\times \mathbb{R}^3 \times \mathbb{R}^{3*},$ given
by $\Phi^{T*}(B,(A,\Pi,\alpha, l))=(BA,\Pi,\alpha, l)$, for any $A,B\in
\textmd{SO}(3), \; \Pi \in \mathfrak{so}^\ast(3), \; \alpha \in
\mathbb{R}^3, \; l \in
\mathbb{R}^{3*}$, which is also free and proper.
Assume that the left $\textmd{SO}(3)$ action $\Phi^{T*}$
is symplectic and admits an
associated $\operatorname{Ad}^\ast$-equivariant momentum map
$\mathbf{J}_Q: T^\ast Q \cong \textmd{SO}(3)\times
\mathfrak{so}^\ast(3) \times \mathbb{R}^3 \times \mathbb{R}^{3*} \to
\mathfrak{so}^\ast(3)$ for the cotangent lift of the left $\textmd{SO}(3)$ action.
If $\Pi \in \mathfrak{so}^\ast(3)$ is a regular value of $\mathbf{J}_Q$,
then the regular point reduced space $(T^\ast Q)_\Pi=
\mathbf{J}^{-1}_Q(\Pi)/\textmd{SO}(3)_\Pi$ is symplectically
diffeomorphic to the coadjoint orbit $\mathcal{O}_\Pi \times
\mathbb{R}^3 \times \mathbb{R}^{3*} \subset \mathfrak{so}^\ast(3)
\times \mathbb{R}^3 \times \mathbb{R}^{3*}$.\\

Let $I=diag(I_1,I_2,I_3)$ be the the matrix of inertia moment of the rigid
body in the body fixed frame, which is a principal body frame,
and $J_i,\; i=1,2,3$ be the
moments of inertia of rotors around their rotation axes. Let
$J_{ik},\; i=1,2,3, \;k=1,2,3,$ be the moments of inertia of the
$i$th rotor with $i=1,2,3,$ around the $k$th principal axis with
$k=1,2,3,$ respectively, and denote
$\bar{I}_i=I_i+J_{1i}+J_{2i}+J_{3i}-J_{ii}, \; i=1,2,3$. Let
$\Omega=(\Omega_1,\Omega_2,\Omega_3)$ be the angular velocity vector of rigid body-rotors
computed with respect to the axes fixed in the carrier body and
$(\Omega_1,\Omega_2,\Omega_3)\in \mathfrak{so}(3)$. Let $\alpha_i,\;
i=1,2,3,$ be the relative angles of rotors and
$\dot{\alpha}=(\dot{\alpha_1},\dot{\alpha_2},\dot{\alpha_3})$ the
relative angular velocity vector of rotor about the principal axes
with respect to a carrier body fixed frame.
For convenience, we assume the total mass of the system $m=1$.\\

Now, by using the local left trivialization, locally,
$T\textmd{SO}(3)\cong \textmd{SO}(3)\times \mathfrak{so}(3)$
and $T\mathbb{R}^3 \cong \mathbb{R}^3 \times \mathbb{R}^{3}$,
then we have that locally,
$TQ \cong \textmd{SO}(3)\times \mathfrak{so}(3)
\times \mathbb{R}^3 \times \mathbb{R}^{3}$.
We consider the Lagrangian of the system
$L(A,\Omega,\alpha,\dot{\alpha}):TQ\cong
\textmd{SO}(3)\times\mathfrak{so}(3)\times\mathbb{R}^3\times\mathbb{R}^3\to
\mathbb{R}$, which is the total kinetic energy of the rigid body
plus the total kinetic energy of the rotor, given by
$$L(A,\Omega,\alpha,\dot{\alpha})=\dfrac{1}{2}[\bar{I}_1\Omega_1^2
+\bar{I}_2\Omega_2^2+\bar{I}_3\Omega_3^2
+J_1(\Omega_1+\dot{\alpha}_1)^2+J_2(\Omega_2+\dot{\alpha}_2)^2
+J_3(\Omega_3+\dot{\alpha}_3)^2],$$ where $A\in \textmd{SO}(3)$,
$\Omega=(\Omega_1,\Omega_2,\Omega_3)\in \mathfrak{so}(3)$,
$\alpha=(\alpha_1,\alpha_2,\alpha_3)\in \mathbb{R}^3$,
$\dot{\alpha}=(\dot{\alpha}_1,\dot{\alpha}_2,\dot{\alpha}_3)\in
\mathbb{R}^3$. If we introduce the conjugate angular momentum, which
is given by
$$\Pi_i= \dfrac{\partial L}{\partial \Omega_i}
=\bar{I}_i\Omega_i+J_i(\Omega_i+\dot{\alpha}_i),\quad l_i =
\dfrac{\partial L}{\partial \dot{\alpha}_i}
=J_i(\Omega_i+\dot{\alpha}_i),\quad i=1,2,3,$$ and by the Legendre
transformation $FL:TQ\cong
\textmd{SO}(3)\times\mathfrak{so}(3)\times\mathbb{R}^3\times\mathbb{R}^3\to
T^\ast Q\cong \textmd{SO}(3)\times
\mathfrak{so}^\ast(3)\times\mathbb{R}^3\times\mathbb{R}^{3*},\quad
(A,\Omega,\alpha,\dot{\alpha})\to(A,\Pi,\alpha,l)$, where
$\Pi=(\Pi_1,\Pi_2,\Pi_3)\in \mathfrak{so}^\ast(3)$,
$l=(l_1,l_2,l_3)\in \mathbb{R}^{3*}$, we have the Hamiltonian
$H(A,\Pi,\alpha,l):T^\ast Q\cong \textmd{SO}(3)\times
\mathfrak{so}^\ast (3)\times\mathbb{R}^3\times\mathbb{R}^{3*}\to
\mathbb{R}$ given by
\begin{align*} H(A,\Pi,\alpha,l)&=\Omega\cdot
\Pi+\dot{\alpha}\cdot l-L(A,\Omega,\alpha,\dot{\alpha})\\
&=\bar{I}_1\Omega_1^2+J_1(\Omega_1^2+\Omega_1\dot{\alpha}_1)
+\bar{I}_2\Omega_2^2+J_2(\Omega_2^2+\Omega_2\dot{\alpha}_2)
+\bar{I}_3\Omega_3^2+J_3(\Omega_3^2\\
&\quad+\Omega_3\dot{\alpha}_3)+J_1(\dot{\alpha}_1\Omega_1+\dot{\alpha}_1^2)
  +J_2(\dot{\alpha}_2\Omega_2+\dot{\alpha}_2^2)
  +J_3(\dot{\alpha}_3\Omega_3+\dot{\alpha}_3^2)\\
&\quad
-\frac{1}{2}[\bar{I}_1\Omega_1^2+\bar{I}_2\Omega_2^2+\bar{I}_3\Omega_3^2
+J_1(\Omega_1+\dot{\alpha}_1)^2+J_2(\Omega_2+\dot{\alpha}_2)^2 +J_3(\Omega_3+\dot{\alpha}_3)^2]\\
&=\frac{1}{2}[\frac{(\Pi_1-l_1)^2}{\bar{I}_1}+\frac{(\Pi_2-l_2)^2}{\bar{I}_2}
+\frac{(\Pi_3-l_3)^2}{\bar{I}_3}+\frac{l_1^2}{J_1}+\frac{l_2^2}{J_2}+\frac{l_3^2}{J_3}].
\end{align*}
From the above expression of the Hamiltonian, we know that
$H(A,\Pi,\alpha,l)$ is invariant under the cotangent lift of the left
$\textmd{SO}(3)$-action $\Phi^{T*}: \textmd{SO}(3)\times T^\ast Q \to
T^\ast Q$. Moreover,
from the Lie-Poisson bracket of rigid body on $\mathfrak{so}^\ast(3)$, that is, for $F,K:
\mathfrak{so}^\ast(3)\to \mathbb{R}, $ we have that
$\{F,K\}_{-}(\Pi)=-\Pi\cdot(\nabla_\Pi F\times \nabla_\Pi K), $ and
the Poisson bracket on
$T^\ast \mathbb{R}^3$, we can get the Poisson bracket on
$\mathfrak{so}^\ast(3)\times\mathbb{R}^3\times\mathbb{R}^{3*}$,
that is, for $F,K:
\mathfrak{so}^\ast(3)\times\mathbb{R}^3\times\mathbb{R}^{3*} \to
\mathbb{R}, $ we have that
\begin{align}
\{F,K\}_{-}(\Pi,\alpha,l)=-\Pi\cdot(\nabla_\Pi F\times \nabla_\Pi
K)+ \{F,K\}_{\mathbb{R}^3}(\alpha,l). \; \label{6.7}
\end{align}
Hence, the Hamiltonian vector field $X_H$ of rigid body-rotor system is given by
\begin{align*}
X_{H}(\Pi)& =\{\Pi,\; H\}_{-}=
-\Pi\cdot(\nabla_\Pi\Pi\times\nabla_\Pi
H)+ \{\Pi,\; H\}_{\mathbb{R}^3}\\
&= -\nabla_\Pi\Pi\cdot(\nabla_\Pi H\times \Pi)+ \sum_{i=1}^3(\frac{\partial \Pi}{\partial \alpha_i}
\frac{\partial H}{\partial l_i}- \frac{\partial
H}{\partial \alpha_i}\frac{\partial \Pi}{\partial l_i})\\
& =(\Pi_1,\Pi_2,\Pi_3)\times (\frac{(\Pi_1- l_1)}{ \bar{I}_1},\;\;
\frac{(\Pi_2- l_2)}{ \bar{I}_2}, \;\; \frac{(\Pi_3-
l_3)}{\bar{I}_3}) \\
&= ( \frac{(\bar{I}_2-\bar{I}_3)\Pi_2\Pi_3-
\bar{I}_2\Pi_2l_3 +
\bar{I}_3\Pi_3l_2}{\bar{I}_2\bar{I}_3}, \;\;
\frac{(\bar{I}_3-\bar{I}_1)\Pi_3\Pi_1-
\bar{I}_3\Pi_3l_1 +
\bar{I}_1\Pi_1l_3}{\bar{I}_3\bar{I}_1}, \\
& \;\;\;\;\;\;
\frac{(\bar{I}_1-\bar{I}_2)\Pi_1\Pi_2-
\bar{I}_1\Pi_1l_2 + \bar{I}_2\Pi_2l_1}{\bar{I}_1\bar{I}_2} ),
\end{align*}
since $\nabla_{\Pi_i}\Pi_i=1,\; \nabla_{\Pi_i}\Pi_j=0, \; i\neq j $
and $\nabla_{\Pi_j} H= (\Pi_j-l_j)/\bar{I}_j , \;  \frac{\partial
\Pi}{\partial \alpha_i}= \frac{\partial H}{\partial
\alpha_i}=0, \; i, j=1,2,3 $.

\begin{align*}
X_{H}(\alpha)& =\{\alpha,\; H\}_{-}=
-\Pi\cdot(\nabla_\Pi\alpha\times\nabla_\Pi
H)+ \{\alpha,\; H\}_{\mathbb{R}^3}\\
&= -\nabla_\Pi \alpha\cdot(\nabla_\Pi H\times \Pi)
+ \sum_{i=1}^3(\frac{\partial \alpha}{\partial \alpha_i}
\frac{\partial H}{\partial l_i}- \frac{\partial
H}{\partial \alpha_i}\frac{\partial \alpha}{\partial l_i})\\
& = ( -\frac{(\Pi_1- l_1)}{ \bar{I}_1}
+\frac{l_1}{J_1},\;\; -\frac{(\Pi_2-
l_2)}{ \bar{I}_2} +\frac{l_2}{J_2},\;\;
-\frac{(\Pi_3- l_3)}{\bar{I}_3}
+\frac{l_3}{J_3} ),
\end{align*}
since $\nabla_{\Pi_i}\alpha=0, $ $\frac{\partial \alpha_i}{\partial
\alpha_i}= 1, \; \frac{\partial \alpha_j}{\partial \alpha_i}=
0, \; j\neq i, \;\; \frac{\partial H}{\partial \alpha_i}=0, $
and $\frac{\partial H}{\partial l_j}= -(\Pi_j-l_j)/\bar{I}_j
+\frac{l_j}{J_j}, \; i, j= 1,2,3. $

\begin{align*}
X_{H}(l)& =\{l,\; H\}_{-}=
-\Pi\cdot(\nabla_\Pi l\times\nabla_\Pi
H)+ \{l,\; H\}_{\mathbb{R}^3}\\
&= -\nabla_\Pi l\cdot(\nabla_\Pi H\times \Pi)+ \sum_{i=1}^3(\frac{\partial l}{\partial \alpha_i}
\frac{\partial H}{\partial l_i}- \frac{\partial
H}{\partial \alpha_i}\frac{\partial l}{\partial l_i})
=(0,0,0),
\end{align*}
since $\nabla_{\Pi_i} l=0,$ and $\frac{\partial l}{\partial \alpha_i}=
\frac{\partial H}{\partial \alpha_i}=0, \; i=1,2,3. $\\

Moreover, if we consider the rigid body-rotor system with a
control torque $u: T^\ast Q \to W $ acting on the rotors,
where the control subset $W\subset T^* Q $ is a fiber submanifold,
and assume that $u\in W $ is invariant under the cotangent lift of left
$\textmd{SO}(3)$-action $\Phi^{T*}$, and
the dynamical vector field of the regular point reducible
controlled rigid body-rotor system $(T^\ast Q,\textmd{SO}(3),\omega_Q,H,u)$
can be expressed by
\begin{align}
\tilde{X}= X_{(T^\ast Q,\textmd{SO}(3),\omega_Q,H,u)}
= X_H+ \textnormal{vlift}(u), \; \label{6.8}
\end{align}
where $\textnormal{vlift}(u)= \textnormal{vlift}(u)\cdot X_H $ is
the change of $X_H$ under the action of the control torque $u$.\\

Since the Hamiltonian $H(A,\Pi,\alpha,l)$ is invariant under
the cotangent lift $\Phi^{T*}$ of the left
$\textmd{SO}(3)$-action, for the point $\Pi_0 =\mu \in \mathfrak{so}^\ast(3)$ is
the regular value of $\mathbf{J}_Q$, we have the $R_p$-reduced Hamiltonian
$h_\mu(\Pi, \alpha,l):\mathcal{O}_\mu
\times\mathbb{R}^3\times\mathbb{R}^{3*} (\subset \mathfrak{so}^\ast
(3)\times\mathbb{R}^3\times\mathbb{R}^{3*})\to \mathbb{R},$ given by
$h_\mu(\Pi,\alpha,l)\cdot \pi_\mu=H(A,\Pi,\alpha,l)|_{\mathcal{O}_\mu
\times\mathbb{R}^3\times\mathbb{R}^{3*}}.$ Moreover, for the $R_p$-reduced
Hamiltonian $h_\mu(\Pi,\alpha,l): \mathcal{O}_\mu
\times\mathbb{R}^3\times\mathbb{R}^{3*} \to \mathbb{R}$, we have the
Hamiltonian vector field
$X_{h_\mu}(K_\mu)=\{K_\mu,h_\mu\}_{-}|_{\mathcal{O}_\mu
\times\mathbb{R}^3\times\mathbb{R}^{3*} }, $ where
$K_\mu(\Pi,\alpha,l): \mathcal{O}_\mu
\times\mathbb{R}^3\times\mathbb{R}^{3*} \to \mathbb{R}.$
Assume that $u\in
W \cap \mathbf{J}^{-1}_Q(\mu)$ and the $R_p$-reduced control torque $u_\mu:
\mathcal{O}_\mu \times\mathbb{R}^3\times\mathbb{R}^{3*} \to W_\mu
(\subset \mathcal{O}_\mu \times\mathbb{R}^3\times\mathbb{R}^{3*}) $ is
given by $u_\mu(\Pi,\alpha,l)\cdot \pi_\mu=
u(A,\Pi,\alpha,l)|_{\mathcal{O}_\mu
\times\mathbb{R}^3\times\mathbb{R}^{3*} }, $ where $\pi_\mu:
\mathbf{J}_Q^{-1}(\mu) \rightarrow \mathcal{O}_\mu
\times\mathbb{R}^3\times\mathbb{R}^{3*}, \; W_\mu= \pi_\mu(W\cap \mathbf{J}^{-1}_Q(\mu)). $
The $R_p$-reduced controlled rigid body-rotor
system is the 4-tuple $(\mathcal{O}_\mu \times \mathbb{R}^3 \times
\mathbb{R}^{3*},\tilde{\omega}_{\mathcal{O}_\mu \times \mathbb{R}^3
\times \mathbb{R}^{3*}}^{-},h_\mu,u_\mu), $ where
$\tilde{\omega}_{\mathcal{O}_\mu \times \mathbb{R}^3 \times
\mathbb{R}^{3*}}^{-}$ is the induced symplectic form
 on $\mathcal{O}_\mu
\times \mathbb{R}^3\times \mathbb{R}^{3*} ,$ such that
$X_{h_\mu}(K_\mu)=\tilde{\omega}_{\mathcal{O}_\mu \times \mathbb{R}^3 \times
\mathbb{R}^{3*}}^{-}(X_{K_\mu}, X_{h_\mu})
=\{K_\mu,h_\mu\}_{-}|_{\mathcal{O}_\mu
\times\mathbb{R}^3\times\mathbb{R}^{3*} }.$
Moreover, assume that the dynamical vector field of the $R_p$-reduced controlled
rigid body-rotor system $(\mathcal{O}_\mu \times \mathbb{R}^3 \times
\mathbb{R}^{3*},\tilde{\omega}_{\mathcal{O}_\mu \times \mathbb{R}^3 \times
\mathbb{R}^{3*}}^{-},h_\mu,u_\mu)$ is expressed by
\begin{align} X_{(\mathcal{O}_\mu \times \mathbb{R}^3 \times
\mathbb{R}^{3*},\tilde{\omega}_{\mathcal{O}_\mu \times \mathbb{R}^3 \times
\mathbb{R}^{3*}}^{-},h_\mu,u_\mu)} = X_{h_\mu} + \mbox{vlift}(u_\mu) , \; \label{6.9}
\end{align}
where $\mbox{vlift}(u_\mu)= \mbox{vlift}(u_\mu)X_{h_\mu} \in
T(\mathcal{O}_\mu \times\mathbb{R}^3 \times
\mathbb{R}^{3*}), $
is the change of $X_{h_\mu}$ under the action of the $R_p$-reduced control torque $u_\mu$.
The dynamical vector fields of the controlled
rigid body-rotor system and the $R_p$-reduced controlled
rigid body-rotor system satisfy the condition
\begin{equation}X_{(\mathcal{O}_\mu \times \mathbb{R}^3 \times
\mathbb{R}^{3*},\tilde{\omega}_{\mathcal{O}_\mu \times \mathbb{R}^3 \times
\mathbb{R}^{3*}}^{-},h_\mu,u_\mu)}\cdot \pi_\mu=T\pi_\mu\cdot X_{(T^\ast
Q,\textmd{SO}(3),\omega_Q,H,u)}\cdot i_\mu. \label{6.10}\end{equation}
See Marsden et al \cite{mawazh10} and Wang \cite{wa18}.\\

In the following we shall derive the geometric constraint
conditions of the $R_p$-reduced symplectic form
$\tilde{\omega}_{\mathcal{O}_\mu \times \mathbb{R}^3 \times
\mathbb{R}^{3*}}^{-}$ for the dynamical vector field
of the regular point reducible controlled rigid body-rotor system,
that is, Type I and Type II of
Hamilton-Jacobi equation for the $R_p$-reduced controlled rigid body-rotor system
$(\mathcal{O}_\mu \times\mathbb{R}^3\times\mathbb{R}^{3*} ,
\omega_{\mathcal{O}_\mu \times\mathbb{R}^3\times\mathbb{R}^{3*} }^{-},h_\mu, u_\mu).$
Assume that $\gamma: \textmd{SO}(3)\times \mathbb{R}^3 \rightarrow
T^*(\textmd{SO}(3)\times \mathbb{R}^3)$ is an one-form on
$\textmd{SO}(3)\times \mathbb{R}^3$,
$\gamma(A, \alpha)=(\gamma_1, \cdots, \gamma_{12})$,
and $\gamma$ is
closed with respect to $T\pi_{\textmd{SO}(3)\times \mathbb{R}^3}:
TT^* (\textmd{SO}(3)\times \mathbb{R}^3) \rightarrow
T(\textmd{SO}(3)\times \mathbb{R}^3). $ For $\mu \in \mathfrak{so}^\ast(3)$
the regular value of $\mathbf{J}_Q$,
$\textmd{Im}(\gamma)\subset \mathbf{J}_Q^{-1}(\mu), $ and it is
$\textmd{SO}(3)_\mu$-invariant, and $\bar{\gamma}=\pi_\mu(\gamma):
\textmd{SO}(3)\times \mathbb{R}^3 \rightarrow \mathcal{O}_\mu \times
\mathbb{R}^3 \times \mathbb{R}^{3*}$. Denote by
$\bar{\gamma}= (\bar{\gamma}_1, \cdots, \bar{\gamma}_9) \in \mathcal{O}_\mu \times
\mathbb{R}^3 \times \mathbb{R}^{3*}(\subset \mathfrak{so}^\ast(3)
\times \mathbb{R}^3 \times \mathbb{R}^{3*}), $ where
$\pi_\mu: \mathbf{J}_Q^{-1}(\mu) \rightarrow \mathcal{O}_\mu \times
\mathbb{R}^3 \times \mathbb{R}^{3*}. $ We choose that
$(\Pi,\alpha,l)\in
\mathcal{O}_\mu\times \mathbb{R}^3 \times \mathbb{R}^{3*}, $ and
$\Pi=(\Pi_1,\Pi_2,\Pi_3)=(\bar{\gamma}_1,\bar{\gamma}_2,\bar{\gamma}_3),
\; \alpha= (\alpha_1,\alpha_2,\alpha_3)
=(\bar{\gamma}_4,\bar{\gamma}_5,\bar{\gamma}_6), \; l=(l_1,l_2,l_3)
=(\bar{\gamma}_7,\bar{\gamma}_8,\bar{\gamma}_9), $ then $h_{\mu}
\cdot \bar{\gamma}: \textmd{SO}(3)\times \mathbb{R}^3 \rightarrow
\mathbb{R} $ is given by
\begin{align*} h_{\mu}(\Pi,\alpha,l) \cdot \bar{\gamma} & =
H(A,\Pi,\alpha,l) |_{\mathcal{O}_\mu \times
\mathbb{R}^3 \times \mathbb{R}^{3*}}\cdot \bar{\gamma}\\
& = \frac{1}{2}[\frac{(\bar{\gamma}_1-
\bar{\gamma}_7)^2}{\bar{I}_1}+\frac{(\bar{\gamma}_2-
\bar{\gamma}_8)^2}{\bar{I}_2} +\frac{(\bar{\gamma}_3-
\bar{\gamma}_9)^2}{\bar{I}_3}+\frac{\bar{\gamma}_7^2}{J_1}
+\frac{\bar{\gamma}_8^2}{J_2}+\frac{\bar{\gamma}_9^2}{J_3}],
\end{align*} and the vector field
\begin{align*}
& X_{h_{\mu}}(\Pi) \cdot \bar{\gamma}=\{\Pi,h_{\mu}
\}_{-}|_{\mathcal{O}_\mu \times \mathbb{R}^3 \times
\mathbb{R}^{3*}}\cdot \bar{\gamma}\\
& = -\Pi\cdot(\nabla_\Pi\Pi\times\nabla_\Pi (h_{\mu})) \cdot
\bar{\gamma}+ \{\Pi,h_{\mu} \}_{\mathbb{R}^3}|_{\mathcal{O}_\mu
\times \mathbb{R}^3 \times
\mathbb{R}^{3*}}\cdot \bar{\gamma}\\
& = -\nabla_\Pi\Pi\cdot(\nabla_\Pi (h_{\mu})\times \Pi)\cdot
\bar{\gamma} + \sum_{i=1}^3(\frac{\partial \Pi}{\partial \alpha_i}
\frac{\partial (h_\mu)}{\partial l_i}- \frac{\partial
(h_\mu)}{\partial
\alpha_i}\frac{\partial \Pi}{\partial l_i})\cdot \bar{\gamma} \\
& =(\Pi_1,\Pi_2,\Pi_3)\times (\frac{(\Pi_1- l_1)}{ \bar{I}_1},\;\;
\frac{(\Pi_2- l_2)}{ \bar{I}_2}, \;\; \frac{(\Pi_3-
l_3)}{\bar{I}_3})\cdot \bar{\gamma}\\
&= ( \frac{(\bar{I}_2-\bar{I}_3)\bar{\gamma}_2\bar{\gamma}_3-
\bar{I}_2\bar{\gamma}_2\bar{\gamma}_9 +
\bar{I}_3\bar{\gamma}_3\bar{\gamma}_8}{\bar{I}_2\bar{I}_3}, \;\;
\frac{(\bar{I}_3-\bar{I}_1)\bar{\gamma}_3\bar{\gamma}_1-
\bar{I}_3\bar{\gamma}_3\bar{\gamma}_7 +
\bar{I}_1\bar{\gamma}_1\bar{\gamma}_9}{\bar{I}_3\bar{I}_1}, \\
& \;\;\;\;\;\;
\frac{(\bar{I}_1-\bar{I}_2)\bar{\gamma}_1\bar{\gamma}_2-
\bar{I}_1\bar{\gamma}_1\bar{\gamma}_8 + \bar{I}_2\bar{\gamma}_2
\bar{\gamma}_7}{\bar{I}_1\bar{I}_2} ),
\end{align*}
since $\nabla_{\Pi_i}\Pi_i=1,\; \nabla_{\Pi_i}\Pi_j=0, \; i\neq j ,$ and $\nabla_{\Pi_j} (h_{\mu})=
(\Pi_j-l_j)/\bar{I}_j , $ and $\frac{\partial
\Pi}{\partial \alpha_i}= \frac{\partial (h_\mu)}{\partial
\alpha_i}=0, \; i, j=1,2,3. $

\begin{align*}
& X_{h_{\mu}}(\alpha) \cdot \bar{\gamma}=\{\alpha,h_{\mu}
\}_{-}|_{\mathcal{O}_\mu \times \mathbb{R}^3 \times
\mathbb{R}^{3*}}\cdot \bar{\gamma}\\
& =-\Pi\cdot(\nabla_\Pi\alpha \times\nabla_\Pi (h_{\mu})) \cdot
\bar{\gamma}+ \{\alpha,h_{\mu}\}_{\mathbb{R}^3}|_{\mathcal{O}_\mu
\times \mathbb{R}^3 \times
\mathbb{R}^{3*}}\cdot \bar{\gamma}\\
& =-\nabla_\Pi\alpha \cdot(\nabla_\Pi (h_{\mu})\times \Pi)\cdot
\bar{\gamma} + \sum_{i=1}^3(\frac{\partial \alpha}{\partial
\alpha_i} \frac{\partial (h_\mu)}{\partial l_i}- \frac{\partial
(h_\mu)}{\partial
\alpha_i}\frac{\partial \alpha}{\partial l_i})\cdot \bar{\gamma} \\
& = ( -\frac{(\bar{\gamma}_1- \bar{\gamma}_7)}{ \bar{I}_1}
+\frac{\bar{\gamma}_7}{J_1},\;\; -\frac{(\bar{\gamma}_2-
\bar{\gamma}_8)}{ \bar{I}_2} +\frac{\bar{\gamma}_8}{J_2},\;\;
-\frac{(\bar{\gamma}_3- \bar{\gamma}_9)}{\bar{I}_3}
+\frac{\bar{\gamma}_9}{J_3} ),
\end{align*}
since $\nabla_{\Pi_i}\alpha=0, $ $\frac{\partial \alpha_j}{\partial
\alpha_i}= 1,\; j=i, $ $\frac{\partial \alpha_j}{\partial \alpha_i}=
0, \; j\neq i, \;\; \frac{\partial (h_\mu)}{\partial \alpha_i}=0,$
and $\frac{\partial (h_\mu)}{\partial l_j}= -(\Pi_j-l_j)/\bar{I}_j
+\frac{l_j}{J_j}, \; i, j= 1,2,3 , $

\begin{align*}
&  X_{h_{\mu}}(l)
\cdot \bar{\gamma}=\{l,h_{\mu} \}_{-}|_{\mathcal{O}_\mu \times
\mathbb{R}^3 \times
\mathbb{R}^{3*}}\cdot \bar{\gamma}\\
& = -\Pi\cdot(\nabla_\Pi l \times\nabla_\Pi (h_{\mu})) \cdot
\bar{\gamma}+ \{l, h_{\mu} \}_{\mathbb{R}^3}|_{\mathcal{O}_\mu
\times \mathbb{R}^3 \times
\mathbb{R}^{3*}}\cdot \bar{\gamma}\\
& = -\nabla_\Pi l \cdot(\nabla_\Pi (h_{\mu})\times \Pi)\cdot
\bar{\gamma} + \sum_{i=1}^3(\frac{\partial l}{\partial \alpha_i}
\frac{\partial (h_\mu)}{\partial l_i}- \frac{\partial
(h_\mu)}{\partial \alpha_i}\frac{\partial l}{\partial l_i})\cdot
\bar{\gamma}=(0,0,0),
\end{align*}
since $\nabla_{\Pi_i} l=0,$ and $\frac{\partial l}{\partial \alpha_i}=
\frac{\partial (h_\mu)}{\partial \alpha_i}=0, \; i=1,2,3. $ \\

On the other hand, from the expressions of the dynamical vector field $\tilde{X}$
and Hamiltonian vector field $X_H$, we have that
\begin{align*}
\tilde{X}(\Pi, \alpha, l)^\gamma & =T\pi_{\textmd{SO}(3)\times \mathbb{R}^3 }\cdot \tilde{X}\cdot\gamma(\Pi, \alpha, l)\\
& =T\pi_{\textmd{SO}(3)\times \mathbb{R}^3 }\cdot (X_H+ \textnormal{vlift}(u))\cdot\gamma (\Pi, \alpha, l)\\
&=T\pi_{\textmd{SO}(3)\times \mathbb{R}^3 }\cdot X_H \cdot\gamma (\Pi, \alpha, l) = X_H\cdot\gamma(\Pi, \alpha, l),
\end{align*}
that is,
\begin{align*}
& \tilde{X}(\Pi)^\gamma  = X_H(\Pi)\cdot\gamma\\
& = ( \frac{(\bar{I}_2-\bar{I}_3)\gamma_5\gamma_6-
\bar{I}_2\gamma_5\gamma_{12} +
\bar{I}_3\gamma_6\gamma_{11}}{\bar{I}_2\bar{I}_3}, \;\;
\frac{(\bar{I}_3-\bar{I}_1)\gamma_6\gamma_4-
\bar{I}_3\gamma_6\gamma_{10} +
\bar{I}_1\gamma_4\gamma_{12}}{\bar{I}_3\bar{I}_1}, \\
& \;\;\;\;\;\;
\frac{(\bar{I}_1-\bar{I}_2)\gamma_4\gamma_5-
\bar{I}_1\gamma_4\gamma_{11} + \bar{I}_2\gamma_5\gamma_{10}}{\bar{I}_1\bar{I}_2} ),
\end{align*}
\begin{align*}
\tilde{X}( \alpha )^\gamma = X_H(\alpha)\cdot\gamma
 = ( -\frac{(\gamma_4- \gamma_{10})}{ \bar{I}_1}
+\frac{\gamma_{10}}{J_1},\;\; -\frac{(\gamma_5-
\gamma_{11})}{ \bar{I}_2} +\frac{\gamma_{11}}{J_2},\;\;
-\frac{(\gamma_6- \gamma_{12})}{\bar{I}_3}
+\frac{\gamma_{12}}{J_3} ),
\end{align*}
\begin{align*}
\tilde{X}( l )^\gamma = X_H(l)\cdot\gamma =(0,0,0).
\end{align*}
Since $\gamma$ is closed with respect to
$T\pi_{\textmd{SO}(3)\times \mathbb{R}^3}: TT^* (\textmd{SO}(3)\times \mathbb{R}^3)
\rightarrow T(\textmd{SO}(3)\times \mathbb{R}^3), $
then $\pi_{\textmd{SO}(3)\times \mathbb{R}^3}^*(\mathbf{d}\gamma)=0.$ We choose that
$(\gamma_4,\gamma_5,\gamma_6)=\Pi=(\Pi_1,\Pi_2,\Pi_3)=
(\bar{\gamma}_1,\bar{\gamma}_2,\bar{\gamma}_3), $
and $(\gamma_7,\gamma_8,\gamma_9)= \alpha= (\alpha_1,\alpha_2,\alpha_3)
=(\bar{\gamma}_4,\bar{\gamma}_5,\bar{\gamma}_6),
\; (\gamma_{10},\gamma_{11},\gamma_{12})= l=(l_1,l_2,l_3)
=(\bar{\gamma}_7,\bar{\gamma}_8,\bar{\gamma}_9). $
Hence
\begin{align*}
T\bar{\gamma}\cdot \tilde{X}(\Pi)^\gamma
= X_{h_\mu}(\Pi) \cdot \bar{\gamma}, \;\;\;\;\;\;
T\bar{\gamma}\cdot \tilde{X}(\alpha)^\gamma
= X_{h_\mu}(\alpha) \cdot \bar{\gamma}, \;\;\;\;\;\;
T\bar{\gamma}\cdot \tilde{X}(l)^\gamma
= X_{h_\mu}(l) \cdot \bar{\gamma}.
\end{align*}
Thus, the Type I of Hamilton-Jacobi equation for the
$R_p$-reduced controlled rigid body-rotor system
$(\mathcal{O}_\mu\times\mathbb{R}^3\times\mathbb{R}^{3*},
\omega_{\mathcal{O}_\mu\times\mathbb{R}^3\times\mathbb{R}^{3*}}^{-},h_\mu, u_\mu)$ holds.\\

Next, for $\mu \in \mathfrak{so}^\ast(3),$
the regular value of $\mathbf{J}_Q$, and
a $\textmd{SO}(3)_\mu$-invariant symplectic map
$\varepsilon: T^* (\textmd{SO}(3)\times \mathbb{R}^3)
\rightarrow T^* (\textmd{SO}(3)\times \mathbb{R}^3),$
assume that $\varepsilon(A,\Pi, \alpha, l)
=(\varepsilon_1,\cdots, \varepsilon_{12}),$
and $\varepsilon(\mathbf{J}_Q^{-1}(\mu))\subset \mathbf{J}_Q^{-1}(\mu). $
Denote by $\bar{\varepsilon}=\pi_\mu(\varepsilon): \mathbf{J}_Q^{-1}(\mu)\rightarrow
\mathcal{O}_\mu\times \mathbb{R}^3 \times \mathbb{R}^{3*}, $ and
$\bar{\varepsilon}=(\bar{\varepsilon}_1, \cdots, \bar{\varepsilon}_9) \in
\mathcal{O}_\mu\times \mathbb{R}^3 \times \mathbb{R}^{3*} (\subset \mathfrak{so}^\ast(3)
\times \mathbb{R}^3 \times \mathbb{R}^{3*}), $ and
$\lambda= \gamma \cdot \pi_{\textmd{SO}(3)\times \mathbb{R}^3}:
T^* (\textmd{SO}(3)\times \mathbb{R}^3)
\rightarrow T^* (\textmd{SO}(3)\times \mathbb{R}^3),$ and $\lambda(A,\Pi, \alpha, l)
=(\lambda_1,\cdots, \lambda_{12}),$
and $\bar{\lambda}=\pi_\mu(\lambda): \mathbf{J}_Q^{-1}(\mu) \rightarrow
\mathcal{O}_\mu\times \mathbb{R}^3 \times \mathbb{R}^{3*}, $ and
$\bar{\lambda}= (\bar{\lambda}_1, \cdots, \bar{\lambda}_9 ) \in
\mathcal{O}_\mu\times \mathbb{R}^3 \times \mathbb{R}^{3*}. $ We choose that
$(\Pi,\alpha,l)\in
\mathcal{O}_\mu\times \mathbb{R}^3 \times \mathbb{R}^{3*}, $ and
$\Pi=(\Pi_1,\Pi_2,\Pi_3)=(\bar{\varepsilon}_1,
\bar{\varepsilon}_2,\bar{\varepsilon}_3),$
and $\alpha= (\alpha_1,\alpha_2,\alpha_3)
=(\bar{\varepsilon}_4,\bar{\varepsilon}_5,\bar{\varepsilon}_6),$ and
$ l=(l_1,l_2,l_3)=(\bar{\varepsilon}_7,\bar{\varepsilon}_8,\bar{\varepsilon}_9),$
then $h_{\mu} \cdot \bar{\varepsilon}: T^*(\textmd{SO}(3)\times \mathbb{R}^3)
\rightarrow \mathbb{R} $ is given by
\begin{align*} h_{\mu}(\Pi,\alpha,l) \cdot \bar{\varepsilon}& =
H(A,\Pi, \alpha, l)|_{\mathcal{O}_\mu\times \mathbb{R}^3
 \times \mathbb{R}^{3*}} \cdot \bar{\varepsilon}\\
& = \frac{1}{2}[\frac{(\bar{\varepsilon}_1-
\bar{\varepsilon}_7)^2}{\bar{I}_1}+\frac{(\bar{\varepsilon}_2-
\bar{\varepsilon}_8)^2}{\bar{I}_2} +\frac{(\bar{\varepsilon}_3-
\bar{\varepsilon}_9)^2}{\bar{I}_3}+\frac{\bar{\varepsilon}_7^2}{J_1}
+\frac{\bar{\varepsilon}_8^2}{J_2}+\frac{\bar{\varepsilon}_9^2}{J_3}],
\end{align*}
and the vector field
\begin{align*}
& X_{h_{\mu}}(\Pi) \cdot \bar{\varepsilon}
=\{\Pi,h_{\mu} \}_{-}|_{\mathcal{O}_\mu\times \mathbb{R}^3 \times \mathbb{R}^{3*}} \cdot
\bar{\varepsilon}\\
& = -\Pi\cdot(\nabla_\Pi\Pi\times\nabla_\Pi (h_{\mu})) \cdot
\bar{\varepsilon}+ \{\Pi,h_{\mu} \}_{\mathbb{R}^3}|_{\mathcal{O}_\mu
\times \mathbb{R}^3 \times
\mathbb{R}^{3*}}\cdot \bar{\varepsilon}\\
&= ( \frac{(\bar{I}_2-\bar{I}_3)\bar{\varepsilon}_2\bar{\varepsilon}_3-
\bar{I}_2\bar{\varepsilon}_2\bar{\varepsilon}_9 +
\bar{I}_3\bar{\varepsilon}_3\bar{\varepsilon}_8}{\bar{I}_2\bar{I}_3}, \;\;
\frac{(\bar{I}_3-\bar{I}_1)\bar{\varepsilon}_3\bar{\varepsilon}_1-
\bar{I}_3\bar{\varepsilon}_3\bar{\varepsilon}_7 +
\bar{I}_1\bar{\varepsilon}_1\bar{\varepsilon}_9}{\bar{I}_3\bar{I}_1}, \\
& \;\;\;\;\;\;
\frac{(\bar{I}_1-\bar{I}_2)\bar{\varepsilon}_1\bar{\varepsilon}_2-
\bar{I}_1\bar{\varepsilon}_1\bar{\varepsilon}_8 + \bar{I}_2\bar{\varepsilon}_2
\bar{\varepsilon}_7}{\bar{I}_1\bar{I}_2} ),
\end{align*}
\begin{align*}
& X_{h_{\mu}}(\alpha) \cdot \bar{\varepsilon}=\{\alpha,h_{\mu}
\}_{-}|_{\mathcal{O}_\mu \times \mathbb{R}^3 \times
\mathbb{R}^{3*}}\cdot \bar{\varepsilon}\\
& =-\Pi\cdot(\nabla_\Pi\alpha \times\nabla_\Pi (h_{\mu})) \cdot
\bar{\varepsilon}+ \{\alpha,h_{\mu}\}_{\mathbb{R}^3}|_{\mathcal{O}_\mu
\times \mathbb{R}^3 \times
\mathbb{R}^{3*}}\cdot \bar{\varepsilon}\\
& = ( -\frac{(\bar{\varepsilon}_1- \bar{\varepsilon}_7)}{ \bar{I}_1}
+\frac{\bar{\varepsilon}_7}{J_1},\;\; -\frac{(\bar{\varepsilon}_2-
\bar{\varepsilon}_8)}{ \bar{I}_2} +\frac{\bar{\varepsilon}_8}{J_2},\;\;
-\frac{(\bar{\varepsilon}_3- \bar{\varepsilon}_9)}{\bar{I}_3}
+\frac{\bar{\varepsilon}_9}{J_3} ),
\end{align*}
\begin{align*}
&  X_{h_{\mu}}(l)
\cdot \bar{\varepsilon}=\{l,h_{\mu} \}_{-}|_{\mathcal{O}_\mu \times
\mathbb{R}^3 \times
\mathbb{R}^{3*}}\cdot \bar{\varepsilon}\\
& = -\Pi\cdot(\nabla_\Pi l \times\nabla_\Pi (h_{\mu})) \cdot
\bar{\varepsilon}+ \{l, h_{\mu} \}_{\mathbb{R}^3}|_{\mathcal{O}_\mu
\times \mathbb{R}^3 \times
\mathbb{R}^{3*}}\cdot \bar{\varepsilon}=(0,0,0).
\end{align*}

On the other hand, from the expressions of the dynamical vector field $\tilde{X}$
and Hamiltonian vector field $X_H$, we have that
\begin{align*}
\tilde{X}(\Pi, \alpha, l)^\varepsilon &
=T\pi_{\textmd{SO}(3)\times \mathbb{R}^3 }\cdot \tilde{X}\cdot\varepsilon(\Pi, \alpha, l)\\
& =T\pi_{\textmd{SO}(3)\times \mathbb{R}^3 }\cdot (X_H
+ \textnormal{vlift}(u))\cdot\varepsilon (\Pi, \alpha, l)\\
& =T\pi_{\textmd{SO}(3)\times \mathbb{R}^3 }
\cdot X_H \cdot\varepsilon (\Pi, \alpha, l)= X_H\cdot\varepsilon(\Pi, \alpha, l),
\end{align*}
that is,
\begin{align*}
& \tilde{X}(\Pi)^\varepsilon  = X_H(\Pi)\cdot\varepsilon\\
& = ( \frac{(\bar{I}_2-\bar{I}_3)\varepsilon_5\varepsilon_6-
\bar{I}_2\varepsilon_5\varepsilon_{12} +
\bar{I}_3\varepsilon_6\varepsilon_{11}}{\bar{I}_2\bar{I}_3}, \;\;
\frac{(\bar{I}_3-\bar{I}_1)\varepsilon_6\varepsilon_4-
\bar{I}_3\varepsilon_6\varepsilon_{10} +
\bar{I}_1\varepsilon_4\varepsilon_{12}}{\bar{I}_3\bar{I}_1}, \\
& \;\;\;\;\;\;
\frac{(\bar{I}_1-\bar{I}_2)\varepsilon_4\varepsilon_5-
\bar{I}_1\varepsilon_4\varepsilon_{11}
+ \bar{I}_2\varepsilon_5\varepsilon_{10}}{\bar{I}_1\bar{I}_2} ),
\end{align*}
\begin{align*}
\tilde{X}( \alpha )^\varepsilon = X_H(\alpha)\cdot\varepsilon
 = ( -\frac{(\varepsilon_4- \varepsilon_{10})}{ \bar{I}_1}
+\frac{\varepsilon_{10}}{J_1},\;\; -\frac{(\varepsilon_5-
\varepsilon_{11})}{ \bar{I}_2} +\frac{\varepsilon_{11}}{J_2},\;\;
-\frac{(\varepsilon_6- \varepsilon_{12})}{\bar{I}_3}
+\frac{\varepsilon_{12}}{J_3} ),
\end{align*}
\begin{align*}
\tilde{X}( l )^\varepsilon = X_H(l)\cdot\varepsilon =(0,0,0).
\end{align*}
Note that
\begin{align*}
T\bar{\gamma}\cdot \tilde{X}(\Pi)^\varepsilon
&= ( \frac{(\bar{I}_2-\bar{I}_3)\bar{\gamma}_2\bar{\gamma}_3-
\bar{I}_2\bar{\gamma}_2\bar{\gamma}_9 +
\bar{I}_3\bar{\gamma}_3\bar{\gamma}_8}{\bar{I}_2\bar{I}_3}, \;\;
\frac{(\bar{I}_3-\bar{I}_1)\bar{\gamma}_3\bar{\gamma}_1-
\bar{I}_3\bar{\gamma}_3\bar{\gamma}_7 +
\bar{I}_1\bar{\gamma}_1\bar{\gamma}_9}{\bar{I}_3\bar{I}_1}, \\
& \;\;\;\;\;\;
\frac{(\bar{I}_1-\bar{I}_2)\bar{\gamma}_1\bar{\gamma}_2-
\bar{I}_1\bar{\gamma}_1\bar{\gamma}_8 + \bar{I}_2\bar{\gamma}_2
\bar{\gamma}_7}{\bar{I}_1\bar{I}_2} ),
\end{align*}
\begin{align*}
T\bar{\gamma}\cdot \tilde{X}(\alpha)^\varepsilon
& = ( -\frac{(\bar{\gamma}_1- \bar{\gamma}_7)}{ \bar{I}_1}
+\frac{\bar{\gamma}_7}{J_1},\;\; -\frac{(\bar{\gamma}_2-
\bar{\gamma}_8)}{ \bar{I}_2} +\frac{\bar{\gamma}_8}{J_2},\;\;
-\frac{(\bar{\gamma}_3- \bar{\gamma}_9)}{\bar{I}_3}
+\frac{\bar{\gamma}_9}{J_3} ),
\end{align*}
\begin{align*}
T\bar{\gamma}\cdot \tilde{X}(l)^\varepsilon=(0,0,0),
\end{align*}
and $$T\bar{\lambda}\cdot \tilde{X} \cdot \varepsilon
=T\pi_\mu\cdot T\lambda \cdot (X_H+ \textnormal{vlift}(u))\cdot\varepsilon
=T\pi_\mu\cdot T\gamma \cdot T\pi_{\textmd{SO}(3)\times \mathbb{R}^3}\cdot (X_H+ \textnormal{vlift}(u))\cdot\varepsilon
=T\bar{\lambda}\cdot X_H \cdot \varepsilon,$$ that is,
\begin{align*}
& T\bar{\lambda}\cdot \tilde{X}(\Pi) \cdot \varepsilon
=T\bar{\lambda}\cdot X_H(\Pi) \cdot \varepsilon \\
&= ( \frac{(\bar{I}_2-\bar{I}_3)\bar{\lambda}_2\bar{\lambda}_3-
\bar{I}_2\bar{\lambda}_2\bar{\lambda}_9 +
\bar{I}_3\bar{\lambda}_3\bar{\lambda}_8}{\bar{I}_2\bar{I}_3}, \;\;
\frac{(\bar{I}_3-\bar{I}_1)\bar{\lambda}_3\bar{\lambda}_1-
\bar{I}_3\bar{\lambda}_3\bar{\lambda}_7 +
\bar{I}_1\bar{\lambda}_1\bar{\lambda}_9}{\bar{I}_3\bar{I}_1}, \\
& \;\;\;\;\;\;
\frac{(\bar{I}_1-\bar{I}_2)\bar{\lambda}_1\bar{\lambda}_2-
\bar{I}_1\bar{\lambda}_1\bar{\lambda}_8 + \bar{I}_2\bar{\lambda}_2
\bar{\lambda}_7}{\bar{I}_1\bar{I}_2} ),
\end{align*}
\begin{align*}
& T\bar{\lambda}\cdot \tilde{X}(\alpha) \cdot \varepsilon
=T\bar{\lambda}\cdot X_H(\alpha) \cdot \varepsilon \\
& = ( -\frac{(\bar{\lambda}_1- \bar{\lambda}_7)}{ \bar{I}_1}
+\frac{\bar{\lambda}_7}{J_1},\;\; -\frac{(\bar{\lambda}_2-
\bar{\lambda}_8)}{ \bar{I}_2} +\frac{\bar{\lambda}_8}{J_2},\;\;
-\frac{(\bar{\lambda}_3- \bar{\lambda}_9)}{\bar{I}_3}
+\frac{\bar{\lambda}_9}{J_3} ),
\end{align*}
\begin{align*}
T\bar{\lambda}\cdot \tilde{X}(l) \cdot \varepsilon
=T\bar{\lambda}\cdot X_H(l) \cdot \varepsilon=(0,0,0),
\end{align*}
Thus, when we choose that $(\Pi,\alpha,l)\in
\mathcal{O}_\mu\times \mathbb{R}^3 \times \mathbb{R}^{3*}, $ and
$(\varepsilon_4,\varepsilon_5,\varepsilon_6)=\Pi=(\Pi_1,\Pi_2,\Pi_3)
=(\bar{\gamma}_1,\bar{\gamma}_2,\bar{\gamma}_3)=
(\bar{\varepsilon}_1,\bar{\varepsilon}_2,\bar{\varepsilon}_3)=
(\bar{\lambda}_1,\bar{\lambda}_2,\bar{\lambda}_3), $
and $(\varepsilon_7,\varepsilon_8,\varepsilon_9)=\alpha= (\alpha_1,\alpha_2,\alpha_3)
=(\bar{\varepsilon}_4,\bar{\varepsilon}_5,\bar{\varepsilon}_6)
=(\bar{\lambda}_4,\bar{\lambda}_5,\bar{\lambda}_6),$ and
$ (\varepsilon_{10},\varepsilon_{11},\varepsilon_{12})=l
=(l_1,l_2,l_3)=(\bar{\varepsilon}_7,\bar{\varepsilon}_8,\bar{\varepsilon}_9)
=(\bar{\lambda}_7,\bar{\lambda}_8,\bar{\lambda}_9),$
we must have that
\begin{align*}
& T\bar{\gamma}\cdot \tilde{X}(\Pi)^\varepsilon=X_{h_{\mu}}(\Pi) \cdot \bar{\varepsilon}
=T\bar{\lambda}\cdot \tilde{X}(\Pi) \cdot \varepsilon, \\
& T\bar{\gamma}\cdot \tilde{X}(\alpha)^\varepsilon
=X_{h_{\mu}}(\alpha) \cdot \bar{\varepsilon}
=T\bar{\lambda}\cdot \tilde{X}(\alpha) \cdot \varepsilon, \\
& T\bar{\gamma}\cdot \tilde{X}(l)^\varepsilon=X_{h_{\mu}}(l) \cdot \bar{\varepsilon}
=T\bar{\lambda}\cdot \tilde{X}(l) \cdot \varepsilon.
\end{align*}
Since the map $\varepsilon: T^* (\textmd{SO}(3)\times \mathbb{R}^3)
\rightarrow T^* (\textmd{SO}(3)\times \mathbb{R}^3)$ is symplectic, then
$T\bar{\varepsilon}\cdot X_{h_{\mu} \cdot \bar{\varepsilon}}
=X_{h_{\mu}} \cdot \bar{\varepsilon}. $
Thus, in this case, we must have that
$\varepsilon$ and $\bar{\varepsilon} $ are the solution of the Type II of
Hamilton-Jacobi equation
$T\bar{\gamma}\cdot \tilde{X}^\varepsilon= X_{h_{\mu}}\cdot \bar{\varepsilon}, $
for the $R_p$-reduced controlled rigid body-rotor system
$(\mathcal{O}_\mu\times \mathbb{R}^3 \times \mathbb{R}^{3*},
\omega_{\mathcal{O}_\mu\times \mathbb{R}^3 \times \mathbb{R}^{3*}}^{-},
h_{\mu}, u_{\mu}) $, if and only if they satisfy
the equation $T\bar{\varepsilon}\cdot(X_{h_{\mu} \cdot \bar{\varepsilon}})
= T\bar{\lambda}\cdot \tilde{X} \cdot\varepsilon. $\\

To sum up the above discussion, we have the following proposition.
\begin{prop}
If the 5-tuple $(T^\ast Q,\textmd{SO}(3),\omega_Q,H,u), $ where $Q=
\textmd{SO}(3)\times \mathbb{R}^3, $ is a regular point reducible
rigid body-rotor system with the control torque $u$ acting on the rotors,
then for a point $\mu \in \mathfrak{so}^\ast(3)$, the regular
value of the momentum map $\mathbf{J}_Q: \textmd{SO}(3)\times
\mathfrak{so}^\ast(3) \times \mathbb{R}^3 \times \mathbb{R}^{3*} \to
\mathfrak{so}^\ast(3)$, the $R_p$-reduced controlled rigid body-rotor system is the 4-tuple
$(\mathcal{O}_\mu \times \mathbb{R}^3 \times
\mathbb{R}^{3*},\tilde{\omega}_{\mathcal{O}_\mu \times \mathbb{R}^3
\times \mathbb{R}^{3*}}^{-},h_\mu,u_\mu), $ where $\mathcal{O}_\mu
\subset \mathfrak{so}^\ast(3)$ is the co-adjoint orbit,
$\tilde{\omega}_{\mathcal{O}_\mu \times \mathbb{R}^3 \times
\mathbb{R}^{3*}}^{-}$ is the induced symplectic form
on $\mathcal{O}_\mu
\times \mathbb{R}^3\times \mathbb{R}^{3*} $, $h_\mu(\Pi,\alpha,l)\cdot \pi_\mu=
H(A,\Pi,\alpha,l)|_{\mathcal{O}_\mu
\times\mathbb{R}^3\times\mathbb{R}^{3*}}$, $u_\mu(\Pi,\alpha,l)\cdot \pi_\mu=
 u(A,\Pi,\alpha,l)|_{\mathcal{O}_\mu
\times\mathbb{R}^3\times\mathbb{R}^{3*}}$. Assume that $\gamma:
\textmd{SO}(3)\times \mathbb{R}^3 \rightarrow
T^*(\textmd{SO}(3)\times \mathbb{R}^3)$ is an one-form on
$\textmd{SO}(3)\times \mathbb{R}^3$,
and $\lambda=\gamma \cdot \pi_{(\textmd{SO}(3)\times \mathbb{R}^3)}:
T^* (\textmd{SO}(3)\times \mathbb{R}^3 )\rightarrow
T^* (\textmd{SO}(3)\times \mathbb{R}^3), $ and $\varepsilon:
T^* (\textmd{SO}(3) \times \mathbb{R}^3 )\rightarrow
T^* (\textmd{SO}(3)\times \mathbb{R}^3) $ is a
$\textmd{SO}(3)_\mu$-invariant symplectic map.
Denote
$\tilde{X}^\gamma = T\pi_{(\textmd{SO}(3)\times \mathbb{R}^3)}\cdot \tilde{X}\cdot \gamma$, and
$\tilde{X}^\varepsilon
= T\pi_{(\textmd{SO}(3)\times \mathbb{R}^3)}\cdot \tilde{X}\cdot \varepsilon$,
where $\tilde{X}=X_{(T^\ast Q,\textmd{SO}(3),\omega_Q,H,u)}$
is the dynamical vector field of the controlled rigid
body-rotor system $(T^\ast Q,\textmd{SO}(3),\omega_Q,H,u)$.
Moreover, assume that $\textmd{Im}(\gamma)\subset \mathbf{J}_Q^{-1}(\mu), $ and it is
$\textmd{SO}(3)_\mu$-invariant,
and $\varepsilon(\mathbf{J}_Q^{-1}(\mu))\subset \mathbf{J}_Q^{-1}(\mu). $
Denote $\bar{\gamma}=\pi_\mu(\gamma):
\textmd{SO}(3)\times \mathbb{R}^3 \rightarrow \mathcal{O}_\mu
\times \mathbb{R}^3\times \mathbb{R}^{3*}, $ and
$\bar{\lambda}=\pi_\mu(\lambda): T^* (\textmd{SO}(3)\times \mathbb{R}^3) \rightarrow
\mathcal{O}_\mu\times \mathbb{R}^3\times \mathbb{R}^{3*}, $ and
$\bar{\varepsilon}=\pi_\mu(\varepsilon): \mathbf{J}_Q^{-1}(\mu)\rightarrow
\mathcal{O}_\mu\times \mathbb{R}^3\times \mathbb{R}^{3*}. $
Then the following two assertions hold:\\
\noindent $(\mathbf{i})$
If the one-form $\gamma: \textmd{SO}(3)\times \mathbb{R}^3 \rightarrow
T^*(\textmd{SO}(3)\times \mathbb{R}^3) $ is closed with respect to
$T\pi_{(\textmd{SO}(3)\times \mathbb{R}^3)}:
TT^* (\textmd{SO}(3)\times \mathbb{R}^3) \rightarrow
T(\textmd{SO}(3)\times \mathbb{R}^3), $
then $\bar{\gamma}$ is a solution of the Type I of Hamilton-Jacobi equation
$T\bar{\gamma}\cdot \tilde{X}^\gamma= X_{h_\mu}\cdot \bar{\gamma}; $\\
\noindent $(\mathbf{ii})$
The $\varepsilon$ and $\bar{\varepsilon} $ satisfy the
Type II of Hamilton-Jacobi equation
$T\bar{\gamma}\cdot \tilde{X}^\varepsilon= X_{h_\mu}\cdot \bar{\varepsilon}, $
if and only if they satisfy
the equation $T\bar{\varepsilon}\cdot(X_{h_\mu \cdot \bar{\varepsilon}})
= T\bar{\lambda}\cdot \tilde{X}\cdot\varepsilon. $ \hskip 0.3cm $\blacksquare$
\end{prop}

When the rigid body does not carry any internal rotor, in this case
$Q=G= \textmd{SO}(3), $ and the rigid body is a regular point
reducible Hamiltonian system $(T^\ast
\textmd{SO}(3),\textmd{SO}(3),\omega, H)$, and hence it is also a
regular point reducible RCH system without the external force and
control. For a point $\mu \in \mathfrak{so}^\ast(3)$, the regular
value of the momentum map $\mathbf{J}: T^\ast \textmd{SO}(3)\to
\mathfrak{so}^\ast(3)$, the Marsden-Weinstein reduced rigid body system is
3-tuple $(\mathcal{O}_{\mu},
\omega_{\mathcal{O}_{\mu}},h_{\mathcal{O}_{\mu}})$, where
$\mathcal{O}_{\mu} \subset \mathfrak{so}^\ast(3)$ is the
co-adjoint orbit, $\omega_{\mathcal{O}_{\mu}}$ is the orbit
symplectic form on $\mathcal{O}_{\mu}$, which is induced by the
rigid body Lie-Poisson bracket on $\mathfrak{so}^\ast(3)$,
$h_{\mathcal{O}_{\mu}}(\Pi)\cdot \pi_{\mathcal{O}_{\mu}}
=H(A,v,\Pi)|_{\mathcal{O}_{\mu}}$. From the above Proposition 6.2
we can obtain the Proposition 5.3 in Wang \cite{wa17}, that is, we give the two
types of Lie-Poisson Hamilton-Jacobi equation for the Marsden-Weinstein
reduced rigid body system $(\mathcal{O}_{\mu},
\omega_{\mathcal{O}_{\mu}},h_{\mathcal{O}_{\mu}})$.
See Marsden and Ratiu \cite{mara99},
Ge and Marsden \cite{gema88}, and Wang \cite{wa17}.

\subsection {Hamilton-Jacobi Equations of Heavy Top with Internal Rotors }

In the following we regard the heavy top with two pairs of symmetric
internal rotors as a regular point reducible RCH system on the
generalization of Euclidean group $\textmd{SE}(3)\times
\mathbb{R}^2$, and give its two types of Hamilton-Jacobi equation
by calculation in detail. Note that our description of
the motion and the equations of heavy top with internal rotors in
this subsection follows some of the notations and conventions in
Marsden \cite{ma92}, Marsden and Ratiu \cite{mara99},  Marsden et
al \cite{mawazh10}, Wang \cite{wa17}.\\

We first describe a heavy top with two pairs of symmetric, "non-mass"
rotors which is called a carrier body, where "non-mass" means that
the mass of a rotor is very very small comparing with the mass of the heavy top.
We mount two pairs of rotors within the top so that each pair's
rotation axis is parallel to the first and the second principal axes
of the heavy top. Since the heavy top is moving in gravitational field,
and the rotor spins under the influence of a torque $u$
acting on the rotor. Then the motion of the controlled
heavy top-rotor system is just the rotation motion with drift,
in this case, the configuration space is
$Q=\textmd{SO}(3)\circledS \mathbb{R}^3\times V
\cong \textmd{SE}(3)\times V$, where $V=S^1\times S^1$, with the first
factor being the position of the carrier body and the second factor
being the angles of rotors. The corresponding phase space is the
cotangent bundle $T^*Q$ and locally,
$T^\ast Q \cong T^\ast \textmd{SE}(3)\times T^\ast V$,
where $T^\ast V = T^\ast (S^1\times S^1)\cong T^\ast \mathbb{R}^2$ locally,
with the canonical symplectic form $\omega_Q$.
By using the local left trivialization, locally,
$T^\ast \textmd{SE}(3)\cong \textmd{SE}(3)\times \mathfrak{se}^\ast(3)$
and $T^*\mathbb{R}^2 \cong \mathbb{R}^2 \times \mathbb{R}^{2*}$,
then we have that locally,
$T^*Q \cong \textmd{SE}(3)\times \mathfrak{se}^\ast(3)
\times \mathbb{R}^2 \times \mathbb{R}^{2*}$.
For convenience, in the following we denote uniformly that, locally,
$Q= \textmd{SE}(3)\times \mathbb{R}^2, $ and
$T^* Q= T^*(\textmd{SE}(3)\times \mathbb{R}^2)
\cong \textmd{SE}(3)\times \mathfrak{se}^\ast(3)
\times \mathbb{R}^2 \times \mathbb{R}^{2*}$.
Assume that Lie group $G=\textmd{SE}(3)$
acts freely and properly on $Q$ by the left translation
on the first factor $\textmd{SE}(3)$ and
the trivial action on the second factor $\mathbb{R}^2$.
Then the action of $\textmd{SE}(3)$ on the phase
space $T^\ast Q$ is by the cotangent lift of the left
$\textmd{SE}(3)$ action on $Q$, that is, $\Phi^{T*}:
\textmd{SE}(3)\times T^\ast Q \cong
\textmd{SE}(3)\times \textmd{SE}(3)\times \mathfrak{se}^\ast(3)
\times \mathbb{R}^2 \times \mathbb{R}^{2*}\to T^* Q \cong \textmd{SE}(3)\times
\mathfrak{se}^\ast(3)\times \mathbb{R}^2 \times \mathbb{R}^{2*},$ given
by $\Phi^{T*}((B,u)((A,v),(\Pi,w),\alpha,l))=((BA,u+Bv),(\Pi,w),\alpha,l)$,
for any $A,B\in \textmd{SO}(3), \; \Pi \in \mathfrak{so}^\ast(3), \;
u,v,w \in \mathbb{R}^3, \; \alpha \in \mathbb{R}^2,
\; l \in \mathbb{R}^{2*}$, and where $\textmd{SO}(3)$ acts
on $\mathbb{R}^3$ in the standard way. the action $\Phi^{T*}$ is also
free and proper. Moreover, assume that the action $\Phi^{T*}$ is
symplectic and admits an associated
$\operatorname{Ad}^\ast$-equivariant momentum map $\mathbf{J}_Q:
T^\ast Q \cong \textmd{SE}(3)\times \mathfrak{se}^\ast(3) \times
\mathbb{R}^2 \times \mathbb{R}^{2*} \to \mathfrak{se}^\ast(3)$
for the cotangent lift of the
left $\textmd{SE}(3)$ action. If $(\Pi,w) \in \mathfrak{se}^\ast(3)$
is a regular value of $\mathbf{J}_Q$, then the $R_p$-reduced
space $(T^\ast Q)_{(\Pi,w)}=
\mathbf{J}^{-1}_Q(\Pi,w)/\textmd{SE}(3)_{(\Pi,w)}$ is symplectically
diffeomorphic to the coadjoint orbit $\mathcal{O}_{(\Pi,w)} \times
\mathbb{R}^2 \times \mathbb{R}^{2*} \subset \mathfrak{se}^\ast(3)
\times \mathbb{R}^2 \times \mathbb{R}^{2*}$.\\

Let $I=diag(I_1,I_2,I_3)$ be the the matrix of inertia moment of the heavy top
in the body fixed frame, which is a principal body frame.
Let $J_k,\; k=1,2$ be the moments of inertia
of rotors around their rotation axes. Let $J_{ki},\; k=1,2,\;
i=1,2,3,$ be the moments of inertia of the $k$-th rotor with $k=1,2$
around the $i$-th principal axis with $i=1,2,3,$ respectively, and
denote $\bar{I}_k=I_k+J_{1k}+J_{2k}-J_{kk}, \; k=1,2$, and
$\bar{I}_3=I_3+J_{13}+J_{23}$. Let
$\Omega=(\Omega_1,\Omega_2,\Omega_3)$
be the angular velocity vector of heavy top-rotors
computed with respect to the axes fixed in the carrier
body and $(\Omega_1,\Omega_2,\Omega_3)\in \mathfrak{so}(3)$. Let
$\theta_k,\; k=1,2,$ be the relative angles of rotors and
$\dot{\theta}=(\dot{\theta_1},\dot{\theta_2})$ the relative angular
velocity vector of rotor about the principal axes with respect to
the carrier body fixed frame.
Let $g$ be the magnitude of the gravitational acceleration and
$h$ be the distance from the origin $O$ to the center of mass of the
system. For convenience, we assume the total mass of the system $m=1$.\\

Now, by the local left trivialization, locally,
$T\textmd{SE}(3)\cong \textmd{SE}(3)\times \mathfrak{se}(3)$
and $T\mathbb{R}^2 \cong \mathbb{R}^2 \times \mathbb{R}^{2}$,
then we have that locally,
$TQ \cong \textmd{SE}(3)\times \mathfrak{se}(3)
\times \mathbb{R}^2 \times \mathbb{R}^{2}$.
We consider the Lagrangian
$L(A,v,\Omega,\Gamma,\theta,\dot{\theta}):TQ\cong
\textmd{SE}(3)\times\mathfrak{se}(3)\times\mathbb{R}^2\times\mathbb{R}^2\to
\mathbb{R}$, which is the total kinetic energy of the heavy top plus
the total kinetic energy of the rotor minus the potential energy of
the system, given by
$$L(A,v,\Omega,\Gamma,\theta,\dot{\theta})
=\dfrac{1}{2}[\bar{I}_1\Omega_1^2+\bar{I}_2\Omega_2^2+\bar{I}_3\Omega_3^2
+J_1(\Omega_1+\dot{\theta}_1)^2+J_2(\Omega_2+\dot{\theta}_2)^2]
-gh\Gamma\cdot\chi,$$
where $(A,v)\in \textmd{SE}(3)$, $(\Omega,\Gamma)\in
\mathfrak{se}(3)$ and $\Omega=(\Omega_1,\Omega_2,\Omega_3)\in
\mathfrak{so}(3)$, $\theta=(\theta_1,\theta_2)\in \mathbb{R}^2$,
$\dot{\theta}=(\dot{\theta}_1,\dot{\theta}_2)\in \mathbb{R}^2$,
$\Gamma\in\mathbb{R}^3$, and the variable $\Gamma$ is regarded as a
parameter with respect to potential energy of the system. If
we introduce the conjugate angular momentum, which is given by
$$\Pi_k= \dfrac{\partial L}{\partial
\Omega_k}=\bar{I}_k\Omega_k+J_k(\Omega_k+\dot{\theta}_k),$$
$$\Pi_3=\dfrac{\partial L}{\partial
\Omega_3}=\bar{I}_3\Omega_3, \quad  l_k=\dfrac{\partial L}{\partial
\dot{\theta}_k}=J_k(\Omega_k+\dot{\theta}_k),\;\;\; k=1,2,$$ and by the
Legendre transformation with the parameter $\Gamma$, that is, $FL:TQ\cong
\textmd{SE}(3)\times\mathfrak{se}(3)\times\mathbb{R}^2\times\mathbb{R}^2\to
T^\ast Q\cong \textmd{SE}(3)\times
\mathfrak{se}^\ast(3)\times\mathbb{R}^2\times\mathbb{R}^{2*},\quad
(A,v,\Omega,\Gamma,\theta,\dot{\theta})\to(A,v,\Pi,\Gamma,\theta,l),$
where $\Pi=(\Pi_1,\Pi_2,\Pi_3)\in \mathfrak{so}^\ast(3)$,
$l=(l_1,l_2)\in \mathbb{R}^{2*}$, we have the Hamiltonian
$H(A,v,\Pi,\Gamma,\theta,l):T^\ast Q\cong \textmd{SE}(3)\times
\mathfrak{se}^\ast (3)\times\mathbb{R}^2\times\mathbb{R}^{2*}\to
\mathbb{R}$ given by

\begin{align*} &H(A,v,\Pi,\Gamma,\theta,l)=\Omega\cdot
\Pi+\dot{\theta}\cdot l-L(A,v,\Omega,\Gamma,\theta,\dot{\theta})\\
&=\bar{I}_1\Omega_1^2+J_1(\Omega_1^2+\Omega_1\dot{\theta}_1)
+\bar{I}_2\Omega_2^2+J_2(\Omega_2^2+\Omega_2\dot{\theta}_2)
+\bar{I}_3\Omega_3^2+J_1(\dot{\theta}_1\Omega_1+\dot{\theta}_1^2)\\
&\quad
  +J_2(\dot{\theta}_2\Omega_2+\dot{\theta}_2^2)
-\frac{1}{2}[\bar{I}_1\Omega_1^2+\bar{I}_2\Omega_2^2+\bar{I}_3\Omega_3^2
+J_1(\Omega_1+\dot{\theta}_1)^2+J_2(\Omega_2+\dot{\theta}_2)^2]
+gh\Gamma\cdot\chi\\
&=\frac{1}{2}[\frac{(\Pi_1-l_1)^2}{\bar{I}_1}+\frac{(\Pi_2-l_2)^2}{\bar{I}_2}
+\frac{\Pi_3^2}{\bar{I}_3}+\frac{l_1^2}{J_1}+\frac{l_2^2}{J_2}]
+gh\Gamma\cdot\chi.
\end{align*}
From the above expression of the Hamiltonian, we know that
$H(A,v,\Pi,\Gamma,\theta,l)$ is invariant under the cotangent lift of the left
$\textmd{SE}(3)$-action $\Phi^{T^*}:\textmd{SE}(3)\times T^\ast Q\to
T^\ast Q$. Moreover,
from the semidirect product Poisson bracket, see Marsden et al.
\cite{mamiorpera07}, we can get the heavy top Lie-Poisson bracket on
$\mathfrak{se}^\ast(3)$, that is, for $F,K: \mathfrak{se}^\ast(3)\to
\mathbb{R}, $ we have that
\begin{equation}
\{F,K\}_{-}(\Pi,\Gamma)=-\Pi\cdot(\nabla_\Pi F\times\nabla_\Pi
K)-\Gamma\cdot(\nabla_\Pi F\times \nabla_\Gamma K-\nabla_\Pi K\times
\nabla_\Gamma F). \label{6.11}
\end{equation}
Hence, from the heavy top Lie-Poisson
bracket on $\mathfrak{se}^\ast(3)$ and the Poisson bracket on
$T^\ast \mathbb{R}^2$, we can get the Poisson bracket
on $\mathfrak{se}^\ast(3)\times \mathbb{R}^2\times \mathbb{R}^{2*} $,
that is, for $F,K: \mathfrak{se}^\ast(3)\times \mathbb{R}^2\times
\mathbb{R}^{2*} \to \mathbb{R}, $ we have that
\begin{align}
& \{F,K\}_{-}(\Pi,\Gamma,\theta,l) \nonumber \\
& = -\Pi\cdot(\nabla_\Pi F\times\nabla_\Pi K)-\Gamma\cdot(\nabla_\Pi
F\times \nabla_\Gamma K-\nabla_\Pi K\times \nabla_\Gamma F)+
\{F,K\}_{\mathbb{R}^2}(\theta,l). \; \label{6.12}
\end{align}
Hence, the Hamiltonian vector fields of heavy top-rotor system are given by
\begin{align*}
& X_{H}(\Pi) = \{\Pi,\; H \}_{-}= -\Pi\cdot(\nabla_\Pi\Pi\times\nabla_\Pi
H) -\Gamma\cdot(\nabla_\Pi\Pi\times\nabla_\Gamma
H-\nabla_\Pi H \times\nabla_\Gamma\Pi) + \{\Pi,\; H \}_{\mathbb{R}^2}\\
& =(\Pi_1,\Pi_2,\Pi_3)\times (\frac{\Pi_1-l_1}{ \bar{I}_1}, \frac{\Pi_2-l_2}{
\bar{I}_2}, \frac{\Pi_3}{\bar{I}_3})
+gh(\Gamma_1,\Gamma_2,\Gamma_3)\times (\chi_1,\chi_2,\chi_3)
+ \sum_{k=1}^2(\frac{\partial \Pi}{\partial \theta_k}
\frac{\partial H}{\partial l_k}- \frac{\partial
H}{\partial \theta_k}\frac{\partial \Pi}{\partial l_k})\\
& = ( \frac{(\bar{I}_2-\bar{I}_3)\Pi_2\Pi_3+
\bar{I}_3\Pi_3l_2}{\bar{I}_2\bar{I}_3}+
gh(\Gamma_2\chi_3-\Gamma_3\chi_2), \;\;
\frac{(\bar{I}_3-\bar{I}_1)\Pi_3\Pi_1-
\bar{I}_3\Pi_3l_1}{\bar{I}_3\bar{I}_1}+gh(\Gamma_3\chi_1-\Gamma_1\chi_3), \\
& \;\;\;\;\;\;
\frac{(\bar{I}_1-\bar{I}_2)\Pi_1\Pi_2-
\bar{I}_1\Pi_1l_2 + \bar{I}_2\Pi_2l_1}{\bar{I}_1\bar{I}_2}
 + gh(\Gamma_1\chi_2-\Gamma_2\chi_1) ),
\end{align*}
since $\nabla_{\Pi_i}\Pi_i=1,\; \nabla_{\Pi_i}\Pi_j=0, \; i\neq j, \;  \nabla_{\Pi_i}\Gamma_j=\nabla_{\Gamma_i}\Pi_j=0$
and $\chi=(\chi_1,\chi_2,\chi_3), \; \nabla_{\Gamma_j}
H= gh\chi_j, \; \nabla_{\Pi_k} H= (\Pi_k-l_k)/\bar{I}_k ,\; \nabla_{\Pi_3}
H= \Pi_3/\bar{I}_3 , \; \frac{\partial
\Pi}{\partial \theta_k}= \frac{\partial H}{\partial
\theta_k}=0, \; i, j=1,2,3, \; k=1,2. $

\begin{align*}
 X_{H}(\Gamma)&
= \{\Gamma,\; H\}_{-} =-\Pi\cdot(\nabla_\Pi\Gamma\times\nabla_\Pi
H)-\Gamma\cdot(\nabla_\Pi\Gamma\times\nabla_\Gamma
H-\nabla_\Pi H \times\nabla_\Gamma\Gamma)+  \{\Gamma,\; H\}_{\mathbb{R}^2} \\
& =\nabla_\Gamma\Gamma\cdot (\Gamma\times\nabla_\Pi
H)+ \sum_{k=1}^2(\frac{\partial \Gamma}{\partial \theta_k}
\frac{\partial H}{\partial l_k}- \frac{\partial
H}{\partial \theta_k}\frac{\partial \Gamma}{\partial l_k})\\
& = (\Gamma_1,\Gamma_2,\Gamma_3)\times (\frac{(\Pi_1- l_1)}{ \bar{I}_1},\;\;
\frac{(\Pi_2- l_2)}{ \bar{I}_2}, \;\; \frac{\Pi_3}{\bar{I}_3}) \\
&= ( \frac{\bar{I}_2\Gamma_2\Pi_3-\bar{I}_3\Gamma_3\Pi_2+
\bar{I}_3\Pi_3l_2}{\bar{I}_2\bar{I}_3}, \;\;
\frac{\bar{I}_3\Gamma_3\Pi_1-\bar{I}_1\Gamma_1\Pi_3-
\bar{I}_3\Pi_3l_1}{\bar{I}_3\bar{I}_1}, \\
& \;\;\;\;\;\;
\frac{\bar{I}_1\Gamma_1\Pi_2-\bar{I}_2\Gamma_2\Pi_1-
\bar{I}_1\Pi_1l_2 + \bar{I}_2\Pi_2l_1}{\bar{I}_1\bar{I}_2} ),
\end{align*}
since $\nabla_{\Gamma_i}\Gamma_i =1, \;
\nabla_{\Gamma_i}\Gamma_j=0, \; i\neq j, \; \nabla_{\Pi_i}\Gamma_j =0, $ and
$\nabla_{\Pi_k} H= (\Pi_k-l_k)/\bar{I}_k ,\; \nabla_{\Pi_3}
H= \Pi_3/\bar{I}_3 , \; \frac{\partial
\Gamma_j}{\partial \theta_k}= \frac{\partial H}{\partial
\theta_k}=0, \; i, j= 1,2,3, \; k=1,2.$

\begin{align*}
 X_{H}(\theta)& = \{\theta,\; H \}_{-}= -\Pi\cdot(\nabla_\Pi\theta\times\nabla_\Pi
H) -\Gamma\cdot(\nabla_\Pi\theta\times\nabla_\Gamma
H-\nabla_\Pi H \times\nabla_\Gamma\theta) + \{\theta,\; H \}_{\mathbb{R}^2}\\
& =\sum_{k=1}^2(\frac{\partial \theta}{\partial \theta_k}
\frac{\partial H}{\partial l_k}- \frac{\partial
H}{\partial \theta_k}\frac{\partial \theta}{\partial l_k})
= ( -\frac{(\Pi_1- l_1)}{ \bar{I}_1}
+\frac{l_1}{J_1},\;\; -\frac{(\Pi_2-
l_2)}{ \bar{I}_2} +\frac{l_2}{J_2} ),
\end{align*}
since $\nabla_{\Pi_i}\theta=\nabla_{\Gamma_i}\theta=0, $
$\frac{\partial \theta_k}{\partial
\theta_k}= 1, \; \frac{\partial \theta_n}{\partial \theta_k}=
0, \; n\neq k, \;\; \frac{\partial H}{\partial \theta_k}=0, $
and $\frac{\partial H}{\partial l_k}= -(\Pi_k-l_k)/\bar{I}_k
+\frac{l_k}{J_k}, \; i= 1,2,3, \; k,n=1,2. $

\begin{align*}
 X_{H}(l)& = \{l,\; H \}_{-}= -\Pi\cdot(\nabla_\Pi l \times\nabla_\Pi
H) -\Gamma\cdot(\nabla_\Pi l \times\nabla_\Gamma
H-\nabla_\Pi H \times\nabla_\Gamma l) + \{l,\; H \}_{\mathbb{R}^2}\\
& = \sum_{k=1}^2(\frac{\partial l}{\partial \theta_k}
\frac{\partial H}{\partial l_k}- \frac{\partial
H}{\partial \theta_k}\frac{\partial l}{\partial l_k})=(0,0),
\end{align*}
since $\nabla_{\Pi_i} l=\nabla_{\Gamma_i} l=0,$ and $\frac{\partial l}{\partial \theta_k}=
\frac{\partial H}{\partial \theta_k}=0, \; i=1,2,3, \; k=1,2. $\\

Moreover, if we consider the heavy top-rotor system with a control torque $u: T^\ast
Q \to W $ acting on the rotors, where the control subset
$W\subset T^* Q $ is a fiber submanifold,
and assume that $u\in W $ is invariant under the cotangent lift
$\Phi^{T*}$ of the left $\textmd{SE}(3)$-action, and
the dynamical vector field of the regular point reducible
controlled heavy top-rotor system $(T^\ast Q,\textmd{SE}(3),\omega_Q,H,u)$
can be expressed by
\begin{align}
\tilde{X}= X_{(T^\ast Q,\textmd{SE}(3),\omega_Q,H,u)}
= X_H+ \textnormal{vlift}(u), \label{6.13}
\end{align}
where $\textnormal{vlift}(u)= \textnormal{vlift}(u)\cdot X_H $
is the change of $X_H$ under the action of the control torque $u$.\\

Since the Hamiltonian
$H(A,v,\Pi,\Gamma,\theta,l)$ is invariant under the cotangent lift $\Phi^{T^*}$ of the left
$\textmd{SE}(3)$-action, for the point $(\Pi_0,\Gamma_0)=(\mu,a)\in
\mathfrak{se}^\ast(3)$ is the regular value of $\mathbf{J}_Q$, we
have the $R_p$-reduced Hamiltonian
$h_{(\mu,a)}(\Pi,\Gamma,\theta,l):\mathcal{O}_{(\mu,a)}\times\mathbb{R}^2
\times \mathbb{R}^{2*} (\subset \mathfrak{se}^\ast (3)\times
\mathbb{R}^2\times \mathbb{R}^{2*}) \to \mathbb{R}$ given by
$h_{(\mu,a)}(\Pi,\Gamma,\theta,l)\cdot \pi_{(\mu,a)}
=H(A,v,\Pi,\Gamma,\theta,l)|_{\mathcal{O}_{(\mu,a)}\times
\mathbb{R}^2\times \mathbb{R}^{2*}}$, where $\pi_{(\mu,a)}:
\mathbf{J}_Q^{-1}(\mu,a) \rightarrow \mathcal{O}_{(\mu,a)}
\times\mathbb{R}^2\times\mathbb{R}^{2*}.$
Moreover, for the $R_p$-reduced
Hamiltonian $h_{(\mu,a)}(\Pi,\Gamma, \theta, l): \mathcal{O}_{(\mu,a)}\times
\mathbb{R}^2\times \mathbb{R}^{2*} \to \mathbb{R}$, we have the
Hamiltonian vector field
$$X_{h_{(\mu,a)}}(K_{(\mu,a)})=\{K_{(\mu,a)},h_{(\mu,a)}\}_{-}|_{\mathcal{O}_{(\mu,a)}
\times \mathbb{R}^2\times \mathbb{R}^{2*}}, $$ where
$K_{(\mu,a)}(\Pi,\Gamma, \theta, l): \mathcal{O}_{(\mu,a)}\times
\mathbb{R}^2\times \mathbb{R}^{2*} \to \mathbb{R}.$
Assume that $u\in W \cap
\mathbf{J}^{-1}_Q(\mu,a)$ and the  $R_p$-reduced control torque $u_{(\mu,a)}:
\mathcal{O}_{(\mu,a)} \times\mathbb{R}^2\times\mathbb{R}^{2*} \to W_{(\mu,a)}
(\subset \mathcal{O}_{(\mu,a)} \times\mathbb{R}^2\times\mathbb{R}^{2*}) $ is
given by $u_{(\mu,a)}(\Pi,\Gamma,\theta,l)\cdot \pi_{(\mu,a)}=
u(A,v,\Pi,\Gamma, \theta,l )|_{\mathcal{O}_{(\mu,a)}
\times\mathbb{R}^2\times\mathbb{R}^{2*} }, $ where
$ W_{(\mu,a)}= \pi_{(\mu,a)}(W\cap \mathbf{J}^{-1}_Q(\mu,a)). $
The $R_p$-reduced controlled heavy top-rotor
system is the 4-tuple $(\mathcal{O}_{(\mu,a)} \times \mathbb{R}^2 \times
\mathbb{R}^{2*},\tilde{\omega}_{\mathcal{O}_{(\mu,a)} \times \mathbb{R}^2
\times \mathbb{R}^{2*}}^{-}, \\ h_{(\mu,a)},u_{(\mu,a)}), $ where
$\tilde{\omega}_{\mathcal{O}_{(\mu,a)} \times \mathbb{R}^2 \times
\mathbb{R}^{2*}}^{-}$ is the induced symplectic form on  $\mathcal{O}_{(\mu,a)}
\times \mathbb{R}^2\times \mathbb{R}^{2*} ,$ such that Hamiltonian vector field
$$ X_{h_{(\mu,a)}}(K_{(\mu,a)})=\tilde{\omega}_{\mathcal{O}_{(\mu,a)} \times \mathbb{R}^2 \times
\mathbb{R}^{2*}}^{-}(X_{K_{(\mu,a)}}, X_{h_{(\mu,a)}})
=\{K_{(\mu,a)},h_{(\mu,a)}\}_{-}|_{\mathcal{O}_{(\mu,a)}
\times\mathbb{R}^2\times\mathbb{R}^{2*} }.$$
Moreover, assume that the dynamical vector field of the $R_p$-reduced controlled
heavy top-rotor system $(\mathcal{O}_{(\mu,a)} \times \mathbb{R}^2 \times
\mathbb{R}^{2*},\tilde{\omega}_{\mathcal{O}_{(\mu,a)} \times
\mathbb{R}^2 \times
\mathbb{R}^{2*}}^{-},h_{(\mu,a)},u_{(\mu,a)})$ is expressed by
\begin{align}
X_{(\mathcal{O}_{(\mu,a)} \times \mathbb{R}^2 \times
\mathbb{R}^{2*},\tilde{\omega}_{\mathcal{O}_{(\mu,a)} \times \mathbb{R}^2 \times
\mathbb{R}^{2*}}^{-},h_{(\mu,a)},u_{(\mu,a)})} = X_{h_{(\mu,a)}} +
\mbox{vlift}(u_{(\mu,a)}), \label{6.14}
\end{align}
where $\mbox{vlift}(u_{(\mu,a)})=
\mbox{vlift}(u_{(\mu,a)})X_{h_{(\mu,a)}} \in T(\mathcal{O}_{(\mu,a)}
\times \mathbb{R}^2 \times
\mathbb{R}^{2*}), $ is the change of $X_{h_{(\mu,a)}}$
under the action of the $R_p$-reduced control torque $u_{(\mu,a)}$.
The dynamical vector fields of the controlled
heavy top-rotor system and the $R_p$-reduced controlled
heavy top-rotor system satisfy the condition
\begin{equation}X_{(\mathcal{O}_{(\mu,a)} \times \mathbb{R}^2 \times
\mathbb{R}^{2*},\tilde{\omega}_{\mathcal{O}_{(\mu,a)} \times \mathbb{R}^2 \times
\mathbb{R}^{2*}}^{-}, h_{(\mu,a)}, u_{(\mu,a)})}\cdot \pi_{(\mu,a)}
=T\pi_{(\mu,a)}\cdot X_{(T^\ast Q,\textmd{SE}(3),\omega_Q,H,u)}\cdot i_{(\mu,a)}.
\label{6.15}\end{equation}
See Marsden et al \cite{mawazh10} and Wang \cite{wa18}.\\

In the following we shall derive the geometric constraint
conditions of the $R_p$-reduced symplectic form
$\tilde{\omega}_{\mathcal{O}_{(\mu,a)} \times \mathbb{R}^2 \times
\mathbb{R}^{2*}}^{-}$ for the
dynamical vector field of the regular point reducible controlled heavy top-rotor system,
that is, Type I and Type II of
Hamilton-Jacobi equation for the $R_p$-reduced controlled heavy top-rotor system
$(\mathcal{O}_{(\mu,a)}\times \mathbb{R}^2\times \mathbb{R}^{2*},
\omega^{-}_{\mathcal{O}_{(\mu,a)}\times \mathbb{R}^2\times \mathbb{R}^{2*}},
h_{(\mu,a)}, u_{(\mu,a)}) .$
Assume that $\gamma: \textmd{SE}(3)\times \mathbb{R}^2 \rightarrow
T^* (\textmd{SE}(3) \times \mathbb{R}^2) $ is an one-form on
$\textmd{SE}(3) \times \mathbb{R}^2 $, and
$\gamma(A,v,\theta)=(\gamma_1, \cdots, \gamma_{16})$,
and $\gamma$ is closed with respect to $T\pi_{\textmd{SE}(3)\times \mathbb{R}^2}:
TT^* (\textmd{SE}(3) \times \mathbb{R}^2) \rightarrow
T(\textmd{SE}(3) \times \mathbb{R}^2). $
For $(\mu,a) \in \mathfrak{se}^\ast(3),$ the regular value of $\mathbf{J}_Q$,
$\textmd{Im}(\gamma)\subset \mathbf{J}_Q^{-1}(\mu,a), $ and it is
$\textmd{SE}(3)_{(\mu,a)}$-invariant, and $\bar{\gamma}=\pi_{(\mu,a)}(\gamma):
\textmd{SE}(3)\times \mathbb{R}^2 \rightarrow \mathcal{O}_{(\mu,a)} \times
\mathbb{R}^2 \times \mathbb{R}^{2*}$. Denote by
$\bar{\gamma}= (\bar{\gamma}_1, \cdots, \bar{\gamma}_{10})
\in \mathcal{O}_{(\mu,a)} \times
\mathbb{R}^2 \times \mathbb{R}^{2*}(\subset \mathfrak{se}^\ast(3)
\times \mathbb{R}^2 \times \mathbb{R}^{2*}), $ where
$\pi_{(\mu,a)}: \mathbf{J}_Q^{-1}(\mu,a) \rightarrow \mathcal{O}_{(\mu,a)} \times
\mathbb{R}^2 \times \mathbb{R}^{2*}. $ We choose that
$(\Pi,\Gamma,\theta,l)\in
\mathcal{O}_{(\mu,a)}\times \mathbb{R}^2 \times \mathbb{R}^{2*}, $ and
$\Pi=(\Pi_1,\Pi_2,\Pi_3)=(\bar{\gamma}_1,\bar{\gamma}_2,\bar{\gamma}_3)$,
$\Gamma=(\Gamma_1,\Gamma_2,\Gamma_3)=(\bar{\gamma}_4,\bar{\gamma}_5,\bar{\gamma}_6)$,
$\theta=(\bar{\gamma}_7, \bar{\gamma}_8)$ and $l=(\bar{\gamma}_9,
\bar{\gamma}_{10}) $. Then $h_{(\mu,a)} \cdot \bar{\gamma}:
\textmd{SE}(3) \times \mathbb{R}^2 \rightarrow \mathbb{R} $ is given
by
\begin{align*} & h_{(\mu,a)}(\Pi,\Gamma,\theta,l) \cdot \bar{\gamma}=
H(A,v,\Pi,\Gamma,\theta,l) |_{\mathcal{O}_{(\mu,a)}
\times \mathbb{R}^2 \times \mathbb{R}^{2*}}
\cdot \bar{\gamma}\\ &
=\frac{1}{2}[\frac{(\bar{\gamma}_1-
\bar{\gamma}_9)^2}{\bar{I}_1}+\frac{(\bar{\gamma}_2-
\bar{\gamma}_{10})^2}{\bar{I}_2}
+\frac{\bar{\gamma}_3^2}{\bar{I}_3}+\frac{\bar{\gamma}_9^2}{J_1}
+\frac{\bar{\gamma}_{10}^2}{J_2}]
+ gh(\bar{\gamma}_4\cdot\chi_1+ \bar{\gamma}_5\cdot\chi_2+
\bar{\gamma}_6\cdot\chi_3), \end{align*} and the vector field
\begin{align*}
& X_{h_{(\mu,a)}}(\Pi) \cdot \bar{\gamma}=\{\Pi,\;
h_{(\mu,a)}\}_{-}|_{\mathcal{O}_{(\mu,a)}
\times \mathbb{R}^2\times \mathbb{R}^{2*}}\cdot \bar{\gamma}\\
& = -\Pi\cdot(\nabla_\Pi\Pi\times\nabla_\Pi (h_{(\mu,a)})) \cdot
\bar{\gamma}-\Gamma\cdot(\nabla_\Pi\Pi\times\nabla_\Gamma
(h_{(\mu,a)})-\nabla_\Pi
(h_{(\mu,a)}) \times\nabla_\Gamma\Pi)\cdot \bar{\gamma}\\
& \;\;\;\; + \{\Pi,\;
h_{(\mu,a)}\}_{\mathbb{R}^2}|_{\mathcal{O}_{(\mu,a)} \times
\mathbb{R}^2\times \mathbb{R}^{2*}}\cdot \bar{\gamma}\\
& = -\nabla_\Pi\Pi \cdot(\nabla_\Pi (h_{(\mu,a)}) \times \Pi)\cdot
\bar{\gamma}- \nabla_\Pi\Pi \cdot (\nabla_\Gamma
(h_{(\mu,a)}) \times \Gamma)\cdot \bar{\gamma}\\
& \;\;\;\; + \sum_{k=1}^2(\frac{\partial \Pi}{\partial \theta_k}
\frac{\partial (h_{(\mu,a)})}{\partial l_k}- \frac{\partial
(h_{(\mu,a)})}{\partial
\theta_k}\frac{\partial \Pi}{\partial l_k})\cdot \bar{\gamma}\\
&=(\Pi_1,\Pi_2,\Pi_3)\times (\frac{(\Pi_1- l_1)}{ \bar{I}_1},
\frac{(\Pi_2- l_2)}{ \bar{I}_2}, \frac{\Pi_3}{\bar{I}_3})\cdot
\bar{\gamma}
+gh(\Gamma_1,\Gamma_2,\Gamma_3)\times (\chi_1,\chi_2,\chi_3)\cdot \bar{\gamma}\\
&= ( \frac{(\bar{I}_2-\bar{I}_3)\bar{\gamma}_2\bar{\gamma}_3+
\bar{I}_3\bar{\gamma}_3\bar{\gamma}_{10}}{\bar{I}_2\bar{I}_3}+
gh(\bar{\gamma}_5\chi_3-\bar{\gamma}_6\chi_2),\\
& \;\;\;\;\;\; \frac{(\bar{I}_3-\bar{I}_1)\bar{\gamma}_3
\bar{\gamma}_1-
\bar{I}_3\bar{\gamma}_3\bar{\gamma}_9}{\bar{I}_3\bar{I}_1} +
gh(\bar{\gamma}_6\chi_1-\bar{\gamma}_4\chi_3),\\
& \;\;\;\;\;\;
\frac{(\bar{I}_1-\bar{I}_2)\bar{\gamma}_1\bar{\gamma}_2-\bar{I}_1\bar{\gamma}_1
\bar{\gamma}_{10}+
\bar{I}_2\bar{\gamma}_2\bar{\gamma}_9}{\bar{I}_1\bar{I}_2} +
gh(\bar{\gamma}_4\chi_2-\bar{\gamma}_5\chi_1) ),
\end{align*}
since $\nabla_{\Pi_i}\Pi_i=1,\; \nabla_{\Pi_i}\Pi_j=0, \; i\neq j, \;
\nabla_{\Pi_i}\Gamma_j=\nabla_{\Gamma_i}\Pi_j=0$
and $\chi=(\chi_1,\chi_2,\chi_3), \; \nabla_{\Gamma_j}
(h_{(\mu,a)})= gh\chi_j, \; \nabla_{\Pi_k} (h_{(\mu,a)}
)= (\Pi_k-l_k)/\bar{I}_k ,\; \nabla_{\Pi_3}
(h_{(\mu,a)})= \Pi_3/\bar{I}_3 , \; \frac{\partial
\Pi}{\partial \theta_k}= \frac{\partial (h_{(\mu,a)})}{\partial
\theta_k}=0, \; i, j=1,2,3, \; k=1,2. $

\begin{align*}
& X_{h_{(\mu,a)}}(\Gamma) \cdot \bar{\gamma}=\{\Gamma,\;
h_{(\mu,a)}\}_{-}|_{\mathcal{O}_{(\mu,a)} \times \mathbb{R}^2\times
\mathbb{R}^{2*}}\cdot \bar{\gamma}\\
& =-\Pi\cdot(\nabla_\Pi\Gamma\times\nabla_\Pi (h_{(\mu,a)})) \cdot
\bar{\gamma}-\Gamma\cdot(\nabla_\Pi\Gamma\times\nabla_\Gamma
(h_{(\mu,a)})-\nabla_\Pi
(h_{(\mu,a)}) \times\nabla_\Gamma\Gamma)\cdot \bar{\gamma}\\
& \;\;\;\; + \{\Gamma,\;
h_{(\mu,a)}\}_{\mathbb{R}^2}|_{\mathcal{O}_{(\mu,a)} \times
\mathbb{R}^2\times \mathbb{R}^{2*}}\cdot \bar{\gamma}\\
& =\nabla_\Gamma\Gamma\cdot(\Gamma\times\nabla_\Pi
(h_{(\mu,a)}))\cdot \bar{\gamma}+ \sum_{k=1}^2(\frac{\partial
\Gamma}{\partial \theta_k} \frac{\partial (h_{(\mu,a)})}{\partial
l_k}- \frac{\partial (h_{(\mu,a)})}{\partial
\theta_k}\frac{\partial \Gamma}{\partial l_k})\cdot \bar{\gamma}\\
&=(\Gamma_1,\Gamma_2,\Gamma_3)\times (\frac{(\Pi_1- l_1)}{
\bar{I}_1}, \frac{(\Pi_2- l_2)}{ \bar{I}_2},
\frac{\Pi_3}{\bar{I}_3})\cdot \bar{\gamma}\\
&= ( \frac{\bar{I}_2\bar{\gamma}_5\bar{\gamma}_3
-\bar{I}_3\bar{\gamma}_6\bar{\gamma}_2+
\bar{I}_3\bar{\gamma}_3\bar{\gamma}_{10}}{\bar{I}_2\bar{I}_3}, \;\;
\frac{\bar{I}_3\bar{\gamma}_6\bar{\gamma}_1
-\bar{I}_1\bar{\gamma}_4\bar{\gamma}_3-
\bar{I}_3\bar{\gamma}_3\bar{\gamma}_9}{\bar{I}_3\bar{I}_1}, \\
& \;\;\;\;\;\;
\frac{\bar{I}_1\bar{\gamma}_4\bar{\gamma}_2-\bar{I}_2\bar{\gamma}_5\bar{\gamma}_1-
\bar{I}_1\bar{\gamma}_1\bar{\gamma}_{10}
+ \bar{I}_2\bar{\gamma}_2\bar{\gamma}_9}{\bar{I}_1\bar{I}_2} ),
\end{align*}
since $\nabla_{\Gamma_i}\Gamma_i =1, \; \nabla_{\Gamma_i}\Gamma_j=0,
\; i\neq j, \; \nabla_{\Pi_i}\Gamma_j =0, $ and
$\nabla_{\Pi_k} (h_{(\mu,a)})= (\Pi_k-l_k)/\bar{I}_k ,\; \nabla_{\Pi_3}
(h_{(\mu,a)})= \Pi_3/\bar{I}_3 , \; \frac{\partial
\Gamma_j}{\partial \theta_k}= \frac{\partial (h_{(\mu,a)})}{\partial
\theta_k}=0, \; i, j= 1,2,3, \; k=1,2.$

\begin{align*}
& X_{h_{(\mu,a)}}(\theta) \cdot \bar{\gamma}=\{\theta,\;
h_{(\mu,a)}\}_{-}|_{\mathcal{O}_{(\mu,a)} \times \mathbb{R}^2\times
\mathbb{R}^{2*}}\cdot \bar{\gamma}\\
& =-\Pi\cdot(\nabla_\Pi \theta \times\nabla_\Pi (h_{(\mu,a)}))\cdot
\bar{\gamma}-\Gamma\cdot(\nabla_\Pi \theta \times\nabla_\Gamma
(h_{(\mu,a)})-\nabla_\Pi
(h_{(\mu,a)}) \times\nabla_\Gamma \theta)\cdot \bar{\gamma}\\
& \;\;\;\; + \{\theta,\;
h_{(\mu,a)}\}_{\mathbb{R}^2}|_{\mathcal{O}_{(\mu,a)} \times
\mathbb{R}^2\times \mathbb{R}^{2*}}\cdot \bar{\gamma}\\
& = \sum_{k=1}^2(\frac{\partial \theta}{\partial \theta_k}
\frac{\partial (h_{(\mu,a)})}{\partial l_k}- \frac{\partial
(h_{(\mu,a)})}{\partial
\theta_k}\frac{\partial \theta}{\partial l_k})\cdot \bar{\gamma}\\
& = (-\frac{(\bar{\gamma}_1- \bar{\gamma}_9)}{ \bar{I}_1}+
\frac{\bar{\gamma}_9}{J_1}, \; -\frac{(\bar{\gamma}_2-
\bar{\gamma}_{10})}{ \bar{I}_2}+ \frac{\bar{\gamma}_{10}}{J_2}),
\end{align*}
since $\nabla_{\Pi_i}\theta=\nabla_{\Gamma_i}\theta=0, $
$\frac{\partial \theta_k}{\partial
\theta_k}= 1, \; \frac{\partial \theta_n}{\partial \theta_k}=
0, \; n\neq k, \;\; \frac{\partial (h_{(\mu,a)})}{\partial \theta_k}=0, $
and $\frac{\partial (h_{(\mu,a)})}{\partial l_k}= -(\Pi_k-l_k)/\bar{I}_k
+\frac{l_k}{J_k}, \; i= 1,2,3, \; k,n=1,2. $

\begin{align*}
& X_{h_{(\mu,a)}}(l) \cdot \bar{\gamma}=\{l,\;
h_{(\mu,a)}\}_{-}|_{\mathcal{O}_{(\mu,a)} \times \mathbb{R}^2\times
\mathbb{R}^{2*}}\cdot \bar{\gamma}\\
& =-\Pi\cdot(\nabla_\Pi l \times\nabla_\Pi (h_{(\mu,a)})) \cdot
\bar{\gamma}-\Gamma\cdot(\nabla_\Pi l \times\nabla_\Gamma
(h_{(\mu,a)})-\nabla_\Pi
(h_{(\mu,a)}) \times\nabla_\Gamma l)\cdot \bar{\gamma}\\
& \;\;\;\; + \{l,\;
h_{(\mu,a)}\}_{\mathbb{R}^2}|_{\mathcal{O}_{(\mu,a)} \times
\mathbb{R}^2\times \mathbb{R}^{2*}}\cdot \bar{\gamma}\\
& = \sum_{k=1}^2(\frac{\partial l}{\partial \theta_k} \frac{\partial
(h_{(\mu,a)})}{\partial l_k}- \frac{\partial (h_{(\mu,a)})}{\partial
\theta_k}\frac{\partial l}{\partial l_k})\cdot \bar{\gamma}=(0,0),
\end{align*}
since $\nabla_{\Pi_i} l=\nabla_{\Gamma_i} l=0,$
and $\frac{\partial l}{\partial \theta_k}=
\frac{\partial (h_{(\mu,a)})}{\partial \theta_k}=0, \; i=1,2,3, \; k=1,2. $\\

On the other hand, from the expressions of the dynamical vector field $\tilde{X}$
and Hamiltonian vector field $X_H$, we have that
\begin{align*}
\tilde{X}(\Pi, \Gamma, \theta, l)^\gamma &
=T\pi_{\textmd{SE}(3)\times \mathbb{R}^2 }\cdot \tilde{X}\cdot\gamma(\Pi, \Gamma, \theta, l)\\
& =T\pi_{\textmd{SE}(3)\times \mathbb{R}^2 }\cdot
(X_H+ \textnormal{vlift}(u))\cdot\gamma (\Pi, \Gamma, \theta, l)\\
&=T\pi_{\textmd{SE}(3)\times \mathbb{R}^2 }\cdot
X_H \cdot\gamma (\Pi, \Gamma, \theta, l) = X_H\cdot\gamma(\Pi, \Gamma, \theta, l),
\end{align*}
that is,
\begin{align*}
\tilde{X}(\Pi)^\gamma & = X_H(\Pi)\cdot\gamma \\
& = ( \frac{(\bar{I}_2-\bar{I}_3)\gamma_8\gamma_9+
\bar{I}_3\gamma_9\gamma_{16}}{\bar{I}_2\bar{I}_3}+
gh(\gamma_{11}\chi_3-\gamma_{12}\chi_2),\\
& \;\;\;\;\;\;
\frac{(\bar{I}_3-\bar{I}_1)\gamma_9\gamma_7-
\bar{I}_3\gamma_9\gamma_{15}}{\bar{I}_3\bar{I}_1}
+gh(\gamma_{12}\chi_1-\gamma_{10}\chi_3), \\
& \;\;\;\;\;\;
\frac{(\bar{I}_1-\bar{I}_2)\gamma_7\gamma_8-
\bar{I}_1\gamma_7\gamma_{16}
+ \bar{I}_2\gamma_8\gamma_{15}}{\bar{I}_1\bar{I}_2}
 + gh(\gamma_{10}\chi_2-\gamma_{11}\chi_1) ),
\end{align*}
\begin{align*}
\tilde{X}(\Gamma)^\gamma & = X_H(\Gamma)\cdot\gamma \\
&= ( \frac{\bar{I}_2\gamma_{11}\gamma_9-\bar{I}_3\gamma_{12}\gamma_8+
\bar{I}_3\gamma_9\gamma_{16}}{\bar{I}_2\bar{I}_3}, \;\;
\frac{\bar{I}_3\gamma_{12}\gamma_7-\bar{I}_1\gamma_{10}\gamma_9-
\bar{I}_3\gamma_9\gamma_{15}}{\bar{I}_3\bar{I}_1}, \\
& \;\;\;\;\;\;
\frac{\bar{I}_1\gamma_{10}\gamma_8-\bar{I}_2\gamma_{11}\gamma_7-
\bar{I}_1\gamma_7\gamma_{16}
+ \bar{I}_2\gamma_8\gamma_{15}}{\bar{I}_1\bar{I}_2} ),
\end{align*}
\begin{align*}
\tilde{X}(\theta)^\gamma & = X_H(\theta)\cdot\gamma
 = ( -\frac{(\gamma_7- \gamma_{15})}{ \bar{I}_1}
+\frac{\gamma_{15}}{J_1},\;\; -\frac{(\gamma_8- \gamma_{16})}{ \bar{I}_2}
+\frac{\gamma_{16}}{J_2} ),
\end{align*}
\begin{align*}
\tilde{X}(l)^\gamma & = X_H(l)\cdot\gamma=(0,0).
\end{align*}
Since $\gamma$ is closed with respect to
$T\pi_{(\textmd{SE}(3) \times \mathbb{R}^2)}:
TT^* (\textmd{SE}(3) \times \mathbb{R}^2)
\rightarrow T(\textmd{SE}(3) \times \mathbb{R}^2), $
then $\pi_{(\textmd{SE}(3) \times \mathbb{R}^2)}^*(\mathbf{d}\gamma)=0.$ We choose that
$(\gamma_7,\gamma_8,\gamma_9)=\Pi=(\Pi_1,\Pi_2,\Pi_3)=
(\bar{\gamma}_1,\bar{\gamma}_2,\bar{\gamma}_3), \;
(\gamma_{10},\gamma_{11},\gamma_{12})=\Gamma=(\Gamma_1,\Gamma_2,\Gamma_3)=
(\bar{\gamma}_4,\bar{\gamma}_5,\bar{\gamma}_6), $
and $(\gamma_{13},\gamma_{14})= \theta= (\theta_1,\theta_2)
=(\bar{\gamma}_7,\bar{\gamma}_8),
\; (\gamma_{15},\gamma_{16})= l=(l_1,l_2)
=(\bar{\gamma}_9,\bar{\gamma}_{10}). $
Hence
\begin{align*}
& T\bar{\gamma}\cdot \tilde{X}(\Pi)^\gamma
= X_{h_{(\mu,a)}}(\Pi) \cdot \bar{\gamma}, \;\;\;\;\;\;
T\bar{\gamma}\cdot \tilde{X}(\Gamma)^\gamma
= X_{h_{(\mu,a)}}(\Gamma) \cdot \bar{\gamma}, \\
& T\bar{\gamma}\cdot \tilde{X}(\theta)^\gamma
= X_{h_{(\mu,a)}}(\theta) \cdot \bar{\gamma}, \;\;\;\;\;\;
T\bar{\gamma}\cdot \tilde{X}(l)^\gamma
= X_{h_{(\mu,a)}}(l) \cdot \bar{\gamma}.
\end{align*}
Thus, the Type I of Hamilton-Jacobi equation for the
$R_p$-reduced controlled heavy top-rotor system
$(\mathcal{O}_{(\mu,a)}\times \mathbb{R}^2\times \mathbb{R}^{2*},
\omega^{-}_{\mathcal{O}_{(\mu,a)}\times \mathbb{R}^2\times \mathbb{R}^{2*}},
h_{(\mu,a)}, u_{(\mu,a)})$ holds.\\

Next, for $(\mu,a) \in \mathfrak{se}^\ast(3),$
the regular value of $\mathbf{J}_Q$, and a $\textmd{SE}(3)_{(\mu,a)}$-invariant symplectic map
$\varepsilon: T^* (\textmd{SE}(3)\times \mathbb{R}^2) \rightarrow
T^*(\textmd{SE}(3)\times \mathbb{R}^2),$
assume that $\varepsilon(A,v,\Pi,\Gamma, \theta, l) =(\varepsilon_1,\cdots,
\varepsilon_{16}),$ and
$\varepsilon(\mathbf{J}_Q^{-1}((\mu,a)))\subset
\mathbf{J}_Q^{-1}((\mu,a)). $ Denote by
$\bar{\varepsilon}=\pi_{(\mu,a)}(\varepsilon):
\mathbf{J}_Q^{-1}((\mu,a)) \rightarrow \mathcal{O}_{(\mu,a)}
\times \mathbb{R}^2\times \mathbb{R}^{2*}, $
and $\bar{\varepsilon}=(\bar{\varepsilon}_1,\cdots,\bar{\varepsilon}_{10})
\in \mathcal{O}_{(\mu,a)}\times \mathbb{R}^2\times \mathbb{R}^{2*}
(\subset \mathfrak{se}^\ast(3)\times \mathbb{R}^2\times \mathbb{R}^{2*}), $
and $\lambda= \gamma \cdot \pi_{\textmd{SE}(3)\times \mathbb{R}^2}:
T^* (\textmd{SE}(3)\times \mathbb{R}^2)
\rightarrow T^* (\textmd{SE}(3)\times \mathbb{R}^2),$
and $\lambda(A,v,\Gamma, \Pi, \theta, l)
=(\lambda_1,\cdots, \lambda_{16}),$ and
$\bar{\lambda}=\pi_{(\mu,a)}(\lambda): T^* (\textmd{SE}(3)\times \mathbb{R}^2)
\rightarrow \mathcal{O}_{(\mu,a)}\times \mathbb{R}^2\times \mathbb{R}^{2*}, $ and
$\bar{\lambda}=
(\bar{\lambda}_1,\cdots,\bar{\lambda}_{10}) \in
\mathcal{O}_{(\mu,a)}\times \mathbb{R}^2\times \mathbb{R}^{2*}
(\subset \mathfrak{se}^\ast(3)\times \mathbb{R}^2\times \mathbb{R}^{2*}). $ We choose that
$(\Pi,\Gamma,\theta,l)\in
\mathcal{O}_{(\mu,a)}\times \mathbb{R}^2 \times \mathbb{R}^{2*}, $ and
$\Pi=(\Pi_1,\Pi_2,\Pi_3)
=(\bar{\varepsilon}_1,\bar{\varepsilon}_2,\bar{\varepsilon}_3),
\; \Gamma = (\Gamma_1,\Gamma_2,\Gamma_3)=
(\bar{\varepsilon}_4,\bar{\varepsilon}_5,\bar{\varepsilon}_6), $
$ \theta = (\theta_1,\theta_2)=(\bar{\varepsilon}_7, \bar{\varepsilon}_8)$ and
$l=(l_1,l_2)=(\bar{\varepsilon}_9, \bar{\varepsilon}_{10}) $.
Then $h_{(\mu,a)} \cdot \bar{\varepsilon}:
T^*(\textmd{SE}(3)\times \mathbb{R}^2) \rightarrow \mathbb{R} $ is given by
\begin{align*}
& h_{(\mu,a)}(\Pi,\Gamma,\theta,l) \cdot \bar{\varepsilon}=
H(A,v,\Pi,\Gamma, \theta, l)|_{\mathcal{O}_{(\mu,a)}\times
\mathbb{R}^2\times \mathbb{R}^{2*}} \cdot \bar{\varepsilon}\\
& =\frac{1}{2}[\frac{(\bar{\varepsilon}_1-
\bar{\varepsilon}_9)^2}{\bar{I}_1}+\frac{(\bar{\varepsilon}_2-
\bar{\varepsilon}_{10})^2}{\bar{I}_2}
+\frac{\bar{\varepsilon}_3^2}{\bar{I}_3}
+\frac{\bar{\varepsilon}_9^2}{J_1}+\frac{\bar{\varepsilon}_{10}^2}{J_2}]
+ gh(\bar{\varepsilon}_4\cdot\chi_1+ \bar{\varepsilon}_5\cdot\chi_2+
\bar{\varepsilon}_6\cdot\chi_3),
\end{align*} and the vector field
\begin{align*}
& X_{h_{(\mu,a)}}(\Pi) \cdot \bar{\varepsilon}
= \{\Pi,h_{(\mu,a)}\}_{-}|_{\mathcal{O}_{(\mu,a)}
\times \mathbb{R}^2\times \mathbb{R}^{2*}}\cdot \bar{\varepsilon}\\
&=(\Pi_1,\Pi_2,\Pi_3)\times (\frac{(\Pi_1- l_1)}{ \bar{I}_1},
\frac{(\Pi_2- l_2)}{ \bar{I}_2}, \frac{\Pi_3}{\bar{I}_3})\cdot
\bar{\varepsilon}
+gh(\Gamma_1,\Gamma_2,\Gamma_3)\times (\chi_1,\chi_2,\chi_3)\cdot \bar{\varepsilon}\\
&= ( \frac{(\bar{I}_2-\bar{I}_3)\bar{\varepsilon}_2\bar{\varepsilon}_3+
\bar{I}_3\bar{\varepsilon}_3\bar{\varepsilon}_{10}}{\bar{I}_2\bar{I}_3}+
gh(\bar{\varepsilon}_5\chi_3-\bar{\varepsilon}_6\chi_2), \\
& \;\;\;\;\;\; \frac{(\bar{I}_3-\bar{I}_1)\bar{\varepsilon}_3
\bar{\varepsilon}_1-
\bar{I}_3\bar{\varepsilon}_3\bar{\varepsilon}_9}{\bar{I}_3\bar{I}_1} +
gh(\bar{\varepsilon}_6\chi_1-\bar{\varepsilon}_4\chi_3),\\
& \;\;\;\;\;\;
\frac{(\bar{I}_1-\bar{I}_2)\bar{\varepsilon}_1\bar{\varepsilon}_2
-\bar{I}_1\bar{\varepsilon}_1 \bar{\varepsilon}_{10}+
\bar{I}_2\bar{\varepsilon}_2\bar{\varepsilon}_9}{\bar{I}_1\bar{I}_2} +
gh(\bar{\varepsilon}_4\chi_2-\bar{\varepsilon}_5\chi_1) ),
\end{align*}
\begin{align*}
& X_{h_{(\mu,a)}}(\Gamma) \cdot
\bar{\varepsilon}
= \{\Gamma,h_{(\mu,a)}\}_{-}|_{\mathcal{O}_{(\mu,a)}
\times \mathbb{R}^2\times \mathbb{R}^{2*}}\cdot \bar{\varepsilon}
=(\Gamma_1,\Gamma_2,\Gamma_3)\times (\frac{(\Pi_1- l_1)}{
\bar{I}_1}, \frac{(\Pi_2- l_2)}{ \bar{I}_2},
\frac{\Pi_3}{\bar{I}_3})\cdot \bar{\varepsilon}\\
&= ( \frac{\bar{I}_2\bar{\varepsilon}_5\bar{\varepsilon}_3
-\bar{I}_3\bar{\varepsilon}_6\bar{\varepsilon}_2+
\bar{I}_3\bar{\varepsilon}_3\bar{\varepsilon}_{10}}{\bar{I}_2\bar{I}_3}, \;\;
\frac{\bar{I}_3\bar{\varepsilon}_6\bar{\varepsilon}_1
-\bar{I}_1\bar{\varepsilon}_4\bar{\varepsilon}_3-
\bar{I}_3\bar{\varepsilon}_3\bar{\varepsilon}_9}{\bar{I}_3\bar{I}_1}, \\
& \;\;\;\;\;\;
\frac{\bar{I}_1\bar{\varepsilon}_4\bar{\varepsilon}_2
-\bar{I}_2\bar{\varepsilon}_5\bar{\varepsilon}_1-
\bar{I}_1\bar{\varepsilon}_1\bar{\varepsilon}_{10} + \bar{I}_2\bar{\varepsilon}_2\bar{\varepsilon}_9}{\bar{I}_1\bar{I}_2} ),
\end{align*}
\begin{align*}
& X_{h_{(\mu,a)}}(\theta) \cdot \bar{\varepsilon}=\{\theta,\;
h_{(\mu,a)}\}_{-}|_{\mathcal{O}_{(\mu,a)} \times \mathbb{R}^2\times
\mathbb{R}^{2*}}\cdot \bar{\varepsilon}\\
& = \sum_{k=1}^2(\frac{\partial \theta}{\partial \theta_k}
\frac{\partial (h_{(\mu,a)})}{\partial l_k}- \frac{\partial
(h_{(\mu,a)})}{\partial
\theta_k}\frac{\partial \theta}{\partial l_k})\cdot \bar{\varepsilon}
= (-\frac{(\bar{\varepsilon}_1- \bar{\varepsilon}_9)}{ \bar{I}_1}+
\frac{\bar{\varepsilon}_9}{J_1}, \; -\frac{(\bar{\varepsilon}_2-
\bar{\varepsilon}_{10})}{ \bar{I}_2}+ \frac{\bar{\varepsilon}_{10}}{J_2}),
\end{align*}
\begin{align*}
& X_{h_{(\mu,a)}}(l) \cdot \bar{\varepsilon}=\{l,\;
h_{(\mu,a)}\}_{-}|_{\mathcal{O}_{(\mu,a)} \times \mathbb{R}^2\times
\mathbb{R}^{2*}}\cdot \bar{\varepsilon}
= \sum_{k=1}^2(\frac{\partial l}{\partial \theta_k} \frac{\partial
(h_{(\mu,a)})}{\partial l_k}- \frac{\partial (h_{(\mu,a)})}{\partial
\theta_k}\frac{\partial l}{\partial l_k})\cdot \bar{\varepsilon}=(0,0).
\end{align*}
On the other hand, from the expressions of the dynamical vector field $\tilde{X}$
and Hamiltonian vector field $X_H$, we have that
\begin{align*}
\tilde{X}(\Pi, \Gamma, \theta, l)^\varepsilon &
=T\pi_{\textmd{SE}(3)\times \mathbb{R}^2 }\cdot \tilde{X}\cdot\varepsilon(\Pi, \Gamma, \theta, l)\\
& =T\pi_{\textmd{SE}(3)\times \mathbb{R}^2 }\cdot (X_H
+ \textnormal{vlift}(u))\cdot\varepsilon (\Pi, \Gamma, \theta, l)\\
& =T\pi_{\textmd{SE}(3)\times \mathbb{R}^2 }\cdot X_H \cdot\varepsilon (\Pi, \Gamma, \theta, l)= X_H\cdot\varepsilon(\Pi, \Gamma, \theta, l),
\end{align*}
that is,
\begin{align*}
\tilde{X}(\Pi)^\varepsilon & = X_H(\Pi)\cdot\varepsilon \\
& = ( \frac{(\bar{I}_2-\bar{I}_3)\varepsilon_8\varepsilon_9+
\bar{I}_3\varepsilon_9\varepsilon_{16}}{\bar{I}_2\bar{I}_3}+
gh(\varepsilon_{11}\chi_3-\varepsilon_{12}\chi_2),\\
& \;\;\;\;\;\;
\frac{(\bar{I}_3-\bar{I}_1)\varepsilon_9\varepsilon_7-
\bar{I}_3\varepsilon_9\varepsilon_{15}}{\bar{I}_3\bar{I}_1}
+gh(\varepsilon_{12}\chi_1-\varepsilon_{10}\chi_3), \\
& \;\;\;\;\;\;
\frac{(\bar{I}_1-\bar{I}_2)\varepsilon_7\varepsilon_8-
\bar{I}_1\varepsilon_7\varepsilon_{16}
+ \bar{I}_2\varepsilon_8\varepsilon_{15}}{\bar{I}_1\bar{I}_2}
 + gh(\varepsilon_{10}\chi_2-\varepsilon_{11}\chi_1) ),
\end{align*}
\begin{align*}
\tilde{X}(\Gamma)^\varepsilon & = X_H(\Gamma)\cdot\varepsilon \\
&= ( \frac{\bar{I}_2\varepsilon_{11}\varepsilon_9
-\bar{I}_3\varepsilon_{12}\varepsilon_8+
\bar{I}_3\varepsilon_9\varepsilon_{16}}{\bar{I}_2\bar{I}_3}, \;\;
\frac{\bar{I}_3\varepsilon_{12}\varepsilon_7-\bar{I}_1\varepsilon_{10}\varepsilon_9-
\bar{I}_3\varepsilon_9\varepsilon_{15}}{\bar{I}_3\bar{I}_1}, \\
& \;\;\;\;\;\;
\frac{\bar{I}_1\varepsilon_{10}\varepsilon_8-\bar{I}_2\varepsilon_{11}\varepsilon_7-
\bar{I}_1\varepsilon_7\varepsilon_{16}
+ \bar{I}_2\varepsilon_8\varepsilon_{15}}{\bar{I}_1\bar{I}_2} ),
\end{align*}
\begin{align*}
\tilde{X}(\theta)^\varepsilon & = X_H(\theta)\cdot\varepsilon
= ( -\frac{(\varepsilon_7- \varepsilon_{15})}{ \bar{I}_1}
+\frac{\varepsilon_{15}}{J_1},\;\; -\frac{(\varepsilon_8- \varepsilon_{16})}{ \bar{I}_2}
+\frac{\varepsilon_{16}}{J_2} ),
\end{align*}
\begin{align*}
\tilde{X}(l)^\varepsilon & = X_H(l)\cdot\varepsilon =(0,0),
\end{align*}
then we have that
\begin{align*}
T\bar{\gamma}\cdot \tilde{X}(\Pi)^\varepsilon
&= ( \frac{(\bar{I}_2-\bar{I}_3)\bar{\gamma}_2\bar{\gamma}_3+
\bar{I}_3\bar{\gamma}_3\bar{\gamma}_{10}}{\bar{I}_2\bar{I}_3}+
gh(\bar{\gamma}_5\chi_3-\bar{\gamma}_6\chi_2),\\
& \;\;\;\;\;\; \frac{(\bar{I}_3-\bar{I}_1)\bar{\gamma}_3
\bar{\gamma}_1-
\bar{I}_3\bar{\gamma}_3\bar{\gamma}_9}{\bar{I}_3\bar{I}_1} +
gh(\bar{\gamma}_6\chi_1-\bar{\gamma}_4\chi_3),\\
& \;\;\;\;\;\;
\frac{(\bar{I}_1-\bar{I}_2)\bar{\gamma}_1\bar{\gamma}_2-\bar{I}_1\bar{\gamma}_1
\bar{\gamma}_{10}+
\bar{I}_2\bar{\gamma}_2\bar{\gamma}_9}{\bar{I}_1\bar{I}_2} +
gh(\bar{\gamma}_4\chi_2-\bar{\gamma}_5\chi_1) ),
\end{align*}
\begin{align*}
T\bar{\gamma}\cdot \tilde{X}(\Gamma)^\varepsilon
&= ( \frac{\bar{I}_2\bar{\gamma}_5\bar{\gamma}_3
-\bar{I}_3\bar{\gamma}_6\bar{\gamma}_2+
\bar{I}_3\bar{\gamma}_3\bar{\gamma}_{10}}{\bar{I}_2\bar{I}_3}, \;\;
\frac{\bar{I}_3\bar{\gamma}_6\bar{\gamma}_1
-\bar{I}_1\bar{\gamma}_4\bar{\gamma}_3-
\bar{I}_3\bar{\gamma}_3\bar{\gamma}_9}{\bar{I}_3\bar{I}_1}, \\
& \;\;\;\;\;\;
\frac{\bar{I}_1\bar{\gamma}_4\bar{\gamma}_2
-\bar{I}_2\bar{\gamma}_5\bar{\gamma}_1-
\bar{I}_1\bar{\gamma}_1\bar{\gamma}_{10}
+ \bar{I}_2\bar{\gamma}_2\bar{\gamma}_9}{\bar{I}_1\bar{I}_2} ),
\end{align*}
\begin{align*}
T\bar{\gamma}\cdot \tilde{X}(\theta)^\varepsilon
& = (-\frac{(\bar{\gamma}_1- \bar{\gamma}_9)}{ \bar{I}_1}+
\frac{\bar{\gamma}_9}{J_1}, \; -\frac{(\bar{\gamma}_2-
\bar{\gamma}_{10})}{ \bar{I}_2}+ \frac{\bar{\gamma}_{10}}{J_2}),
\end{align*}
\begin{align*}
T\bar{\gamma}\cdot \tilde{X}(l)^\varepsilon=(0,0).
\end{align*}
Note that $$T\bar{\lambda}\cdot \tilde{X} \cdot \varepsilon
=T\pi_{(\mu,a)}\cdot T\lambda \cdot (X_H+ \textnormal{vlift}(u))\cdot\varepsilon
=T\pi_{(\mu,a)}\cdot T\gamma \cdot T\pi_{\textmd{SE}(3)\times \mathbb{R}^2}\cdot (X_H+ \textnormal{vlift}(u))\cdot\varepsilon
=T\bar{\lambda}\cdot X_H \cdot \varepsilon,$$ that is,
\begin{align*}
T\bar{\lambda}\cdot \tilde{X}(\Pi) \cdot \varepsilon &
=T\bar{\lambda}\cdot X_H(\Pi) \cdot \varepsilon\\
&= ( \frac{(\bar{I}_2-\bar{I}_3)\bar{\lambda}_2\bar{\lambda}_3+
\bar{I}_3\bar{\lambda}_3\bar{\lambda}_{10}}{\bar{I}_2\bar{I}_3}+
gh(\bar{\lambda}_5\chi_3-\bar{\lambda}_6\chi_2),\\
& \;\;\;\;\;\; \frac{(\bar{I}_3-\bar{I}_1)\bar{\lambda}_3
\bar{\lambda}_1-
\bar{I}_3\bar{\lambda}_3\bar{\lambda}_9}{\bar{I}_3\bar{I}_1} +
gh(\bar{\lambda}_6\chi_1-\bar{\lambda}_4\chi_3),\\
& \;\;\;\;\;\;
\frac{(\bar{I}_1-\bar{I}_2)\bar{\lambda}_1\bar{\lambda}_2-\bar{I}_1\bar{\lambda}_1
\bar{\lambda}_{10}+
\bar{I}_2\bar{\lambda}_2\bar{\lambda}_9}{\bar{I}_1\bar{I}_2} +
gh(\bar{\lambda}_4\chi_2-\bar{\lambda}_5\chi_1) ),
\end{align*}
\begin{align*}
T\bar{\lambda}\cdot \tilde{X}(\Gamma) \cdot \varepsilon &
=T\bar{\lambda}\cdot X_H(\Gamma) \cdot \varepsilon\\
&= ( \frac{\bar{I}_2\bar{\lambda}_5\bar{\lambda}_3
-\bar{I}_3\bar{\lambda}_6\bar{\lambda}_2+
\bar{I}_3\bar{\lambda}_3\bar{\lambda}_{10}}{\bar{I}_2\bar{I}_3}, \;\;
\frac{\bar{I}_3\bar{\lambda}_6\bar{\lambda}_1
-\bar{I}_1\bar{\lambda}_4\bar{\lambda}_3-
\bar{I}_3\bar{\lambda}_3\bar{\lambda}_9}{\bar{I}_3\bar{I}_1}, \\
& \;\;\;\;\;\;
\frac{\bar{I}_1\bar{\lambda}_4\bar{\lambda}_2
-\bar{I}_2\bar{\lambda}_5\bar{\lambda}_1-
\bar{I}_1\bar{\lambda}_1\bar{\lambda}_{10}
+ \bar{I}_2\bar{\lambda}_2\bar{\lambda}_9}{\bar{I}_1\bar{I}_2} ),
\end{align*}
\begin{align*}
T\bar{\lambda}\cdot \tilde{X}(\theta) \cdot \varepsilon
=T\bar{\lambda}\cdot X_H(\theta) \cdot \varepsilon
= (-\frac{(\bar{\lambda}_1- \bar{\lambda}_9)}{ \bar{I}_1}+
\frac{\bar{\lambda}_9}{J_1}, \; -\frac{(\bar{\lambda}_2-
\bar{\lambda}_{10})}{ \bar{I}_2}+ \frac{\bar{\lambda}_{10}}{J_2}),
\end{align*}
\begin{align*}
T\bar{\lambda}\cdot \tilde{X}(l) \cdot \varepsilon
=T\bar{\lambda}\cdot X_H(l) \cdot \varepsilon=(0,0).
\end{align*}
Thus, when we choose that
$(\Pi,\Gamma,\theta,l)\in
\mathcal{O}_{(\mu,a)}\times \mathbb{R}^2 \times \mathbb{R}^{2*}, $ and
 $(\varepsilon_7,\varepsilon_8,\varepsilon_9)=\Pi
 =(\Pi_1,\Pi_2,\Pi_3)=(\bar{\gamma}_1,\bar{\gamma}_2,\bar{\gamma}_3)=
(\bar{\varepsilon}_1,\bar{\varepsilon}_2,\bar{\varepsilon}_3)=
(\bar{\lambda}_1,\bar{\lambda}_2,\bar{\lambda}_3), $ and
$(\varepsilon_{10},\varepsilon_{11},\varepsilon_{12})=\Gamma
=(\Gamma_1,\Gamma_2,\Gamma_3)
=(\bar{\gamma}_4,\bar{\gamma}_5,\bar{\gamma}_6)=
(\bar{\varepsilon}_4,\bar{\varepsilon}_5,\bar{\varepsilon}_6)=
(\bar{\lambda}_4,\bar{\lambda}_5,\bar{\lambda}_6), $
and $(\varepsilon_{13},\varepsilon_{14})=\theta
= (\theta_1,\theta_2)=(\bar{\gamma}_7,\bar{\gamma}_8)
=(\bar{\varepsilon}_7,\bar{\varepsilon}_8)=(\bar{\lambda}_7,\bar{\lambda}_8),
\; (\varepsilon_{15},\varepsilon_{16})=l
=(l_1,l_2)=(\bar{\gamma}_9,\bar{\gamma}_{10})
=(\bar{\varepsilon}_9,\bar{\varepsilon}_{10})
=(\bar{\lambda}_9,\bar{\lambda}_{10}). $ we must have that
\begin{align*}
& T\bar{\gamma}\cdot \tilde{X}(\Pi)^\varepsilon
=X_{h_{(\mu,a)}}(\Pi) \cdot \bar{\varepsilon}
=T\bar{\lambda}\cdot \tilde{X}(\Pi) \cdot \varepsilon, \\
& T\bar{\gamma}\cdot \tilde{X}(\Gamma)^\varepsilon
=X_{h_{(\mu,a)}}(\Gamma) \cdot \bar{\varepsilon}
=T\bar{\lambda}\cdot \tilde{X}(\Gamma) \cdot \varepsilon,\\
& T\bar{\gamma}\cdot \tilde{X}(\theta)^\varepsilon
=X_{h_{(\mu,a)}}(\theta) \cdot \bar{\varepsilon}
=T\bar{\lambda}\cdot \tilde{X}(\theta) \cdot \varepsilon, \\
& T\bar{\gamma}\cdot \tilde{X}(l)^\varepsilon
=X_{h_{(\mu,a)}}(l) \cdot \bar{\varepsilon}
=T\bar{\lambda}\cdot \tilde{X}(l) \cdot \varepsilon.
\end{align*}
Since the map $\varepsilon: T^* (\textmd{SE}(3)\times \mathbb{R}^2)
\rightarrow T^* (\textmd{SE}(3)\times \mathbb{R}^2)$ is symplectic, then
$T\bar{\varepsilon}\cdot X_{h_{(\mu,a)} \cdot \bar{\varepsilon}}
=X_{h_{(\mu,a)}} \cdot \bar{\varepsilon}. $
Thus, in this case, we must have that
$\varepsilon$ and $\bar{\varepsilon} $ are the solution of the Type II of
Hamilton-Jacobi equation
$T\bar{\gamma}\cdot \tilde{X}^\varepsilon= X_{h_{(\mu,a)}}\cdot \bar{\varepsilon}, $
for the $R_p$-reduced controlled heavy top-rotor system
$(\mathcal{O}_{(\mu,a)}\times \mathbb{R}^2\times \mathbb{R}^{2*},
\omega^{-}_{\mathcal{O}_{(\mu,a)}\times \mathbb{R}^2\times \mathbb{R}^{2*}},
h_{(\mu,a)}, u_{(\mu,a)})$, if and only if they satisfy
the equation $T\bar{\varepsilon}\cdot(X_{h_{(\mu,a)} \cdot \bar{\varepsilon}})
= T\bar{\lambda}\cdot \tilde{X}\cdot\varepsilon. $\\

To sum up the above discussion, we have the following proposition.
\begin{prop}
If the 5-tuple $(T^\ast Q,\textmd{SE}(3),\omega_Q,H,u), $ where $Q=
\textmd{SE}(3)\times \mathbb{R}^2, $ is a regular point reducible
heavy top-rotor system with the control torque $u$ acting on the rotors,
then for a point $(\mu,a) \in \mathfrak{se}^\ast(3)$, the regular
value of the momentum map $\mathbf{J}_Q: \textmd{SE}(3)\times
\mathfrak{se}^\ast(3) \times \mathbb{R}^2 \times \mathbb{R}^{2*} \to
\mathfrak{se}^\ast(3)$, the $R_p$-reduced controlled heavy top-rotor system is the 4-tuple
$(\mathcal{O}_{(\mu,a)} \times \mathbb{R}^2 \times
\mathbb{R}^{2*},\tilde{\omega}_{\mathcal{O}_{(\mu,a)} \times \mathbb{R}^2
\times \mathbb{R}^{2*}}^{-},h_{(\mu,a)},u_{(\mu,a)}), $ where $\mathcal{O}_{(\mu,a)}
\subset \mathfrak{se}^\ast(3)$ is the co-adjoint orbit,
$\tilde{\omega}_{\mathcal{O}_{(\mu,a)} \times \mathbb{R}^2 \times
\mathbb{R}^{2*}}^{-}$ is the induced symplectic form on $\mathcal{O}_{(\mu,a)}
\times \mathbb{R}^2\times \mathbb{R}^{2*} $,
$h_{(\mu,a)}(\Pi,\Gamma,\theta,l)\cdot \pi_{(\mu,a)}= \\
H(A,v,\Pi,\Gamma,\theta,l)|_{\mathcal{O}_{(\mu,a)}
\times\mathbb{R}^2\times\mathbb{R}^{2*}}$,
$u_{(\mu,a)}(\Pi,\Gamma,\theta,l)\cdot \pi_{(\mu,a)}
= u(A,v,\Pi,\Gamma,\theta,l)|_{\mathcal{O}_{(\mu,a)}
\times\mathbb{R}^2\times\mathbb{R}^{2*}}$. Assume that
$\gamma: \textmd{SE}(3)\times \mathbb{R}^2 \rightarrow
T^*(\textmd{SE}(3)\times \mathbb{R}^2)$ is an one-form on
$\textmd{SE}(3)\times \mathbb{R}^2$,
and $\lambda=\gamma \cdot \pi_{(\textmd{SE}(3)\times \mathbb{R}^2)}:
T^* (\textmd{SE}(3)\times \mathbb{R}^2 )\rightarrow
T^* (\textmd{SE}(3)\times \mathbb{R}^2), $ and $\varepsilon:
T^* (\textmd{SE}(3) \times \mathbb{R}^2 )\rightarrow
T^* (\textmd{SE}(3)\times \mathbb{R}^2) $ is a
$\textmd{SE}(3)_{(\mu,a)}$-invariant symplectic map.
Denote
$\tilde{X}^\gamma = T\pi_{(\textmd{SE}(3)\times
\mathbb{R}^2)}\cdot \tilde{X}\cdot \gamma$, and
$\tilde{X}^\varepsilon = T\pi_{(\textmd{SE}(3)\times
\mathbb{R}^2)}\cdot \tilde{X}\cdot \varepsilon$,
where $\tilde{X}=X_{(T^\ast Q,\textmd{SE}(3),\omega_Q,H,u)}$
is the dynamical vector field of
the controlled heavy top-rotor system $(T^\ast Q,\textmd{SE}(3),\omega_Q,H,u)$.
Moreover, assume that $\textmd{Im}(\gamma)\subset \mathbf{J}_Q^{-1}(\mu,a), $
and it is $\textmd{SE}(3)_{(\mu,a)}$-invariant,
and $\varepsilon(\mathbf{J}_Q^{-1}(\mu,a))\subset \mathbf{J}_Q^{-1}(\mu,a). $
Denote $\bar{\gamma}=\pi_{(\mu,a)}(\gamma):
\textmd{SE}(3)\times \mathbb{R}^2 \rightarrow
\mathcal{O}_{(\mu,a)}\times \mathbb{R}^2\times \mathbb{R}^{2*}, $ and
$\bar{\lambda}=\pi_{(\mu,a)}(\lambda): T^* (\textmd{SE}(3)\times \mathbb{R}^2) \rightarrow
\mathcal{O}_{(\mu,a)}\times \mathbb{R}^2\times \mathbb{R}^{2*}, $ and
$\bar{\varepsilon}=\pi_{(\mu,a)}(\varepsilon): \mathbf{J}_Q^{-1}(\mu,a)\rightarrow
\mathcal{O}_{(\mu,a)}\times \mathbb{R}^2\times \mathbb{R}^{2*}. $
Then the following two assertions hold:\\
\noindent $(\mathbf{i})$
If the one-form $\gamma: \textmd{SE}(3)\times \mathbb{R}^2 \rightarrow
T^*(\textmd{SE}(3)\times \mathbb{R}^2) $ is closed with respect to
$T\pi_{(\textmd{SE}(3)\times \mathbb{R}^2)}:
TT^* (\textmd{SE}(3)\times \mathbb{R}^2) \rightarrow
T(\textmd{SE}(3)\times \mathbb{R}^2), $
then $\bar{\gamma}$ is a solution of the Type I of Hamilton-Jacobi equation
$T\bar{\gamma}\cdot \tilde{X}^\gamma= X_{h_{(\mu,a)}}\cdot \bar{\gamma}; $\\
\noindent $(\mathbf{ii})$
The $\varepsilon$ and $\bar{\varepsilon} $ satisfy the Type II of Hamilton-Jacobi equation
$T\bar{\gamma}\cdot \tilde{X}^\varepsilon= X_{h_{(\mu,a)}}\cdot \bar{\varepsilon}, $
if and only if they satisfy
the equation $T\bar{\varepsilon}\cdot(X_{h_{(\mu,a)} \cdot \bar{\varepsilon}})
= T\bar{\lambda}\cdot \tilde{X}\cdot\varepsilon. $ \hskip 0.3cm $\blacksquare$
\end{prop}

When the heavy top does not carry any internal rotor, in this case
$Q=G= \textmd{SE}(3), $ and the heavy top is a regular point
reducible Hamiltonian system $(T^\ast
\textmd{SE}(3),\textmd{SE}(3),\omega, H)$, and hence it is also a
regular point reducible RCH system without the external force and
control. For a point $(\mu,a)\in \mathfrak{se}^\ast(3)$, the regular
value of the momentum map $\mathbf{J}: T^\ast \textmd{SE}(3)\to
\mathfrak{se}^\ast(3)$, the Marsden-Weinstein reduced heavy top system is
3-tuple $(\mathcal{O}_{(\mu,a)},
\omega_{\mathcal{O}_{(\mu,a)}},h_{\mathcal{O}_{(\mu,a)}})$, where
$\mathcal{O}_{(\mu,a)} \subset \mathfrak{se}^\ast(3)$ is the
co-adjoint orbit, $\omega_{\mathcal{O}_{(\mu,a)}}$ is orbit
symplectic form on $\mathcal{O}_{(\mu,a)}$, which is induced by the
heavy top Lie-Poisson bracket on $\mathfrak{se}^\ast(3)$,
$h_{\mathcal{O}_{(\mu,a)}}(\Pi,\Gamma)\cdot \pi_{\mathcal{O}_{(\mu,a)}}
=H(A,v,\Pi,\Gamma)|_{\mathcal{O}_{(\mu,a)}}$. From the above Proposition 6.3
we can obtain the Proposition 5.5 in Wang \cite{wa17}, that is, we give the two
types of Lie-Poisson Hamilton-Jacobi equation for the Marsden-Weinstein
reduced heavy top system $(\mathcal{O}_{(\mu,a)},
\omega_{\mathcal{O}_{(\mu,a)}},h_{\mathcal{O}_{(\mu,a)}})$.
See Marsden and Ratiu \cite{mara99},
Ge and Marsden \cite{gema88}, and Wang \cite{wa17}.\\

From Marsden et al \cite{mawazh10} we know that
as two regular point reducible RCH systems, the controlled rigid body with internal
rotors and the controlled heavy top with internal rotors are RpCH-equivalent.
If the equivalent map $\varphi: Q_1=\textmd{SO}(3)\times \mathbb{R}^3 \rightarrow
M (\subset Q_2= \textmd{SE}(3)\times \mathbb{R}^2)$
is a diffeomorphism,
where $M$ is a submanifold of $\textmd{SE}(3))\times \mathbb{R}^2$,
and the cotangent lift map $\varphi^*: T^*M \rightarrow T^* Q_1$ is symplectic.
Then from Theorem 4.7 there is an
induced RCH-equivalent map $\varphi^*_{\mu/G}: \mathcal{O}_{(\mu,a)}
\times \mathbb{R}^2 \times \mathbb{R}^{2*}
\rightarrow \mathcal{O}_{\mu_1} \times \mathbb{R}^3\times \mathbb{R}^{3*}, $
where $ \mathcal{O}_{\mu_1}\subset \mathfrak{so}^\ast(3)$ and
$\mathcal{O}_{(\mu,a)} \subset \mathfrak{se}^\ast(3)$ are coadjoint
orbits, respectively. Moreover, from RCH-equivalence and Theorem 4.8 we know that,
if the one-form $\gamma_2: M \rightarrow T^* M$ is closed with
respect to $T\pi_{M}: TT^* M \rightarrow TM, $ and
$\bar{\gamma}_2=\pi_{(\mu,a)}(\gamma_2): M \rightarrow (T^* M)_{(\mu,a)} $
is a solution of the Type I of Hamilton-Jacobi equation for
the $R_p$-reduced controlled heavy top-rotor system
$(\mathcal{O}_{(\mu,a)} \times \mathbb{R}^2 \times
\mathbb{R}^{2*},\tilde{\omega}_{\mathcal{O}_{(\mu,a)} \times \mathbb{R}^2
\times \mathbb{R}^{2*}}^{-},h_{(\mu,a)},u_{(\mu,a)}), $
then $\gamma_1=\varphi^* \cdot \gamma_2\cdot \varphi: Q_1 \rightarrow T^* Q_1, $ is
a solution of the Type I of Hamilton-Jacobi equation for the controlled
rigid body-rotor system $(T^\ast (\textmd{SO}(3)\times \mathbb{R}^3),\textmd{SO}(3),
\omega_{\textmd{SO}(3)\times \mathbb{R}^3},H,u) $,
and $\bar{\gamma}_1=\pi_{\mu_1}(\gamma_1): Q_1 \rightarrow (T^* Q_1)_{\mu_1}$ is
a solution of the Type I of Hamilton-Jacobi equation for the $R_p$-reduced controlled
rigid body-rotor system $(\mathcal{O}_{\mu_1} \times \mathbb{R}^3 \times
\mathbb{R}^{3*},\tilde{\omega}_{\mathcal{O}_{\mu_1} \times \mathbb{R}^3
\times \mathbb{R}^{3*}}^{-},h_{\mu_1},u_{\mu_1}) $.
Moreover, if the $\textmd{SE}(3)_{(\mu,a)}$-invariant symplectic map
$\varepsilon_2: T^*M \rightarrow T^* M$ and
$\bar{\varepsilon}_2=\pi_{(\mu,a)}(\varepsilon_2):
\mathbf{J_2}^{-1}(\mu,a) (\subset T^* M) \rightarrow (T^* M)_{(\mu,a)} $
satisfy the Type II of Hamilton-Jacobi equation for
the $R_p$-reduced controlled heavy top-rotor system
$(\mathcal{O}_{(\mu,a)} \times \mathbb{R}^2 \times
\mathbb{R}^{2*},\tilde{\omega}_{\mathcal{O}_{(\mu,a)} \times \mathbb{R}^2
\times \mathbb{R}^{2*}}^{-},h_{(\mu,a)},u_{(\mu,a)})$, then $\varepsilon_1=
\varphi^* \cdot \varepsilon_2\cdot \varphi_*: T^*Q_1 \rightarrow T^* Q_1 $ and
$\bar{\varepsilon}_1=\pi_{\mu_1}(\varepsilon_1):
\mathbf{J_1}^{-1}(\mu_1) (\subset T^*Q_1) \rightarrow (T^* Q_1)_{\mu_1} $
satisfy the Type II of Hamilton-Jacobi equation for the $R_p$-reduced controlled
rigid body-rotor system $(\mathcal{O}_{\mu_1} \times \mathbb{R}^3 \times
\mathbb{R}^{3*},\tilde{\omega}_{\mathcal{O}_{\mu_1} \times \mathbb{R}^3
\times \mathbb{R}^{3*}}^{-},h_{\mu_1},u_{\mu_1}).$\\

The theory of controlled mechanical system is a very important subject,
its research gathers together some separate areas of research such as
mechanics, differential geometry and nonlinear control theory, etc.,
and the emphasis of this research on geometry is motivated by the
aim of understanding the structure of equations of motion of the
system in a way that helps both analysis and design. Thus, it is
natural to study the controlled mechanical systems by combining with the
analysis of dynamical systems and the geometric reduction theory of
Hamiltonian and Lagrangian systems. Following
the theoretical development of geometric mechanics, a lot
of important problems about this subject are being explored and
studied, see Wang \cite{wa18}.
In this paper, we give precisely the geometric constraint conditions of
canonical symplectic form and its regular reduced symplectic forms for the
dynamical vector fields of the RCH system and its regular reducible RCH systems.
we also give the geometric constraint conditions of
the magnetic symplectic form for the controlled magnetic Hamiltonian vector field.
These conditions are called the two types of Hamilton-Jacobi equations, which are the
development of classical Hamilton-Jacobi equation given by Abraham and Marsden
\cite{abma78}, also see Wang \cite{wa17}. In particular,
it is worthy of noting that, in the applications of the Hamilton-Jacobi
theoretical result for the regular point reducible RCH system,
the motions of the controlled rigid body-rotor system
and the controlled heavy top-rotor system are different, and
the configuration spaces, the Hamiltonian functions, the actions of Lie groups,
the $R_p$-reduced symplectic forms and the $R_p$-reduced systems of
the controlled rigid body-rotor system
and the controlled heavy top-rotor system are also different. But,
the two types of Hamilton-Jacobi equations given by calculation in detail
are same, that is, the internal rules are same by comparing
Proposition 6.2 and Proposition 6.3.
Moreover, Wang in \cite{wa20a, wa13e} applies the research work
to give explicitly the two types of Hamilton-Jacobi
equations for the regular reduced controlled spacecraft-rotor system and
the regular reduced controlled underwater vehicle-rotor system,
and in de Le\'{o}n and Wang \cite{lewa15}, the authors
give also the two types of Hamilton-Jacobi equations for the nonholonomic
Hamiltonian system and the nonholonomic reducible Hamiltonian system
on a cotangent bundle, by using the distributional Hamiltonian system
and the reduced distributional Hamiltonian system. These researches
reveal from the geometrical point of view the internal relationships of
constraints, symplectic forms and dynamical vector fields of a
mechanical system and its regular reduced systems.
It is the key thought of the researches of geometrical mechanics
of the professor Jerrold E. Marsden to explore and reveal the deeply internal
relationship between the geometrical structure of phase space and the dynamical
vector field of a mechanical system. It is also our goal of pursuing and inheriting.

\end{document}